\documentclass[a4paper]{article}

\usepackage{amsfonts,amssymb,amsthm}
\usepackage[totalwidth=17cm,totalheight=24cm]{geometry}

\newcommand{\C}{{\mathbb C}}
\newcommand{\N}{{\mathbb N}}
\newcommand{\R}{{\mathbb R}}
\newcommand{\Z}{{\mathbb Z}}
\newcommand{\cL}{{\mathcal L}}
\newcommand{\cH}{{\mathcal H}}
\newcommand{\cM}{{\mathcal M}}
\newcommand{\cT}{{\mathcal T}}

\newcommand{\St}{\tilde{S}}
\newcommand{\be}{\begin{equation}}
\newcommand{\ee}{\end{equation}}

\newtheorem{prop}{Proposition}[section]
\newtheorem{df}[prop]{Definition}
\newtheorem{lemma}[prop]{Lemma}
\newtheorem{thm}[prop]{Theorem}
\newtheorem{cor}[prop]{Corollary}

\newtheorem{remark}[prop]{Remark}
\newtheorem{example}[prop]{Example}
\newtheorem{question}[prop]{Question}

\newcommand{\II}{I\hspace{-0.1cm}I}
\newcommand{\III}{I\hspace{-0.1cm}I\hspace{-0.1cm}I}
\newcommand{\dr}{\partial}
\newcommand{\hess}{\mbox{Hess}}
\newcommand{\tr}{\mbox{tr}}
\newcommand{\mess}{\mbox{Mess}}
\newcommand{\ext}{\mbox{Ext}}
\newcommand{\isom}{\mbox{Isom}}

\begin{document}

\title{Minimal surfaces and particles in 3-manifolds}

\date{January 2006 (v2)}

\author{Kirill Krasnov
\thanks{School of Mathematical Sciences, University of Nottingham, Nottingham,
  NG7 2RD, UK.  
\texttt{krasnov@maths.nott.ac.uk}}
and 
Jean-Marc Schlenker
\thanks{
Laboratoire Emile Picard, UMR CNRS 5580,
Institut de Math{\'e}matiques, 
Universit{\'e} Paul Sabatier,
31062 Toulouse Cedex 9,
France.
\texttt{http://www.picard.ups-tlse.fr/\~{ }schlenker}.}
}

\maketitle

\begin{abstract}

We use minimal (or CMC) surfaces to describe 3-dimensional
hyperbolic, anti-de Sitter, de Sitter or Minkowski manifolds. 
We consider whether these manifolds admit ``nice'' foliations and explicit
metrics, and whether the space of these metrics has a simple description
in terms of Teichm\"uller theory. In the hyperbolic settings both questions
have positive answers for a certain subset of the quasi-Fuchsian
manifolds: those containing a closed surface with principal curvatures at
most $1$. We show that this subset is parameterized by an open domain of the 
cotangent bundle of Teichm\"uller space. These results are extended 
to ``quasi-Fuchsian'' manifolds with conical singularities 
along infinite lines, known in the physics literature as ``massive, 
spin-less particles''. 

Things work better for globally hyperbolic
anti-de Sitter manifolds: the parameterization by the cotangent
of Teichm\"uller space works for all manifolds. 
There is another description of this moduli space as the product
two copies of Teichm\"uller space due to Mess. Using the 
maximal surface description, we propose a new parameterization by two 
copies of Teichm\"uller space, alternative to that of Mess, and extend 
all the results to manifolds with conical singularities along time-like 
lines. Similar results are obtained for de Sitter or Minkowski manifolds. 

Finally, for all four settings, we show that the symplectic 
form on the moduli space of 3-manifolds that comes from 
parameterization by the cotangent bundle of Teichm\"uller space is the same 
as the 3-dimensional gravity one. 

\bigskip

\begin{center} {\bf R{\'e}sum{\'e}} \end{center}

On utilise des surfaces minimales (ou CMC) pour d\'ecrire des 
vari\'et\'es de dimension 3 hyperboliques, anti-de Sitter, de Sitter,
ou Minkowski. Les principales questions consid\'er\'ees sont si ces
vari\'et\'es admettent des ``bon'' feuilletages et des m\'etriques 
explicites, et si les espaces de modules de ces m\'etriques on des descriptions
simples en termes de th\'eorie de Teichm\"uller. Dans le cadre hyperbolique,
ces deux questions ont des r\'eposes positives si on se restreint \`a un
sous-ensemble des m\'etriques quasi-fuchsiennes~: celles qui contiennent
une surface ferm\'ee dont les courbure principales sont partout dans 
$(-1,1)$. L'espace de ces vari\'et\'es est param\'etr\'e par un ouvert 
dans le cotangent \`a l'espace de Teichm\"uller. Ces r\'esultats s'\'etendent
\'a des vari\'et\'es ``quasi-fuchsiennes'' ayant des singularit\'es
coniques le long de lignes infinies, connues dans la litt\'erature
physique sous le nom de ``particules massives sans spin''.

Les choses se passent mieux pour les vari\'et\'es anti-de Sitter globalement 
hyperboliques~: la param\'etrisation par le cotangent de l'espace de
Teichm\"uller fonctionne pour toutes les vari\'et\'es. Il y a une autre
description de cet espace de module, d\^u \`a Mess, comme le produit de
deux copies de l'espace de Teichm\"uller. En utilisant la description 
par les surfaces maximales, on propose une autre param\'etrisation par
le produit de deux copies de l'espace de Teichm\"uller, et on \'etend 
ces r\'esultats \`a des vari\'et\'es ayant des singularit\'es coniques
le long de courbes de type temps. Des r\'esultats similaires sont obtenus
pour les vari\'et\'es de Sitter ou Minkowski.

Finalement, pour chacun des quatre cadres, on montre que la forme
symplectique sur l'espace des modules de 3-vari\'et\'es provenant de 
la param\'etrisation par le cotangent \`a l'espace de Teichm\"uller
est identique \`a celle qui provient de la gravit\'e en dimension 3.

\end{abstract}

\maketitle

\vspace{0.4cm}


\vspace{0.4cm}


\bigskip

\tableofcontents

\section{Introductions, results}

The results presented in this paper could be of interest from two different
viewpoints, 
the geometry of constant curvature 3-manifolds and 3d gravity. So we provide
two introductions, hoping that most potentially interested readers can find
adequate motivations in one or in the other. Readers more interested in
physics aspects 
are advised to skip the maths-oriented introduction and start by 
reading the ``physics'' introduction.

\subsection{A mathematical introduction}

\paragraph{Canonical foliations of quasi-Fuchsian manifolds.}

Quasi-Fuchsian hyperbolic 3-manifolds are topologically simple objects which
can exhibit a rich and interesting geometry \cite{thurston-notes}. One of the
main motivation for this paper -- which remains mostly unfulfilled -- is to
give explicit forms of the hyperbolic metrics on such manifolds, and to
understand to what extend they have canonical foliations by surfaces. This
leads us to consider ``simple'' quasi-Fuchsian 3-manifolds, for which such
a foliation exist. 

\begin{df}
A quasi-Fuchsian hyperbolic 3-manifold is {\bf almost-Fuchsian} if it contains
a closed surface with principal curvatures everywhere in $(-1,1)$.
\end{df}

Simple arguments going back to Uhlenbeck \cite{uhlenbeck}, 
given in section 2, show that almost-Fuchsian manifolds
contain a unique closed 
minimal surface, and that they have a canonical foliation by
a family of equidistant surfaces which includes the minimal one. 

Moreover those minimal surfaces can be described in complex analytic terms
--- using ideas going back to H. Hopf \cite{hopf}, 
and used recently by Taubes \cite{taubes} ---
which leads to an identification of the space of almost-Fuchsian manifolds
with an open subset of the cotangent space of Teichm\"uller space. This in
turns leads to a simple and explicit form of the hyperbolic metric on
almost-Fuchsian manifolds, first found by V. Fock \cite{Fock}.

Almost-Fuchsian manifolds therefore appear as particularly nice and
well-behaved among quasi-Fuchsian manifolds. Things would be quite simple if
all quasi-Fuchsian manifolds were actually almost-Fuchsian, however we show in
subsection 2.3, using a recent argument of Haas and Thurston on hyperbolic
manifolds which fiber over the circle, that this is not the case.

\paragraph{Other settings.}

The description of 3-manifolds as foliated by surfaces equidistant to a
surface of vanishing mean curvature works even better 
in contexts other than that
of quasi-Fuchsian manifolds. Because of this, the results stated in the hyperbolic
setting should be considered as prototypes --- easily understandable for
geometers interested in hyperbolic geometry --- of ``better'' results in the
Lorentzian contexts. Thus, apart from the 
case of hyperbolic 3-manifolds mentioned above, we consider the
cases of Lorentzian manifolds with curvature of both signs (or zero). In
other words, we also consider anti-de Sitter (or AdS) 3-manifolds (that is
Lorentzian signature, negative curvature), de Sitter (or dS)
3-manifolds (Lorentzian signature, positive curvature), as well as
flat Lorentzian signature 3-manifolds (the Minkowski case). 

\paragraph{GHMC AdS manifolds.}

This is the best case scenario for our setup. In the AdS case, 
the description of a 3-manifold by a foliation equidistant to a surface
of vanishing mean curvature works always, and there is no need for
the notion of an almost-Fuchsian manifold of the hyperbolic setting. 
More precisely, we shall work in the context of globally hyperbolic maximally compact
(or GHMC) anti-de Sitter (or AdS) 3-manifolds. Those were introduced
by G. Mess \cite{mess}, who showed that they share several key properties of
quasi-Fuchsian hyperbolic manifolds. 

It was recently shown by Barbot, B\'eguin and Zeghib
\cite{BBZ-cras,barbot-zeghib} that those manifolds
contain a unique closed space-like maximal surface. Moreover, the ideas of
Taubes on minimal surfaces in hyperbolic 3-manifolds, transfered to the
AdS context, show directly that those
maximal surfaces are parametrized by the cotangent bundle of Teichm\"uller space, so that
$T^*\cT_g$ parameterizes the space of GHMC AdS manifolds of genus $g$. Our Theorem \ref{tm:ads}
is a statement to this effect.

Mess \cite{mess} had earlier shown that the GHMC AdS manifolds of genus $g$
are parametrized by two copies of $\cT_g$, which yields a natural map from
$T^*\cT_g$ to $\cT_g\times \cT_g$ which might be interesting in its own
right. In section 3, we recover this Mess parameterization using (smooth)
space-like surfaces, and give another, similar but distinct parameterization
based on maximal surfaces.

\paragraph{Minkowski and de Sitter 3-manifolds.}

Our description of a 3-manifold as foliated by surfaces equidistant to the
minimal (resp. maximal) one works well also in the context of Minkowski
3-dimensional 
manifolds, with a small twist: it is necessary to consider constant mean 
curvature, rather than maximal, surfaces. 
Similarly, one should consider CMC surfaces in dS manifolds. Both cases
lead
to parameterizations of the moduli space by the cotangent bundle of
the Teichm\"uller space.

\paragraph{Hyperbolic cone-manifolds.}

Another aim of the paper is to extend the results on foliations and on
the ``Hamiltonian'' description using the cotangent bundle of
Teichm\"uller space to the case of 3-dimensional manifolds containing
conical singularities along lines. 
Physically, conical singularities along lines 
in 3 dimensions describe matter sources,
so their inclusion is necessary to go from a physically void theory with only
global topological degrees of freedom to a theory with rich local dynamics.
  
Recently, Bromberg \cite{bromberg2}
has considered another kind of singular convex co-compact hyperbolic
manifolds, with singular locus a disjoint union of closed curves. It is
perhaps necessary to emphasise the difference with the objects considered
here, in which the singular locus is a union of open curves ending on the
boundary at infinity. 

The singular quasi-Fuchsian manifolds considered here, when the cone angles
are between $0$ and $\pi$, have
interesting geometric properties generalizing those found in the non-singular
case. Thus, it
is  for instance possible to define their convex core, and its boundary has an
induced metric which is a hyperbolic metric with conical singularities of
angles equal to the singular angle of the manifold, and with a bending
lamination which ``avoids'' the singularities. 
This is an important motivation for us, however we do not
elaborate on the geometry of the convex core here and refer the reader to a future publication.

The constructions outlined above, concerning (smooth) quasi-Fuchsian
manifolds, extend to singular quasi-Fuchsian manifolds, as defined here. There
is a notion of singular almost-Fuchsian manifolds, which is defined as those
containing a minimal surfaces (which is ``orthogonal'' to the singular locus,
in a way defined in section 3) with principal curvatures everywhere in
$(-1,1)$. Those singular almost-Fuchsian manifolds contain no other minimal
surface, and their hyperbolic metric admits a canonical foliation by
equidistant surfaces. Their moduli space has a natural parametrization, for
fixed singular angles, by 
an open subset in $T^*\cT_{g,n}$, the cotangent bundle 
of the Teichm\"uller space 
with $n$ marked points on a surface of genus $g$.
This statement constitutes our
Theorem \ref{tm:hyper-cone-main}.

\paragraph{Teichm\"uller theory with marked points.}

Just as the smooth almost-Fuchsian manifolds can be associated to points in a
subset of the cotangent space of Teichm\"uller space, the singular
almost-Fuchsian manifolds with $n$ singular curves 
can be constructed from points in a subset of the
cotangent bundle of the Teichm\"uller space with $n$ marked points. 
This fact, together with the corresponding AdS results, is interesting for
reasons 
related to quantum gravity. Namely, singular lines can be interpreted as
particles and the fact that cotangent spaces are particularly simple to
quantize 
gives hope for an explicit quantum theory. 
In this respect it is interesting to note that, in the different contexts 
considered here, the symplectic forms induced on the moduli spaces of
3-dimensional 
metrics by the parametrization by the cotangent of Teichm\"uller space are the
same as the symplectic forms of 3-dimensional gravity; 
see the ``physics introduction'' below in this section
for more on these physics motivations behind our constructions.

The manifolds refered to, and, more generally, the singular quasi-Fuchsian
manifolds, 
have a boundary at infinity which is naturally endowed with a conformal
structure with marked points, which correspond to the endpoints of the
singular curves. It would be interesting to know whether those two conformal
structures with marked points, along with the total angles around the
singular curves, uniquely determine a singular quasi-Fuchsian manifold, as
happens in the non-singular case by the Ahlfors-Bers theorem. We do not
achieve such a result here. However, the analogous AdS case is simpler and
does yield some statements.

Quasi-fuchsian 3-manifolds can actually be considered as relevant tools to
study Teichm\"uller space, as is demonstrated for instance by the recent work
of Takhtajan and Teo \cite{takhtajan-teo} and of T. Hodge
\cite{hodge-these}. In this respect, one can hope that quasi-fuchsian
cone-manifolds will play a similar role for Teichm\"uller space with some
marked points. It also appears that GHMC AdS manifolds can provide an
alternative to quasi-fuchsian hyperbolic manifolds in this respect, either
with or without conical singularities.

\paragraph{AdS manifolds with conical singularities.}

As in the hyperbolic setting, the ideas described can be extended naturally
from GHMC AdS manifolds to a wider setting of ``GHMC AdS cone-manifolds'' with
singularities along time-like lines, this is the theme of section 5.  
Similarly to the hyperbolic case, such cone-manifolds 
are parametrized by total angles around the singular curves and by the
cotangent bundle of $\cT_{g,n}$, where $n$ is the number of singular curves. What is
different in the AdS case is that it is all of the cotangent bundle 
that arises this way, not
just a subset. Thus, the AdS situation is much simpler than the hyperbolic
case.

As in the hyperbolic case, the AdS cone-manifolds that we consider here share,
when the angles at the singular lines are in $(0,\pi)$, some key properties of
the corresponding smooth GHMC AdS manifolds. The properties of the convex core
and of its boundary, in particular, remain similar to those in the
non-singular setting. However, in order to keep the paper within
reasonable limits we decided not to consider the geometry of the
convex core here. 

Given an AdS cone-manifold $M$ with singularities, it is possible to associate to
$M$ a pair of hyperbolic metrics, and therefore a point in
$\cT_{g,n}\times \cT_{g,n}$.  This is done by a construction similar to the one given 
for the smooth case in section 3, and which is a more differential-geometric 
reformulation of a construction of Mess \cite{mess}. 
Thus, as for smooth manifolds, the parameterization of singular
GHMC AdS manifolds by $T^*\cT_{g,n}$ provides a natural map from
$T^*\cT_{g,n}$  to $\cT_{g,n}\times\cT_{g,n}$. It is not clear at this point
whether this map is one-to-one.

\paragraph{De Sitter and Minkowski cone-manifolds.}

The results on AdS cone-manifolds, in particular the description that uses the
cotangent bundle of Teichm\"uller space, mostly extend to the setting of singular
GHMC de Sitter and Minkowski manifold. This is done in section 6.

\paragraph{Infinite area surfaces and Schottky manifolds.}

Rather than considering ``spatial'' surfaces $\Sigma$ which are closed, one
can  
consider surfaces that have ends of infinite area. This case is of great
importance 
for physical applications --- indeed, in physics one is typically interested in
spacetimes 
that have non-compact spatial slices and thus have ``infinity'' where
``observers'' can live. 
It turns out that this
case can be treated by a simple extension of the formalism, namely, by 
adding the ``conical singularities'' with the ``angle deficit'' larger than
$2\pi$. However, the resulting hyperbolic manifolds are very different. 
In particular, the conformal boundary at infinity consists of a single
component. 
Under the condition that the principal curvatures are in
$(-1,1)$, those manifolds are foliated by equidistant surfaces, and contain a
unique (complete) minimal surface with a constrained behaviour at infinity. So
we obtain a natural parameterization of such manifolds by an open subset of the
cotangent bundle of a Teichm\"uller space with some marked points. When
conical singularities are allowed on the minimal surface, the resulting
hyperbolic 3-manifolds also have conical singularities along lines. 

This construction also extends to the AdS and Minkowski setting, and provide
``Schottky type'' manifolds, possibly with conical singularities along
time-like curves, which are parametrized in a simple way by Teichm\"uller
data. 
There are, however, some subtle technical differences between the ``Schottky''
case and  
the cases treated in this paper. For this reason we decided to leave the
open surface case for a separate publication.

\subsection{Physics introduction} 

\paragraph{Hamiltonian formulation of 3d gravity.}

The mathematical physics community has
become genuinely interested in the subject of 3d gravity after Witten
\cite{witten-jones} has shown that the theory is exactly soluble and, to
a certain extent, explicitly quantizable. Witten considered the case
of zero cosmological constant Lorentzian signature gravity and,
using the Chern-Simons formulation of the theory, has shown that the
phase space of the theory is the cotangent bundle over the 
Teichm\"uller space of the spatial slice. The main idea of the
argument is so simple that it is worth being repeated here. 

One starts by formulating the theory of flat 3d gravity as the
CS theory of ${\rm ISO}(1,2)$ connections. One then considers 
3-manifolds of topology $M=\Sigma\times \R$ so that
the Hamiltonian formulation of the theory is possible. The phase
space of CS theory (which is also the space of its classical solutions)
is the space of flat ${\rm ISO}(1,2)$ connections on $\Sigma$. In the
situation  
at hand this is the same as the cotangent bundle over the space ${\cal M}$ of
flat 
${\rm SO}(1,2)={\rm SL}(2,\R)$ connections. One of the components of
${\cal M}$ is the Teichm\"uller space $\cT_\Sigma$ of $\Sigma$, and it is
shown to be the physically right component to consider, so the
phase space becomes the cotangent bundle $T^* \cT_\Sigma$.

A related development has been that by Moncrief \cite{Moncrief} who showed
that the phase space of 3d gravity is $T^* \cT_\Sigma$ working entirely
in the geometrodynamics setting. He considered the Hamiltonian 
description corresponding to a foliation of $M$ by surfaces $\Sigma_t$ of 
constant mean curvature. The constant mean curvature condition
can be viewed as a gauge choice. Accordingly, one gets a true
Hamiltonian, not just the Hamiltonian constraint. The Hamiltonian
$H(t)$ turns out to be the area of the constant mean curvature slices
$\Sigma_t$, and generates a flow in $T^* \cT_\Sigma$. Moncrief had
also pointed out that similar methods could be used to study the
Hamiltonian reduction of 3d gravity with non-zero 
cosmological constant. 

This paper gives an explicit description of various classes
of 3-manifolds in terms of data on zero (more generally constant) mean curvature surfaces.
In particular, we provide a very concrete description of a metric in 
the 3-manifold by considering equidistant foliations. These results can be used to give a
Hamiltonian description of the corresponding gravity theories. 
Thus, we propose a Hamiltonian formulation alternative to that
in \cite{Moncrief}, which is based not on constant mean curvature but
on foliations equidistant to a minimal (or, more generally, constant mean curvature) 
surface. As such a description is
available and requires almost no modifications for several 
different types of 3-manifolds, we automatically get a
very similar Hamiltonian description for all these cases. Thus,
in all the cases considered in this paper, namely the hyperbolic,
AdS, dS and Minkowski cases, the phase space is shown 
to be (possibly a subset of) the
cotangent bundle over the Teichm\"uller space of the minimal (or, more generally, constant
mean curvature) surface. We believe that the Hamiltonian
description proposed here will be instrumental in the future quantizations
of 3d gravity.

One basic reason why the phase space turns out to be always related to the
cotangent bundle $T^* \cT_\Sigma$ is simple enough to be explained here in the introduction. 
In the Hamiltonian formulation, the phase space of gravity is essentially the space of
pairs $(I,\II)$, where $I,\II$ are the first and second fundamental forms of the
spatial slice. In 3d gravity the spatial slice $\Sigma$ is 2-dimensional,
and the induced metric $I$ can be written as $I=e^\varphi |dz|^2$, where
$e^\varphi$ is some conformal factor, and $|dz|^2$ is the (locally) flat metric in the 
conformal class of $I$. In other words, the first fundamental form can be
parameterized by a point $c\in \cT_\Sigma$ in the Teichm\"uller space, together
with the Liouville field $\varphi$. The Liouville field $\varphi$ is not arbitrary, and
is required to satisfy the Gauss equation. In turn, the second fundamental
form $\II$ of $\Sigma$ can be shown to give rise to a holomorphic
quadratic differential $t$ on $\Sigma$. Thus, consider: $\II' = \II - (H/2) I$,
where $H={\rm Tr}(I^{-1} \II)$ is the mean curvature. The new form $\II'$ is
thus traceless. In the complex-analytic description, this means that 
$\II'$ is the real part $\II' = tdz^2 + \bar{t} d\bar{z}^2$ of a quadratic 
differential $t$ on $\Sigma$. It is then a simple exercise to see that, 
whenever $H$ is constant on $\Sigma$, the Codazzi equation $d^{\nabla^I} \II=0$ implies 
that the quadratic differential $t$ is holomorphic. Thus, a constant mean
curvature surface always gives rise to a point $(c,t)\in T^* \cT_\Sigma$.
Under some conditions one can invert this map, and use $T^* \cT_\Sigma$ to
parametrize the phase space. Most of this paper deals precisely with the
question of when the described map is invertible. This is related to the question of
existence and uniqueness of solutions for the Gauss equation for the
Liouville field $\varphi$. Whenever the map 
(constant mean curvature surface) $\to T^* \cT_\Sigma$ is invertible then
(possibly a subspace of) the cotangent bundle can be used to describe the gravity phase
space. It is then an easy exercise to check that the gravitational
symplectic structure is always that of the cotangent bundle. The corresponding
computation is carried out in section 7. This
provides an explanation for the fact that the phase space of gravity
(in the formulation based on minimal, or, more generally, constant
mean curvature surfaces) is always (possibly a subset of) the cotangent
bundle over the Teichm\"uller space. Curiously enough, in most cases,
in particular in all cases when the signature is Lorentzian, the
phase space happens to be all of the cotangent bundle. This 
confirms one's physical intuition that the Hamiltonian description
that uses the initial data really only makes sense when the signature 
is Lorentzian. When dealing with Euclidean manifolds it is much
more natural to consider the associated Dirichlet problem, 
and encode the geometry of the manifold by data on its boundary.
The fact that we do get a phase space-like description of the
hyperbolic case is shadowed by the fact that only a small subset
of quasi-Fuchsian manifolds is described this way. This
should be compared with the Bers description, where
the data are specified on the conformal boundaries of the
space, and which is applicable to all the convex co-compact 
hyperbolic manifolds, and to the comparable parametrization discovered by Mess
\cite{mess} for GHMC AdS manifolds.

The argument given in the beginning of the previous paragraph is so 
general that it seems to imply that
the phase space is {\it always} related to $T^* \cT_\Sigma$, namely for all
combinations of signature and sign of the curvature. This is not so, there are
two cases that are exceptions. Namely, the phase space of the 
Euclidean signature gravity with zero or positive curvature is
{\it not} related to the cotangent bundle over the Teichm\"uller space
of the spatial slice. The basic reason for this is that it is not
possible to embed a surface of genus $g>1$ into $\R^3$ or $S^3$ in
such a way that the mean curvature is constant.\footnote{Note, however,
that it {\it is} possible to have an immersion of a $g>1$ as constant
mean curvature surfaces in these spaces, but this immersion fails to be an embedding.} 
This explains why these two settings are different.
We do not consider them in the present paper.

For the convenience of the reader we would like to summarize all the facts about 
3d gravity phase space in the following tables. Some of these facts are well-known,
some other follow from the results of this work. The first table lists the
groups of isometry of the corresponding maximally symmetric spaces.

\begin{table}[ht]\begin{center}
\begin{tabular}{|c|c|c|c|}
\hline 
{} & Positive Curvature & Zero Curvature & Negative Curvature \\
\hline 
Euclidean & $SO(4)\sim SU(2)\times SU(2)$ & $SO(3)\ltimes \R^3$ & $SO(1,3)\sim SL(2,\C)$ \\
\hline 
Lorentzian & $SO(1,3)\sim SL(2,\C)$ & $SO(1,2)\ltimes \R^3$ & $SO(2,2)\sim 
SL(2,\R)\times SL(2,\R)$ \\ 
\hline 
\end{tabular}
\caption{Isometry groups}\end{center}
\end{table}

The isometry groups listed are also the groups appearing in the Chern-Simons formulation
of the corresponding theories. Since the phase space of Chern-Simons theory on a
manifold of topology $\Sigma\times \R$ is (possibly a component of) the space
${\rm Hom}(\pi_1(\Sigma),G)/G$, where $G$ is the gauge group in question, we
immediately get one possible description of the phase space in each case. In some cases
this description can be further simplified. Thus, in the case of zero curvature, Euclidean signature,
using the fact that the isometry group is a semi-direct product, one arrives at
the phase space being $T^* {\cal A}_\Sigma$, where 
${\cal A}_\Sigma$ is the space of flat ${\rm SU}(2)$ connections on $\Sigma$.
In the zero curvature Lorentzian case one arrives at the phase space being
$T^* \cT_\Sigma$, where $\cT_\Sigma$, the Teichm\"uller space of $\Sigma$ is (a component of the) space of flat
${\rm SL}(2,\R)$ connections on $\Sigma$.
In the positive curvature Euclidean case the phase space becomes the product:
${\cal A}_\Sigma\times {\cal A}_\Sigma$. The other cases are as follows. In the
negative curvature Euclidean case the phase space of almost-Fuchsian 
(AF) manifolds is a subset of $T^* \cT_\Sigma$ by the results of this paper. 
A complete description of the moduli space (phase space) is provided in
this case by the Bers uniformization by two copies of Teichm\"uller space. 
This description also extends to the positive curvature Lorentzian case by the well-known
hyperbolic-dS duality. However, in this case the initial data description is
also possible, and is given by all of $T^* \cT_g$. This is a bit surprising, 
in view of the duality to the hyperbolic case. The reason why all of
the cotangent bundle appears and not just a subset of it is related to
the fact that one uses not minimal, but constant mean curvature surfaces
in the dS setting. In the 
case of negative curvature Lorentzian signature the phase space is the product
of two copies of Teichm\"uller space by the result of Mess (and, in a
different way, 
by the results of this paper), and is also the cotangent bundle of the
Teichm\"uller space 
by our results. This is summarized by the following table:
\begin{table}[h]\begin{center}
\begin{tabular}{|c|c|c|c|}
\hline 
{} & Positive Curvature & Zero Curvature & Negative Curvature \\
\hline 
Euclidean & ${\cal A}_\Sigma\times{\cal A}_\Sigma$ & $T^* {\cal A}_\Sigma$ & 
$\cT_\Sigma\times \cT_\Sigma$, AF -- subset of $T^* \cT_\Sigma$ \\
\hline 
Lorentzian & $T^* \cT_\Sigma\sim \cT_\Sigma\times \cT_\Sigma$ &
$T^* \cT_\Sigma$  
& $T^* \cT_\Sigma \sim \cT_\Sigma \times \cT_\Sigma$ \\
\hline 
\end{tabular}
\caption{Phase spaces}\end{center}
\end{table}

\paragraph{Analytic continuations.} 

Our results give an interesting perspective on the question of ``analytic
continuation'' between different constant curvature 3-manifolds. Such an
analytic continuation was considered by one of the authors in \cite{Holography}
as a continuation between the negative curvature Lorentzian and Euclidean
signature cases, for the case of Fuchsian manifolds. It was then studied
infinitesimally away from the Fuchsian case in \cite{Continuation}. More
recently, a general ``analytic continuation'' procedure was proposed by
Benedetti and Bonsante in \cite{Benedetti}. Here we would like to describe
a procedure that follows from the results of the present paper, 
and compare it to that of \cite{Benedetti}.

The fact that there is a relation between spaces of different signatures and/or different
signs of the curvature is well-known. The simplest example of such a relation
is the well-known duality between the hyperbolic and dS manifolds. The dS manifold dual
to a quasi-Fuchsian hyperbolic manifold $M$ is obtained by considering 
all the geodesic planes in the covering
space of $M$ that do not intersect the cover of the convex core. These geodesic
planes correspond to a set of points in dS space on which the action of the group
is properly discontinuous. The quotient dS manifold has a conformal boundary 
consisting of two components (as $M$), with the same conformal structures. Unlike
$M$, however, the dual dS space is not connected. It consists of two components each
with its own conformal boundary and each ending on a singular graph. The
universal cover of this graph consists of points in $dS_3$ such that the dual planes 
in $H^3$ are tangent to the universal cover of the boundary of the convex core.

Let us consider other possible examples of such relations. 
The fact that the Lorentzian cases of zero and negative curvature have the same
phase space suggests that there must be a map between such manifolds. In fact, there
is a more general relation involving the zero curvature Lorentzian case together
with dS, AdS and the hyperbolic cases. The corresponding relations have been
exhibited by Benedetti and Bonsante \cite{Benedetti}. The key tool used in their
work is that of measured geodesic laminations ${\cal ML}$ on a hyperbolic surface $\Sigma$. One
then notices that measured geodesic laminations naturally appear in all the
contexts mentioned. Indeed, in the setting of hyperbolic manifolds, measured
geodesic laminations appear on boundaries of the convex core. The convex core boundary
is a totally geodesic pleated surface, bent along a geodesic lamination, and the
bending angles provide the transverse measure. Measured laminations appear similarly
in the negative curvature Lorentzian setting, again on the boundaries of the
convex core. The relevance of measured laminations to the dS setting can be expected
by duality. Finally, in the Lorentzian flat context measured laminations appear
by considering {\it flat regular domains}, see \cite{Benedetti} for details. The basic
idea here is that a lamination provides a way to deform the surface $\Sigma$ that
is embedded inside the future light cone of the Minkowski space $M^{1,2}$ as a
quotient $M^{1,2}/\Gamma$, where $\Gamma\subset{\rm SL}(2,\R)$ is the Fuchsian
group that uniformizes $\Gamma$. In all the cases considered one can view the
measured geodesic lamination on $\Sigma$ as specifying a set of initial data for
evolution in the corresponding space. In hyperbolic and AdS spaces this evolution starts
from the boundary of the convex core, and in the flat Lorentzian case it starts
from the level surface of the so called {\it canonical time} function. 

The above description is quite reminiscent of the one in this paper, except that instead
of measured geodesic laminations we use the cotangent bundle over Teichm\"uller space,
and instead of convex core boundaries we use a surface of zero, or constant, mean curvature.
Using this $T^* \cT_\Sigma$ structure one can similarly define a 
mapping, or ``analytic continuation'' between the corresponding spaces. The analytic
continuation arising is, however, rather different. First of all, there is no unique
map from flat Lorentzian manifolds to the non-flat setting anymore. Indeed, the phase
space $T^* \cT_\Sigma$ of the flat Lorentzian case is obtained by considering
a foliation by constant mean curvature surfaces. There is no preferred surface in 
this foliation, and thus, there is no preferred point in $T^* \cT_\Sigma$ that 
could be taken as the data for reconstruction of the non-flat spaces. However, one can take an
arbitrary value of $H$; the corresponding map gets therefore parameterized
by a real number $H\in(0,\infty)$. 

On the other hand, the other three settings are connected by canonical maps. Indeed, 
the hyperbolic and dS settings are related by
the usual duality. The results of this paper provide a natural map to the 
the setting of negative curvature Lorentzian manifolds. Indeed,
given a GHMC AdS manifold, there is a unique minimal surface in it, and 
it thus gives rise to a unique point in $T^* \cT_\Sigma$. One could try to use
this data for an analytic continuation to the hyperbolic setting. But,
as we know, only a bounded domain of $T^* \cT_\Sigma$ is relevant in
that setting. Thus, the analytic continuation would not be defined for
all GHMC AdS manifolds. However, there is another natural possibility.
Indeed, in both cases a complete description of the corresponding 
moduli spaces of manifolds is given by $\cT_\Sigma \times\cT_\Sigma$.
In the AdS case this is either due to Mess map, or, due to another
similar map that is described in this paper. In the hyperbolic or
dS settings this is due to Bers simultaneous uniformization. Thus,
the analytic continuation that we would like to propose relates
AdS with hyperbolic and dS settings by using the $\cT_\Sigma \times\cT_\Sigma$
description. 

Importantly, there are two different such analytic continuations that can be
defined. In one of them, to describe AdS manifolds by $\cT_\Sigma \times\cT_\Sigma$
one uses the metrics $I^\#_\pm = I((E\pm JB)\cdot,(E\pm JB)\cdot)$. Here
$I$ is the induced metric on some space-like surface $S$ in $M$, $B$ is the shape
operator of $S$, and $J$ is the operator $J^2=-1$ that gives the
complex structure of $S$. These two metrics are hyperbolic (i.e. constant curvature $-1$), and
give two points in the Teichm\"uller space $\cT_\Sigma$. It turns out that it does not  
matter which surface is used to get them. The construction of Mess referred to one
of the two boundaries of the convex core. 
The map that we introduce in this paper is quite similar,
but instead of an arbitrary surface, one uses a surface of vanishing mean curvature,
and considers $I^*_\pm = I((E\pm B)\cdot,(E\pm B)\cdot)$. Note that there
is no longer the operator $J$ in these metrics. These two metrics are also
hyperbolic, and also give two points in $\cT_\Sigma$. These are different two
points as compared to those given by the Mess map. Both pairs of points in $\cT_\Sigma$
are canonically defined given a GHMC AdS manifold, and both can be used to get
a hyperbolic or a dS manifold. It would be of great interest to compare these two
maps with the analytic continuation of Benedetti and Bonsante. It would also
be important to understand the difference between them, and, eventually,
understand which of these two maps is the physically ``right'' one. To
this end, it is probably of importance to extend the analysis of this paper
to the case of open $\Sigma$, for it is this case that describes the
physically interesting spacetimes containing black holes.

\paragraph{Point particles.}

It is an old and well-known problem to describe the gravity theory in 3 dimensions when 
conical defects are allowed to be present. Conical defects correspond to point
particles, and their inclusion is important to make the theory more realistic.
In the case when there are no conical defects the holonomy group describes
the 3-manifold completely. In particular, the manifold $M$ can be obtained as
the quotient of the (domain where the action is properly discontinuous) of the
maximally symmetric space by the holonomy group $\Gamma$, which is discrete.
When there are point particles, the group $\Gamma$ is no longer discrete, and
the 3-manifold is no longer the quotient. This makes the techniques that work
in the no-particle case essentially useless. 

The description of constant curvature 3-manifolds that we give in this paper is
based on equidistant foliations. A large class of 3-manifolds gets parameterized
by (in general a subset of) the cotangent bundle over the Teichm\"uller space of
the spatial slice. It is natural to try to extend this formulation to the
situation when point particles are present. We show that this can indeed be done, so that the
phase space in case point particles are present is (possibly a subset of) the
cotangent bundle $T^* \cT_{g,n}$ of the Teichm\"uller space of a surface with marked
points. We extend this description to all the cases considered in this paper, 
i.e. to hyperbolic, AdS, dS and Minkowski manifolds. 
However, there is a price to pay to have such a description: it turns out that one 
must restrict the range of angle deficits that
are allowed to $(\pi,2\pi)$. Thus, in particular, the small deficit angles,
which describe the ``almost'' no particle case is not covered by this
description. The reason for this 
is that when the total angle is large (angle deficit is small), the principal
curvatures of the surface turn out to diverge at the singular points. 
This makes the description based on equidistant foliations essentially
useless, because the equidistant foliation will break quite close to the
minimal surface. Thus, it turns out that the description based on
the first and second fundamental forms of a constant mean curvature surface
does not work when the total angle at the singular points is greater than
$\pi$.

However, it turns out that there is a certain ``dual'' description that ``works'' 
for the case total angles between $\pi$ and $2\pi$. 
Let us take a moment to explain what
this dual description is. In the original description, to which
most of the present paper is devoted, one uses the first and second
fundamental forms of a constant mean curvature surface $S$ as the data.
These two forms have to satisfy the Gauss and Codazzi equations,
and leave one with the free data of a point in Teichm\"uller space
of the surface $S$, together with a holomorphic quadratic differential on $S$.
However, one can instead use the dual description that is based on the
so-called third fundamental form. The later is defined as: 
$\III(x,y)=I(Bx,By)$,
where $x,y$ are two arbitrary tangent vectors to $S$. The third and second
fundamental forms satisfy similar Gauss and Codazzi equations, and again
leave one with the free data being a point in Teichm\"uller space together
with a holomorphic quadratic differential. The third fundamental form
plays the role of the ``dual'' metric in the sense that will be explained
in more details in the main text. Now the point is that, similarly to
what happened in the $(I,\II)$ description, given a point
in $T^* \cT_{g,n}$ one can reconstruct the third fundamental form $\III$.  
One can now view the third fundamental form 
as a metric induced on a ``dual'' surface, and the principal
curvatures of this dual surface remain finite at singular
points. Using these data, one can then reconstruct what
can be called the dual manifold. The relation with the
``original'' manifold is not simple, however, and is not
dealt with in the present paper. 

Thus, the results of this paper give a satisfactory description of 
cone manifolds when all the total angles are in $[0,\pi]$. When the
angles are in $[\pi,2\pi]$ the above description does not work, one has
to work with $\III,\II$ forms and the ``dual'' 3-manifold. Note, however, that both 
of these descriptions become
invalid when particles of both types are present at the same time. 
We do not consider this case in the present paper.

To summarize, as our results show, there is a big difference between 
particles with the total angle less and greater than $\pi$. In a certain
sense, these two types of particles behave like different species,
and one has to use a different description to deal with each type.
There is a more or less satisfactory description of the cone manifolds 
if there is only a single type of particles present. It remains to
be seen how to describe the general situation. Or, more speculatively,
it may be that the dual types of particles are on the same
footing as particles and anti-particles of the quantum field
theory, and one never has to consider them together at the classical
level. Thus, it may be that, as far as the classical theory is concerned,
the description we give in this paper is sufficient, and both
types of particles should only be considered together in 
an appropriate quantum theory of 3d gravity. Whether
this speculation is indeed realized remains to be seen.

\section{Minimal surfaces in quasi-Fuchsian manifolds} 

\subsection{Minimal surfaces in germs of manifolds}

Throughout this paper $\Sigma$ is a closed, orientable surface of
genus at least 2, unless specified otherwise.
 
\paragraph{Minimal surfaces in hyperbolic manifolds.}

Taubes recently introduced a convenient notion of ``the 
minimal surface in a germ
of a hyperbolic manifold''. Considered on $\Sigma$, this is a couple $(g,h)$,
where $g$ is a smooth Riemannian metric and $h$ is a symmetric bilinear form
on $T\Sigma$, such that $h$ is the second fundamental form of a minimal
isometric embedding of $\Sigma$ in a (possibly non-complete) hyperbolic
3-manifold. By the ``fundamental theorem of surface theory'' (see
e.g. \cite{spivak}) this is equivalent to the following 3 conditions: 
\begin{enumerate}
\item The trace of $h$ with respect to $g$ vanishes: $\tr_gh=0$.
\item $h$ satisfies the Codazzi equation with respect to the Levi-Civit\`a
  connection $\nabla$ of $g$: $d^\nabla h=0$.
\item The determinant of $h$ with respect to $g$ satisfies the Gauss equation: 
$K_g=-1+\det_g(h)$. 
\end{enumerate}
In condition (2), $d^\nabla$ is the exterior derivative, associated to
$\nabla$, on the bundle of 1-forms with values in the tangent space of
$\Sigma$, and $h$ is considered as such a vector-valued 1-form.

Given an immersion of $\Sigma$ in a hyperbolic 3-manifold, its shape operator
$B:T\Sigma\rightarrow T\Sigma$ is the unique self-adjoint operator such that: 
$$ \forall m\in \Sigma, \forall x,y\in T_m\Sigma, \II(x,y) =
I(Bx,y)=I(x,By)~. $$
The immersion is minimal if and only if $B$ is traceless.

Following the notations in \cite{taubes}, we call $\cH_g$ the space of 
minimal surfaces in germs of hyperbolic manifolds, i.e. the space
of pairs  $(g,h)$ satisfying the three conditions above.

\paragraph{The index of minimal surfaces.}

Let $u:\Sigma \rightarrow M$ be a minimal immersion of $\Sigma$ in a 
hyperbolic 3-dimensional manifold. The first-order deformations of $u(\Sigma)$
are determined by vector fields tangent to $M$ defined along $u(\Sigma)$. 
However, vector fields which are tangent to $u(\Sigma)$ act trivially, 
so that it is sufficient to consider vector fields which are orthogonal to 
$u(\Sigma)$. Since $\Sigma$ is orientable, it is possible to choose a 
unit normal vector field $N$ on $u(\Sigma)$, and vector fields orthogonal to
$u(\Sigma)$ are then of the form $fN$, for a function $f:\Sigma\rightarrow \R$.

The hypothesis that $u$ is a minimal immersion translates as the fact that 
the area does not vary, at first order, under such a deformation. It is
therefore 
possible to consider the second variation of the area, which is given by a 
well-known integral formula (a more general case can be found e.g. in 
\cite{berard-docarmo-santos}).

\begin{lemma} \label{lm:second}
The second variation of the area under the first-order deformation $fN$ is
given by: 
\begin{equation}\label{sec-var}
A''(f) = \frac{1}{2} \int_\Sigma \left( f\Delta f + 2(1+K_e) f^2 \right) da~, 
\end{equation}
where the area form $da$ and the Laplace operator 
$\Delta f=-g^{ab}\nabla_a \nabla_b f$
come from the metric $g$ induced by
$u$ on $\Sigma$, and $K_e$ is the extrinsic curvature, i.e. the 
product of two principal curvatures of $\Sigma$.
\end{lemma}

A minimal surface is called {\bf stable} when the integral is strictly
positive  
for any non-zero function $f$. If a minimal surface is stable, then it is
a local minimum of the area. Conversely, for any local minimum of the 
area, the quantity $A''(f)$ is non-negative for any function $f$.

\paragraph{Some Teichm\"uller theory.}

In all the paper, we call $\cT_g$ the Teichm\"uller space of genus $g$, 
for $g\geq 2$, i.e. the space of conformal structures on a surface of
genus $g$, considered up to the diffeomorphisms isotopic to the identity.
It will be convenient to set $\cT_\Sigma:=\cT_g$, where $g$ is the genus of
$\Sigma$.

Since $\Sigma$ is orientable, conformal structures on $\Sigma$ are in
one-to-one 
correspondence with complex structures on $\Sigma$. Given such a complex
structure, 
it is possible to consider the holomorphic quadratic differentials 
(called HQD below) on $\Sigma$. Given a conformal structure $c\in \cT_\Sigma$, 
the space of HQD on $\Sigma$ for $c$ is canonically identified with the 
cotangent bundle of $\cT_\Sigma$ at $c$, see e.g. \cite{ahlfors}. 

\paragraph{Minimal surfaces from HQDs}

There is a striking relation between HQDs and minimal surfaces, basically
going back to Hopf \cite{hopf} but recently  exposed by Taubes \cite{taubes}.
It is valid for minimal surfaces in constant curvature 3-manifolds.

\begin{lemma} \label{lm:hqd}
Let $g$ be a Riemannian metric on $\Sigma$, and let $h$ be a 
bilinear symmetric form on $T\Sigma$. Then:
\begin{enumerate}
\item The trace of $h$ with respect to $g$, $\tr_g(h)$, is zero if 
and only if $h$ is the real part of a quadratic differential $q$ over 
$\Sigma$.
\item If (1) holds, then $q$ is holomorphic if and only if $h$ satisfies the 
Codazzi equation, $d^\nabla h=0$.
\item If (1) and (2) hold, then $(g,h)$ is a minimal surface in a germ
of hyperbolic manifold if and only if the Gauss equation is satisfied, 
i.e. $K_g=-1 + \det_g(h)$. 
\end{enumerate}
\end{lemma}

A direct but remarkable consequence of point (2) is that, if $\tr_g(h)=0$,
then the Codazzi equation is invariant under conformal deformations of $g$.

This lemma
shows that, for all $(g,h)\in \cH_g$, $h$ is the real part of a QHD on $\Sigma$
with the complex structure coming from $g$, and this defines a natural map
\begin{equation}\label{h-tt}
\phi:\cH_g\rightarrow T^*\cT_g~.
\end{equation}

The third condition in the previous lemma can be realized, in some cases, 
by a conformal change in $g$. 

\begin{lemma} \label{lm:confchange}
Let $(g,h)$ be as in Lemma \ref{lm:hqd}, suppose that conditions (1) and 
(2) of that lemma are satisfied. Let $g':=e^{2u}g$. Then conditions (1) 
and (2) are also satisfied for $(g',h)$, and condition (3) holds if and
only if:
\begin{equation} \label{eq:*}
\Delta u = - e^{2u} - K_g + e^{-2u} det_g(h)~. 
\end{equation}
\end{lemma}

\begin{proof}
An elementary scaling argument shows that: $\det_{g'}(h) = e^{-4u}\det_g(h)$. 
Moreover, a well-known formula (see e.g. \cite{Be}, ch. 1) expresses the 
curvature $K_{g'}$ of $g'$ as: 
$$ K_{g'} = e^{-2u} (\Delta u + K_g)~. $$
Those two equations can be used to formulate condition (3) as claimed.
\end{proof}

\paragraph{CMC surfaces.}

It is interesting to note that Lemma \ref{lm:hqd}, and subsequently Lemma
\ref{lm:confchange}, extend from minimal to constant mean curvature (here
abbreviated as CMC) surfaces. Here we defined the {\it mean curvature} of a
surface as $H=\tr(B)/2$ (some definitions differ by the coefficient $1/2$),
and a CMC surface is a surface with $H$ equal to a constant. A {\it CMC
$H$ surface in a germ of hyperbolic manifold} is defined, as for minimal
surfaces, as a couple $(g,h)$, where $g$ is a metric on $\Sigma$, $h$ is a
symmetric bilinear form on $T\Sigma$, and $g$ and $h$ satisfy the same
equation as the first and second fundamental form of a CMC surface in $H^3$:
\begin{itemize}
\item $\tr_g(h)=2H$.
\item The Codazzi equation $d^\nabla h=0$
\item The Gauss equation: $K_g=-1+\det_g(h)$.
\end{itemize}

\begin{lemma} \label{lm:cmc}
Let $(g,h)$ be as in Lemma \ref{lm:hqd}, and let $h_0:=h-Hg$. Then:
\begin{enumerate}
\item The trace of $h_0$ with respect to $g$, $\tr_g(h_0)$, is zero if 
and only if $h_0$ is the real part of a quadratic differential $q$ over 
$\Sigma$.
\item If (1) holds, then $q$ is holomorphic if and only if $h_0$ satisfies the 
Codazzi equation, $d^\nabla h_0=0$.
\item If (1) and (2) hold, then $(g,h)$ is a CMC 
$H$ surface in a germ
of hyperbolic manifold if and only if:
$$ K_g = (H^2-1) + det_g(h_0)~. $$ 
\item If $g':=e^{2u}g$, then point (3) holds for $(g',h)$ if and only if $u$
  is a solution of: 
\begin{equation} \label{eq:*cmc}
\Delta u = (H^2 - 1) e^{2u} - K_g + e^{-2u} det_g(h_0)~. 
\end{equation}
\end{enumerate}
\end{lemma}

\begin{proof}
Points (1) and (2) are already present in Lemma \ref{lm:hqd}, because
$d^\nabla g=0$. For point (3) note that the Gauss equation can be written as: 
$$ K_g = -1 + det_g(Hg+h_0) = -1 + H^2 + H\tr_g(h_0) + det_g(h_0)~, $$
while $\tr_g(h_0)=0$ if and only if $\tr_g(h)=2H$, so if and only if $(g,h)$
is a CMC $H$ surface in a germ of hyperbolic manifold. 

Point (4) is proved exactly as Lemma \ref{lm:confchange}. 
\end{proof}

\subsection{Almost-Fuchsian manifolds}

\paragraph{Minimal surfaces in quasi-Fuchsian manifolds.}

Consider a quasi-Fuchsian metric on $M=\Sigma\times \R$. Since the metric is 
convex co-compact, it is not difficult to check that surfaces isotopic to
$\Sigma \times \{ 0\}$ which have area bounded by some constant $A_0$ 
are always contained in a compact subset $K(A_0)$. It follows, using
standard (but non-trivial) results on minimal surfaces, that $M$
contains at least one minimal, area-minimizing surface isotopic to
$\Sigma \times \{ 0\}$. It is conceivable, however, that $M$ contains
many such surfaces, and we will see below that this can actually happen.

\paragraph{Almost-Fuchsian manifolds.}

It is interesting to consider a sub-class of the quasi-Fuchsian manifolds, 
defined in terms of the minimal surfaces which they contain. 

\begin{df}
A quasi-Fuchsian hyperbolic 3-manifold is {\bf almost-Fuchsian} if it contains
a closed embedded surface with principal curvatures in $(-1,1)$.  
\end{df}

According to a recent result of Ben Andrews (see \cite{rubinstein}, ch. 8) any 
closed surface with principal curvatures in $(-1,1)$ in a complete hyperbolic 
3-manifold can be deformed to a minimal surface with principal curvatures in 
$(-1,1)$. So, in the definition above, we could have replaced ``a closed
embedded surface'' by ``a closed {\it minimal} embedded surface''. This leads
to the strongly related definition.

\begin{df}
Let $(g,h)\in \cH_g$ be a minimal surface in a germ of hyperbolic manifold. It
is {\it almost-Fuchsian} if the principal curvatures of its shape operator are
everywhere in $(-1,1)$. The set of such almost-Fuchsian minimal surfaces is
denoted by $\cH_g^{af}$.
\end{df}

An important property of almost-Fuchsian manifolds, noted by Uhlenbeck
\cite{uhlenbeck}, is that
they admit a (smooth) foliation by surfaces equidistant from the minimal 
surface appearing in their definition. This can be seen as a consequence of 
an elementary (and well-known, see e.g. \cite{gray}) 
statement on equidistant foliations from surfaces in $H^3$.

\begin{lemma} \label{lm:27}
Let $S\subset H^3$ be a complete, oriented, smooth surface with principal
curvatures in $(-1,1)$.
For all $r\in \R$, the set of
points $S_r$ at oriented distance $r$ from $S$ 
is a smooth embedded surface. The closest-point
projection defines a smooth map $u_r:S_r\rightarrow S$, 
and, identifying
the two surfaces by $u_r$, the induced metric $I_r$ on $S_r$ is:
\begin{equation}\label{metric}
I_r(x,y) = I((\cosh(r)E+\sinh(r)B)x,(\cosh(r)E+\sinh(r)B)y)~,
\end{equation}
where $E$ is the identity operator. The shape operator of $S_r$ is: 
\begin{equation}\label{shape}
B_r:=(\cosh(r) E +\sinh(r)B)^{-1}(\sinh(r)E+\cosh(r)B)~. 
\end{equation}
\end{lemma}

\begin{proof} Because of the importance of this lemma for what follows, we
give an explicit proof. The shape operator $B_r$ must satisfy the following
differential equation:
$$ \frac{d B_r}{d r} = E - B_r^2~.$$
On the other hand, $B=I^{-1} \II$, and by definition:
$$ \II_r = \frac{1}{2} \frac{d I_r}{dr}~.$$
These equations imply the following differential equation for the induced
metric: 
\begin{equation}\label{eqn-metric}
(I_r)'' - \frac{1}{2} (I_r)' I_r^{-1} (I_r)' = 2I_r,
\end{equation}
where the prime denotes the derivative with respect to $r$. This equation can
be 
solved by differentiating it one more time. One gets: $(I_r)'''=4I_r$ or,
equivalently: $(\II_r)''=4\II_r$. The most general solution of this, and
the corresponding expression for the induced metric is:
\begin{eqnarray}
\II_r = \alpha e^{2r} + \beta e^{-2r}, \\ \nonumber
I_r = \alpha e^{2r} - \beta e^{-2r} + \gamma,
\end{eqnarray}
where $\alpha,\beta$ are arbitrary operators and $\gamma$ is the 
integration ``constant''. All these quantities can be expressed
in terms of the first $I$ and second $\II$ fundamental forms of the surface $S$,
which we place at $r=0$. We have:
$$ \II=\alpha+\beta, \qquad I = \alpha-\beta + \gamma. $$
The equation (\ref{eqn-metric}) implies an additional relation between 
$\alpha, \beta, \gamma$. One gets:
$$ \gamma = \frac{1}{2} (I- \II I^{-1} \II).$$
Actually, $\gamma$ is given by the same expression in terms of the
fundamental forms for any $r$ and is independent of $r$ because of
(\ref{eqn-metric}). Collecting all these facts, the induced
metric and the second fundamental forms of $S_r$ can be written as:
\begin{eqnarray}
I_r = (\cosh(r) I  + \sinh(r) \II) I^{-1} (\cosh(r) I  + \sinh(r) \II), \\ \nonumber
\II_r = (\cosh(r) I  + \sinh(r) \II) I^{-1} (\cosh(r) \II  + \sinh(r) I).
\end{eqnarray}
The expressions (\ref{metric}) and (\ref{shape}) follow.

If the principal curvatures of $S$ are everywhere in $(-1,1)$, then $\cosh(r)E
+\sinh(r)B$ is always non-singular, so that $\Sigma_r$ is always smooth, and
the projection from $S$ to $S_r$ along geodesics orthogonal to $S$ is a
diffeomorphism. 
\end{proof}

\begin{cor}
Let $M$ be an almost-Fuchsian manifold diffeomorphic to $\Sigma\times \R$, and
let $S\subset M$ be a minimal surface isotopic to $\Sigma \times\{0\}$, with
principal curvatures in $(-1,1)$. Let $I$ and $B$ be the induced metric and
shape operator of $S$, respectively. Then $M$ is isometric to $\Sigma
\times\R$ with the metric: 
\be\label{metric-1}
dr^2 + I((\cosh(r)E+\sinh(r)B)\cdot,(\cosh(r)E+\sinh(r)B)\cdot)~. 
\ee
It is foliated by the smooth surfaces $\Sigma_r:=\Sigma\times \{r\}$,
which have shape operator given (under the natural identification of
$\Sigma_r$ with $S$ by projection) by:
$$ B_r:=(\cosh(r) E + \sinh(r)B)^{-1}(\sinh(r)E+\cosh(r)B)~. $$
\end{cor}

It is not difficult to extract from the expression of the metric an explicit
expression of the conformal structure on each connected component of the
boundary at infinity of $M$. We will give these expressions below.

\begin{cor} \label{cr:unique-min}
Let $M$ be an almost-Fuchsian hyperbolic manifold, then $M$ contains only one
closed minimal surface (which is almost-Fuchsian).
\end{cor}

\begin{proof}
Let $\Sigma_r, r\in \R$ be the foliation of $M$ given by the previous
corollary. Since $S=\Sigma_0$ is minimal, its principal curvatures are
opposite numbers $k$ and $-k$, for some $k\geq 0$. 
It follows from the expression of $B_r$ that the mean curvature of
$\Sigma_r$ is equal to:
\begin{eqnarray*}
2H_r & = & \tr((\cosh(r) E + \sinh(r)B)^{-1}(\sinh(r)E+\cosh(r)B)) \\
& = & \frac{\sinh(r) + k\cosh(r)}{\cosh(r) + k\sinh(r)} + \frac{\sinh(r) -
  k\cosh(r)}{\cosh(r) - k\sinh(r)} \\
& = & \frac{2 (1-k^2)\tanh(r)}{1-k^2\tanh^2(r)}~. 
\end{eqnarray*}
Since $\tanh(r)\in (-1,1)$ and $k\in [0,1)$, the denominator is always
positive, while the numerator has the same sign as $r$. So the surfaces
$\Sigma_r$ has strictly positive mean curvature for $r>0$ and strictly
negative curvature for $r<0$.

Let $S'\subset M$ be a closed minimal surface. Let $r_1$ be the largest $r$
such that $S'\cap \Sigma_r\neq 0$. Then $S'$ is tangent to $\Sigma_{r_1}$ 
on the side of decreasing $r$, and, since $S'$ is minimal and $H_r$ has the
same sign as $r$, this is possible only if $r_1\leq 0$. Conversely, if $r_2$ is
the smallest $r$ such that $S'\cap \Sigma_r\neq 0$, then $r_2\geq 0$. So
$r_1=r_2=0$, and $S'=S$.
\end{proof}

\paragraph{Almost-Fuchsian manifolds from HQD.}

The definition of a map $\phi:\cH_g\rightarrow T^*\cT_g$ 
was given above (\ref{h-tt}). This map is in general neither injective
nor surjective. Let us first show that it is non-surjective. Given
a point in $T^* \cT_\Sigma$ one selects an arbitrary metric $g$ in the
corresponding conformal class and tries to satisfy the Gauss equation
by searching for the conformal factor $u$ that satisfies (\ref{eq:*}).
Let us integrate this equation over the minimal surface. We get:
$$ \int_\Sigma (e^{2u} + k^2 e^{-2u})da = - \int_\Sigma K_g da = 2\pi(2g-2).$$
But the integrand on the right hand side is greater or equal to $2k$. Thus, we
get: 
$$ \int_\Sigma k da\leq 2\pi(g-1).$$
Now, as we know from Lemma \ref{lm:hqd}, the second fundamental form of
$\Sigma$ 
is the real part of a holomorphic quadratic differential $t dz^2$. As we will
show 
in the next subsection, this implies that 
$k da = \sqrt{t\bar{t}} |dz|^2$, so the above inequality expressed in terms of
the quadratic differential $t$ becomes:
\begin{equation}
\int_\Sigma \sqrt{t\bar{t}} |dz|^2 \leq 2\pi(g-1).
\end{equation}
This clearly shows that solutions of (\ref{eq:*}) can only exist in an
open subset of $T^* \cT_\Sigma$. Thus, $\phi$ is non-surjective. It is also
generally non-injective, because even when there is a solution of  
equation (\ref{eq:*}) it is in general non-unique. That this indeed happens is
suggested by the following lemma (compare with a similar lemma of Taubes
\cite{taubes}). 

\begin{lemma} 
The map $\phi$ fails to be a homeomorphism when the quadratic
form that is given by the second variation (\ref{sec-var}) has zero
eigenvalues. 
\end{lemma}

\begin{proof} 
Let us consider the second variation (\ref{sec-var}) of the area of 
the surface with the metric $g'=e^{2u} g$, where $u$ satisfies
(\ref{eq:*}). Using the 
fact that $\Delta'=e^{-2u}\Delta, K_e'=e^{-4u} K_e$ we get:
$$ A''(f) = \frac{1}{2} \int_S \left( f\Delta f + 2 e^{2u}(1+ e^{-4u} K_e)f^2
\right).$$ 
This quadratic form has zero eigenvalues iff the equation:
$$ \Delta f + 2(e^{2u}+ e^{-2u} K_e) f =0$$
has non-trivial solutions. However, because $K_e=\det_g(h)$ the above equation
is the same as the 
one obtained by considering the first 
variation $u\to u+f$ of the equation (\ref{eq:*}). This implies the statement
of the lemma. 
\end{proof}

In spite of all these nasty features of the map $\phi$, it does behave rather
nicely on the subspace 
of almost-Fuchsian minimal surfaces. In particular, it is injective on
$\cH_g^{af}$, as we shall now show. 

\begin{lemma} \label{unique} The map $\phi$ is injective when restricted to
minimal surfaces in germs of almost-Fuchsian manifolds. In other words,
let $(g,h),(g',g')\in \cH_g^{af}$ be such that
$\phi(g,h)=\phi(g',h')$. Then $(g,h)=(g',h')$.
\end{lemma}

\begin{proof}
By definition, the element of $T^*\cT_g$ is the same, which means that, $g'$
is conformal to $g$ while $h=h'$. So there exists $u:\Sigma\rightarrow \R$
such that $g'=e^{2u}g$, and $u$ is a solution of (\ref{eq:*}), which can be
written as:
$$ \Delta u = -e^{2u} + 1 + k^2 - k^2 e^{-2u} = (1-e^{2u}) + k^2(1-e^{-2u})~,
$$ 
where $k$ is the largest principal curvature of $(g,h)$ at each point of
$\Sigma$. Let:
$$ f(u):= (1-e^{2u}) + k^2(1-e^{-2u})~. $$
Then $f$ is the sum of two strictly concave functions, so it is concave, and
moreover $f(0)=0$ and $f'(0)=2(k^2-1)<0$, because $(g,h)$ is almost-Fuchsian.
It follows that $f(u)$ is negative at all points of $\Sigma$ where $u>0$.

At points where $u$ attains its maximum, $\Delta u\geq 0$, so that $f(u)\geq 0$
and therefore $u\leq 0$ -- and therefore $u\leq 0$ everywhere on $\Sigma$. But
the same argument can be applied with $(g,h)$ and $(g',h')$ exchanged, and $u$
replaced by $-u$, and it shows that $-u\leq 0$. So $u=0$ on $\Sigma$, and
$g'=g$. 
\end{proof}

This, along with Corollary \ref{cr:unique-min}, shows that almost-Fuchsian
hyperbolic 
manifolds are parametrized by an open subset of the cotangent space of
Teichm\"uller space. This can be summed up as follows.

\begin{thm} \label{tm:hyperbolic}
There exists an open subset $\Omega_g\subset T^*\cT_g$ such that:
\begin{itemize}
\item For each $\tau\in T^*\cT_g$, there exists a unique almost-Fuchsian
manifold $M\in \cM_g^{af}$ 
containing a closed minimal surface $S$ such that, if $I$ and
$\II$ are the first and second fundamental forms of $S$, then $\tau =
  \phi(I,\II)$. 
\item Conversely, any almost-Fuchsian manifold $M\in \cM_g^{af}$ 
contains a unique closed
minimal surface $S$, and $\phi(I,\II)\in \Omega_g$, where $I$ and $\II$ are 
the induced metric and second fundamental forms of $S$.
\end{itemize}
\end{thm}

\paragraph{The Fock metric on almost-Fuchsian manifolds.} 

The fact that for a minimal surface the second fundamental form 
is the real part of some holomorphic quadratic differential
can be further exploited to write the metric (\ref{metric-1}) in a nice 
complex-analytic way. This way of writing (\ref{metric-1}) was discovered
by Fock \cite{Fock}. 

Let us choose the Fuchsian uniformization of the minimal surface $S$ by the
hyperbolic plane $H_2$, and let $z$ be the usual complex coordinate in $H_2$.
The group $\Gamma$ of deck transformations acts properly discontinuously on
$H_2$ and $S=H_2/\Gamma$. The metric on $S$ can be written as: $g=e^\varphi
|dz|^2$, 
where the Liouville field $\varphi=\varphi(z,\bar{z})$ has the following
transformation properties: 
\be
\varphi(\gamma\circ z,\overline{\gamma\circ z}) = \varphi(z,\bar{z}) -
\log{|\gamma'|^2},  
\qquad \forall\gamma\in\Gamma.
\ee 
Here $\gamma\circ z$ is the action of an element of $\Gamma$ on $H_2$, and
$\gamma'$ is 
the derivative with respect to $z$ of $\gamma\circ z$ considered as a function
of $z$. 

Let us also introduce a holomorphic quadratic differential $t$ on $S$, namely a
holomorphic function on $H_2$ with the following transformation properties:
\be
t(\gamma\circ z) = t(z) (\gamma')^{-2}, \qquad \forall\gamma\in\Gamma.
\ee
The second fundamental form on $S$ in terms of $t$ is then $h=tdz^2 +
\bar{t}d\bar{z}^2$. 
The metric (\ref{metric-1}) takes the following simple form, first discovered
by Fock \cite{Fock}: 
\be\label{fock}
ds^2 = dr^2 + e^\varphi |\cosh(r)dz + \sinh(r)e^{-\varphi} 
\bar{t} d\bar{z}|^2~.
\ee 
It is clear from the transformations properties of $\varphi, t$ that the metric
on surfaces $r=const$ is invariant under $\Gamma$ and thus descends to a
metric on the equidistant surfaces $S_r$.

The equation (\ref{eq:*}) takes the following form:
\be\label{eq:**}
2 \partial_{z\bar{z}} \varphi = e^\varphi + e^{-\varphi} t\bar{t}.
\ee
Is is easy to verify that the principal curvatures of the minimal surface $S$ in
the parameterization used are: $\lambda_{1,2}=\pm e^{-\varphi} \sqrt{t\bar{t}}=\pm k$,
where $k$ is introduced for future convenience.
Therefore, the extrinsic curvature is given by: $K_e =-k^2$.
The area form on surfaces $S_r$ is given by:
\be
da_r = (\cosh^2(r) - \sinh^2(r) k^2)da.
\ee
For almost-Fuchsian manifolds $k<1$ everywhere on $S$ and $da_r$ never goes to zero. 

It is easy to write down the expression for the metrics 
on two conformal boundaries of $M$. As is easy to see, the metric on the 
$S_r$ surfaces is asymptotic to $(1/4) e^{2r} I^*$,
with the asymptotic metric given by:
\be
I^*_{1,2} = e^\varphi |dz \pm e^{-\varphi} \bar{t} d\bar{z}|^2.
\ee 
It is clear from this expression that the deformation of the complex structure of $S$
that resulted by following the equidistant $S_r$ surfaces is that corresponding to 
the Beltrami differential:
\be\label{Beltrami}
\mu(z,\bar{z})= e^{-\varphi} \bar{t},
\ee
so that the asymptotic metrics are: $I^*_{1,2}=e^{\varphi} |dz\pm \mu d\bar{z}|^2$.
For almost-Fuchsian manifolds $|\mu|^2<1$ everywhere on $S$, and exactly for such
Beltrami differentials the Beltrami equation $f_{\bar{z}} = \mu f_z$ has a unique
(up to conjugation by ${\rm PSL}(2,\R)$) solution in the class of quasi-conformal
mappings of $H_2$ into itself. 

\paragraph{A map: $T^* \cT_\Sigma \to \cT_\Sigma\times \cT_\Sigma$.}

Theorem \ref{tm:hyperbolic} states that almost-Fuchsian manifolds are parametrized
by an open subset in $T^* \cT_\Sigma$. On the other hand, by Bers simultaneous uniformization, 
quasi-Fuchsian manifolds $M$ are in one-to-one correspondence with points of
$\cT_\Sigma\times \cT_\Sigma$. Thus, in the almost-Fuchsian case one can combine the inverse of the
map $\phi$ with the Bers map to get a map:
$$\Phi: T^* \cT_\Sigma \to \cT_\Sigma\times \cT_\Sigma~.$$

Using the above complex-analytic description, this map can be 
described as follows. One starts with a point in $T^* \cT_\Sigma$.
The base point is described by the corresponding Fuchsian group $\Gamma$, and the
point in the fiber is a quadratic differential $t$ for $\Gamma$. To find the metric on
$\Sigma$ and in $M$ one has to solve the equation (\ref{eq:**}). For small enough $t$
the solution exists, and its uniqueness is guaranteed by the lemma \ref{unique}. Once the
metric on $\Sigma$ is found by solving (\ref{eq:**}), one immediately obtains the
metric in the whole $M$ via (\ref{fock}). The complex structure at each conformal infinity is the one
obtained by a quasi-conformal deformation $f_{\bar{z}} = \mu f_z$, with the
Beltrami differential given by plus-minus (\ref{Beltrami}). The map
$\Phi$ is non-trivial and involves solving the equation (\ref{eq:**}). 

\paragraph{An upper bound for the convex core volume.} 

As a simple application of the equidistant to a minimal surface description we
derive an upper bound for the volume of the convex core of an almost-Fuchsian manifold.
The principal curvatures of the $S_r$ surfaces are easily found to be given by:
\be
\lambda_{1,2} = \frac{\cosh(r)\sinh(r)(1-k^2)\pm k}{\cosh^2(r) - \sinh^2(r) k^2}.
\ee
It is worth noting that the principal directions on $S$ are just the horizontal
and vertical trajectories of $t$, as is not hard to verify.

Both of the principal curvatures become positive for large enough $r$, so for
such large $r$ 
the surface $S_r$ is convex. Let us find the value of $r$ when the equidistant
surface 
is already convex. This happens when the principal curvatures on $S_r$ are everywhere non-negative.
Let $k_{max}$ denote the maximum of the principal curvature on $S$. Then the principal
curvatures are everywhere non-negative for
$\cosh(r_{max})\sinh(r_{max})=k_{max}/(1-k_{max}^2)$ or:
\be
e^{2r_{max}}=\frac{1+k_{max}}{1-k_{max}}.
\ee
An upper bound for the convex core volume is given by the volume in between the
two surfaces $r=\pm r_{max}$. One gets:
$$ V(r_{max},-r_{max}) = \int_{-r_{max}}^{r_{max}} dr \int_S da_r = 
\frac{1}{2}\sinh(2r_{max})\int_S (1-k^2)da - r_{max}\int_S (1+k^2)da.$$
After some simple algebra we get:
\be
V(core)<
\frac{2k_{max}}{1-k_{max}^2} A(S)-2\pi(2g-2)\left( \frac{k_{max}}{1-k_{max}^2} + 
\frac{1}{2}\ln{\frac{1+k_{max}}{1-k_{max}}} \right).
\ee
Here $A(S), g$ is the area and genus of the minimal surface correspondingly. 
To write this final formula we have used the fact that 
$\int_S (1+k^2)da = -2\pi(2-2g)$, which follows from
(\ref{eq:*}). For small $k_{max}$ the quantity on the right hand side goes as:
$2k_{max}(A(S)-2\pi(2g-2)) + O(k_{max}^2)$. Note, however, that the quantity in
brackets here vanishes in the Fuchsian case $k_{max}=0$. Thus, the convex core
volume is $o(k_{max})$ for small $k_{max}$.
As our bound shows, the quantity $k_{max}$ can be used as a measure of how 
far one is from the Fuchsian case.

\subsection{Non-almost-Fuchsian manifolds}

Unfortunately, the considerations made above for almost-Fuchsian manifolds do
not apply more generally to all quasi-Fuchsian manifold. We point out in this
subsection that there are some quasi-Fuchsian but not almost-Fuchsian
manifolds, and that things then don't work as well as could be hoped -- this
is in striking contrast with the situation for AdS, dS or Minkowski manifolds, as
described in the following sections.

\paragraph{Hyperbolic 3-manifolds which fiber over the circle.}

Thurston (see \cite{thurston-notes,otal-hyperbolisation,kapovich}) gave a
general theorem on the existence of hyperbolic metrics on ``most'' 3-manifolds
which fiber over the circle. It was at some point conjectured that all
hyperbolic manifolds which fiber over the circle are foliated by minimal
surfaces. This conjecture turned out to be false, as pointed out recently by
Hass and Thurston and by Rubinstein.

Since the argument is both very nice and fairly simple, and since it was
apparently not published, we give here a very broad and imprecise outline,
which we hope can indicate to the reader how to build a complete proof. 
Let $M$ be a
closed manifold which fibers over the circle, suppose that it has a foliation
by minimal surfaces $S_t, t\in S^1$. Then the maximum principle (basically the
same argument as used in the proof of Corollary \ref{cr:unique-min}) shows that
any closed minimal surface $S\subset M$ has to be one of the $S_t$. But $M$
contains at least one area-minimizing surface, and all the $S_t$ have the same
area (since they are all critical points of the area), so all the $S_t$ are
area-minimizing. 

Now let $c$ be a closed curve in $M$ which is isotopic to a closed curve in
one of the $S_t$, let $M_c$ be the finite-volume hyperbolic manifold obtained
by ``drilling'' $c$ out of $M$; then $M_c$ still fibers over the circle 
(see \cite{thurston-notes,gromov-bourbaki}). 
Let $M_n$ be a sequence
of hyperbolic manifolds obtained from $M_c$ by Dehn surgery with slopes going
to infinity (see \cite{thurston-notes}). Then the $M_n$ contain a sequence of
compact subsets $K_n\subset M_n$ converging to any compact subset of
$M_c$. Moreover, the $M_n$ also fiber over the circle, so that they should
have a foliation by minimal surfaces, which would have to be
area-minimizing. One of those area-minimizing surfaces would have to go
arbitrarily far (for $n$ large) in the part which is ``beyond'' the approximation 
of the cusp, and it is not difficult to convince oneself that it could then
not be area-minimizing, leading to a contradiction. 

\paragraph{Quasi-Fuchsian manifolds with many minimal surfaces.}

Since hyperbolic manifolds which fiber over the circle do not, in general,
have a foliation by minimal surface, there is one such manifold, $M$, which
contains an isolated area-minimizing closed surface, say $S$. Let
$\overline{M}$ be the infinite cyclic cover of $M$ with fundamental group
$\pi_1(S)$, considered as a subgroup of $\pi_1(M)$. 
$\overline{M}$ contains a sequence of
area-minimizing surfaces $S_n, n\in \Z$, which are the lifts of $S$ in
$\overline{M}$. 

The proof of the hyperbolization theorem for manifolds which fiber over the
circle (more precisely the double limit theorem) shows that there is a
sequence of quasi-Fuchsian manifolds $N_n, n\in \N$, converging to
$\overline{M}$. For $n$ large enough, there is a (finite) subset of the $S_n$
which have corresponding surfaces in $M_n$ which are local minima of the
area. So, for $n$ large enough, $M_n$ contains several minimal surfaces, and
it can not be almost-Fuchsian by Corollary \ref{cr:unique-min}. 

\paragraph{What goes wrong.} We thus see that non almost-Fuchsian manifolds
are much more complex. There is more than one minimal surface. One can choose
any of the minimal surfaces and base the equidistant foliation on it, but
because the principal curvatures are no longer in $(-1,1)$ this foliation
becomes singular. It is illuminating to see where this happens. The metric
induced on $S_r$ becomes singular whenever $\cosh^2(r)/\sinh^2(r)=k^2$. This
first happens at:
\be
e^{2r_{max}}=\frac{k_{max}+1}{k_{max}-1}.
\ee 
Note that when the induced metric becomes singular only one of the 
principal curvatures (the negative one) diverges. The other, positive one
remains finite and 
is equal to: $\lambda_2=(1/2)(k+1/k)$. It is clear from this that the foliation
breaks down when still inside the convex core. 

It is also interesting to discuss what happens with the map $\phi$. We have
seen that 
this map fails to be a homeomorphism when the Hessian of the area has an
null eigenvalue at a surface $S$. Thus, such 
surfaces correspond to folds of the map $\phi$. Note that
one typically gets such an unstable minimal surface between two stable ones,
so by a continuity argument there are also surfaces for which the Hessian
of the area has a null eigenvalue, so there 
are indeed folds of $\phi$. All this indicates that for non almost-Fuchsian
cases there is no uniqueness of solutions of (\ref{eq:**}), and the map $\phi$ is 
not invertible.

\section{GHMC AdS manifolds}

\subsection{Some background information on AdS}

Let $\R^4_2$ be $\R^4$ with the scalar product given by
$$ \langle x,x\rangle_{2,2} = x_1^2+x_2^2-x_3^2-x_4^2~. $$ 
The 3-dimensional 
Anti-de Sitter space --- called AdS here, and denoted by $AdS^3$ ---
can be defined as the domain in $\R P^3$ defined by the homogeneous equation: 
$\langle x,x\rangle_{2,2} <0$, 
with the scalar product coming under projectivization from
$\langle,\rangle_{2,2}$. 

It is a Lorentz manifold of dimension 3, with constant curvature $-1$. It is
not simply connected. Note that some definitions differ from the one used here
in that either a double cover or the universal cover is considered instead.

Another possible definition of $AdS^3$ is as one side of a quadric $Q$, of
signature $(1,1)$ in $\R P^3$ (all such quadrics are projectively equivalent),
with the ``Hilbert'' metric of $Q$, see e.g. \cite{shu} for this kind of
description. 

The space-like geodesics in $AdS^3$ are infinite geodesics, and the space-like
totally geodesic planes are isometric to the hyperbolic plane. The time-like
geodesics are closed curves of length $\pi$. The reader interested in the
elementary geometric properties of $AdS^3$ can find them e.g. in \cite{O}.
Here we will call $\nabla$ its Levi-Civit\`a connection.

\paragraph{The isometry group of $AdS^3$.}

It follows directly from the definition of $AdS^3$ given above that its
isometry group is $O(2,2)$. We are mostly interested here in the group of
orientation and time-orientation preserving isometries of $AdS^3$, denoted by
$\isom_+(AdS^3)$. It is interesting to note that it has a subgroup
of index two which is isomorphic to $PSL(2,\R)\times PSL(2,\R)$, see
\cite{mess}. Indeed, by the projective definition mentioned above,
$\isom(AdS^3)$ is the group of projective transformation leaving invariant a
quadric $Q$ of signature $(1,1)$ in $\R P^3$. But $Q$ is foliated by two
families of projective lines. Each of those two families can be identified
to $\R P^1$ (for instance by considering the intersections with one line of
the other family) and the isometries of $AdS^3$ which send lines of one family
to other lines of the same family act projectively on each of the two families
of lines. This defines a natural map from a subgroup of index 2 of
$\isom_+(AdS^3)$ --- the elements acting on each family of lines --- to
$PSL(2,\R)\times PSL(2, \R)$, and it is not difficult to show that it is an
isomorphism.  

\paragraph{GHMC AdS manifolds.} 

It is quite natural to define an {\it AdS manifold} as a manifold with a
Lorentz metric which is locally isometric to $AdS^3$. An important subclass
are those manifolds which are ``globally hyperbolic maximal compact'' (written
as ``GHMC'' below), which means that: 
\begin{itemize}
\item they contain a closed orientable space-like surface $S$,
\item each complete time-like geodesic intersects $S$ exactly once,
\item they can not be non-trivially isometrically embedded in any AdS manifold
  satisfying the first two properties.
\end{itemize}
Clearly, any GHMC AdS manifold is topologically the product of a closed
surface by an interval. 
The reader will find important results in \cite{mess}, which has served as an
important motivation for our work; actually one point of this section is to
recover some results of \cite{mess} by different methods. One result of
\cite{mess} is that, for GHMC AdS manifolds, the closed space-like surface has
genus at least 2, and that the holonomy is contained in the index 2 subgroup
of $O(2,2)$ isomorphic to $PSL(2,\R)\times PSL(2, \R)$. 

\begin{df}
As above we consider a closed surface $\Sigma$ of genus $g\geq 2$. We call 
$\cM_{g,AdS}$ the space of GHMC AdS
metrics on $\Sigma\times \R$, considered up to isotopy. 
\end{df}

\paragraph{Surfaces in AdS.}

The local differential geometry of surfaces in $H^3_1$ is mostly the same as
in hyperbolic space, with some minor differences which will be important later
on. Given a surface smooth $S\subset H^3_1$, it is said to be {\it space-like}
if the induced metric on $S$, still called $I$, 
is Riemannian (i.e. positive definite). Given
such a space-like surface and a unit normal vector field $N$ on $S$, we define
the shape operator $B:TS\rightarrow TS$ of $S$ as:
$$ Bx := - \nabla_xN~, $$
and its second fundamental form as $\II(x,y)=I(Bx,y)=I(x,By)$. Then $B$
satisfies the Codazzi equation, $d^\nabla B=0$, and a modified form of the
Gauss equation:
$$ K_g = -1 -\det(B) = -1 - det_g(h)~. $$
Conversely, an AdS form of the ``Fundamental Theorem of surface theory''
holds: on a simply connected surface $S$, given a metric $g$ and a symmetric
bilinear form $h$ which satisfy the Codazzi and the ``modified'' Gauss
equation, there is a unique immersion of $S$ in $AdS^3$ such that the induced
metric is $g$ and the second fundamental form is $h$.

The mean curvature of $S$ is still defined as the trace of $B$ divided by $2$.
Surfaces with zero mean curvature are called {\it maximal} surfaces, since the
second variation of the area is non-positive on the complement of a finite
dimensional space. 

\subsection{Maximal surfaces in germs of AdS manifolds}

\paragraph{Definition and first properties.}

The motivations leading to the definition of $\cH_g$ above apply also in the
AdS context. One is led to the following definition of a ``space of maximal
space-like surfaces in germs of AdS manifolds''.

\begin{df}
We call $\cH_{g,AdS}$ the space of couples $(g,h)$ such that $g$ is a smooth
metric on $\Sigma$ and $h$ is a symmetric bilinear form on $TS$, such that: 
\begin{itemize}
\item The trace of $h$ with respect to $g$, $\tr_g(h)$, vanishes.
\item $h$ satisfies the Codazzi equation with respect to the Levi-Civit\`a
  connection $\nabla$ of $g$: $d^\nabla h=0$.
\item The determinant of $h$ with respect to $g$ satisfies the corresponding 
Gauss equation: 
$K_g=-1-\det_g(h)$. 
\end{itemize}
\end{df}

An interesting difference with the hyperbolic case already appears: maximal
surfaces in germs of AdS manifolds always correspond to maximal surfaces in
GHMC AdS manifolds; moreover, this AdS manifold is unique. 
By contrast, minimal surfaces in hyperbolic manifolds do
not always correspond to minimal surfaces in quasi-Fuchsian
manifolds. 

\begin{lemma} \label{lm:g-ads}
Let $(g,h)\in \cH_{g,AdS}$. There exists a unique GHMC AdS manifold $M$
containing a maximal surface $S$ for which there exists a diffeomorphism
$u:\Sigma \rightarrow S$ such that the pull-back by $u$ of the first and
second fundamental forms of $S$ are equal to $g$ and $h$, respectively. 
\end{lemma}

\begin{proof}
By the ``fundamental theorem'' mentioned above, there exists a unique
immersion of the universal cover $\tilde{\Sigma}$ of $\Sigma$ in 
$AdS^3$ such that the
induced metric is $g$ and the second fundamental form is $h$. By the
uniqueness this immersion is equivariant under an action of $\pi_1\Sigma$. 
Taking the quotient by this action of a neighbourhood of the image of
$\tilde{\Sigma}$ in $AdS^3$ yields a non-complete AdS manifold $M'$ 
containing an embedded
maximal surface $S$ -- corresponding to quotient of the image of
$\tilde{\Sigma}$ -- with the required properties. 

By construction, $M'$ is globally hyperbolic and contains a compact space-like
surface, so it is contained in a GHMC manifold $M$, which is completely
determined by $g$ and $h$ and therefore unique.
\end{proof}

Another fundamental property of GHMC AdS manifolds is that they contain 
a unique space-like embedded maximal surface. This is an immediate
consequence of recent results \cite{BBZ-cras} on the existence of a 
(unique) foliation of those manifolds by constant mean curvature 
surfaces.

\paragraph{The second variation of the area.}

Maximal surfaces in AdS manifolds, as minimal surfaces in hyperbolic
 manifolds, are critical points of the area function. It is
therefore 
possible to consider the second variation of the area, which is given by a 
well-known integral formula very similar to the hyperbolic formula, although
one sign is different --- which will make things much simpler below.
Again we consider a first-order variation of a maximal surface given by an
orthogonal vector field of the form $fN$, where $N$ is a unit normal vector
field. 

\begin{lemma} \label{lm:second-ads}
The second variation of the area under the first-order deformation $fN$ is
given by: 
\begin{equation}\label{sec-var-ads}
A''(f) = - \frac{1}{2} \int_\Sigma \left( f\Delta f + 2(1-K_e) f^2 \right)
da~,  
\end{equation}
where the area form $da$ and the Laplace operator 
$\Delta f=-g^{ab}\nabla_a \nabla_b f$
come from the metric $g$ induced by
$u$ on $\Sigma$, and $K_e$ is the extrinsic curvature, i.e. the 
product of two principal curvatures of $\Sigma$.
\end{lemma}

A direct consequence is that maximal surfaces are always local maxima of the
area, because $K_e\leq 0$ since the mean curvature vanishes.

\paragraph{Maximal surfaces and HQD.}

Again as in the hyperbolic setting, there is a relation between points in the
cotangent space of Teichm\"uller space and maximal surfaces in germs of AdS
manifolds, but the AdS situation is much simpler. The next lemma mostly
repeats Lemma \ref{lm:hqd}, except in its last line where only one sign
changes.  

\begin{lemma} \label{lm:hqd-ads}
Let $g$ be a Riemannian metric on $\Sigma$, and let $h$ be a 
bilinear symmetric form on $T\Sigma$. Then:
\begin{enumerate}
\item The trace of $h$ with respect to $g$, $\tr_g(h)$, is zero if 
and only if $h$ is the real part of a quadratic differential $q$ over 
$\Sigma$.
\item If (1) holds, then $q$ is holomorphic if and only if $h$ satisfies the 
Codazzi equation, $d^\nabla h=0$.
\item If (1) and (2) hold, then $(g,h)$ is a maximal surface in a germ
of AdS manifold if and only if the Gauss equation is satisfied, 
i.e. $K_g=-1 - \det_g(h)$. 
\end{enumerate}
\end{lemma}

The proof is exactly the same as the proof of Lemma \ref{lm:hqd}, so we do not
repeat it here. However the behaviour of the equation in point (3) under
conformal changes of metric is much simpler than in section 2. 

\begin{lemma} \label{lm:conf-ads}
Let $(g,h)$ be as in Lemma \ref{lm:hqd-ads}, suppose that conditions (1) and 
(2) of that Lemma are satisfied. Let $g':=e^{2u}g$. Then conditions (1) 
and (2) are also satisfied for $(g',h)$, and condition (3) holds if and
only if:
\begin{equation} \label{eq:*-ads}
\Delta u = - e^{2u} - K_g - e^{-2u} det_g(h)~. 
\end{equation}
Moreover, if $h$ is not identically zero then 
this equation has a unique solution $u$.
\end{lemma}

The proof of the first point is the same as the proof of Lemma
\ref{lm:confchange}. Equation (\ref{eq:*-ads}) is known to have a unique
solution when $k$ is identically zero, because in this case it reduces to the
better-known prescribed metric equation, see \cite{troyanov}. So from here on
we suppose that $k$ is not identically zero, and, since $h$ is the real part
of a holomorphic quadratic differential, it follows that $k$ has only isolated
zeros. 

To solve equation (\ref{eq:*-ads}), consider the following functional:
$$ 
\begin{array}{lrcl}
F: & L^2(\Sigma) & \rightarrow & \R \\
   & u & \mapsto & 
\frac{1}{2}\int_\Sigma \left( \| \nabla u\|^2 + e^{2u} + 2K_gu + k^2 e^{-2u}
   \right) da~,    
\end{array} $$
where $-k^2=det_g(h)$  and the norm of $\nabla u$ and the area form are
relative to 
$g$. A simple computation shows that: 
\begin{eqnarray*}
\delta F(u) =  
\int_\Sigma \delta u (\Delta u + e^{2u} + K_g - k^2 e^{-2u}) da~,
\end{eqnarray*}
so $u$ is a solution of (\ref{eq:*-ads}) if and only if it is a critical point
of $F$. 

Moreover, $F$ is the sum of four convex functionals -- the last term is convex
because $\det_g(h)\leq 0$ since $\tr_g(h)=0$ -- so it is convex, and
actually strictly convex since several of the four summands are strictly
convex. So $F$ has at most one critical point, which can only be a minimum.

The proof that equation (\ref{eq:**}) actually has a solution is based on the
following technical statement.

\begin{prop} \label{pr:ineg}
Under the hypothesis of Lemma \ref{lm:hqd-ads}, if $K\in (-1,0)$ on $\Sigma$
and if $k$ is not identically zero, there exists a constant
$\epsilon_0>0$ such that, for all $u\in L^2(\Sigma)$, $F(u)\geq
\epsilon_0 \| u\|^2_{L^2}$. 
\end{prop}

\begin{proof}
Consider the function: 
$$ f_{K,k}(u) := e^{2u} +2Ku +k^2e^{-2u}~. $$
Since $K\in (-1,0)$, $e^{2u}+2Ku\geq 0$ for all $u\in \R$, so that $f_{K,k}$
is positive on $\R$. Moreover, at all points where $k\neq 0$,
$f_{K,k}(u)/u^2\rightarrow \infty$ as
$u\rightarrow \pm \infty$. So there exist a non-empty open subset 
$\Omega\subset \Sigma$ (which is the complement of the neighbourhood of the
points where $k$ vanishes) and a constant $\epsilon_1>0$ such that:
$$ \forall u\in L^2(\Sigma), 
\forall x\in \Omega, e^{2u(x)} +2K(x)u(x) +k(x)^2e^{-2u(x)} \geq \epsilon_1
u(x)^2~. $$
It follows that:
$$ \forall u\in L^2(\Sigma), F(u) \geq \epsilon_1 \int_\Omega u^2da + 
\int_\Sigma \| du\|^2da~, $$
so the proposition will follow if we can prove that there exists 
$\epsilon_2>0$ such that:
$$ \forall u\in L^2(\Sigma),\int_\Omega u^2da + 
\int_\Sigma \| du\|^2da \geq \epsilon_2 \int_{\Omega'} u^2 da~. $$
Let $\Omega':=\Sigma\setminus \Omega$, let $A$ and $A'$ be the areas of
$\Omega$ and $\Omega'$, and let:
$$ u_0 = \frac{1}{A+A'} \int_\Sigma uda~. $$
By the Poincar\'e inequality there exists a constant $\epsilon_3>0$
such that:
$$ \forall u\in L^2(\Sigma), \epsilon_3 \int_\Sigma (u-u_0)^2 da \leq 
\int_\Sigma \| du\|^2 da~, $$
so the proposition will be established if we can show that there
exists a constant $\epsilon_4>0$ such that: 
\be \label{eq:target}
\forall u\in L^2(\Sigma),\int_\Omega u^2da + 
\int_\Omega (u-u_0)^2da + \int_{\Omega'} 
(u-u_0)^2da \geq \epsilon_4 \int_{\Omega'} u^2 da~.
\ee

At this point it is helpful to simplify slightly the notations by 
setting (without loss of generality):
$$ \int_{\Omega'} u^2 da = 1, \epsilon := \int_{\Omega'} (u-u_0)^2 da~.$$
Then:
$$ \epsilon := \int_{\Omega'} (u-u_0)^2 da = \int_{\Omega'} u^2 da 
- 2u_0\int_{\Omega'} u da + A'u_0^2~, $$
and, since $\left| \int_{\Omega'}uda \right|\leq \sqrt{A'}$ by the
Cauchy-Schwarz inequality:
\begin{equation}\label{eq:add} 
\epsilon \geq 1 - 2u_0\sqrt{A'} + A'u_0^2 = (1-u_0\sqrt{A'})^2~. 
\end{equation}
This means that, if $\epsilon\leq 1/4$, then $u_0\geq 1/2\sqrt{A'}$. 
However, the same computations as the one we have just done, with $\Omega'$ 
replaced by $\Omega$, shows that: 
\begin{eqnarray*}
\int_\Omega (u-u_0)^2 da & = & \int_{\Omega} u^2 da 
- 2u_0\int_{\Omega} u da + A u_0^2 \\
& \geq & \int_{\Omega} u^2 da 
- 2u_0\sqrt{A\int_\Omega u^2da} + A u_0^2 \\
& \geq & \left( \sqrt{A} u_0 - \sqrt{\int_\Omega u^2 da}\right)^2~, 
\end{eqnarray*}
and this shows that there exists $\epsilon_5>0$ such that, if 
$u_0\geq 1/2\sqrt{A'}$, then:
\begin{itemize}
\item either $\int_\Omega u^2 da\geq \epsilon_5$,
\item or $\int_\Omega (u-u_0)^2 da\geq \epsilon_5$.
\end{itemize}
However, if $u_0<1/2\sqrt{A'}$, then it follows from equation (\ref{eq:add})
that: 
$$ \epsilon = \int_{\Omega'}(u-u_0)^2 da > \frac{1}{4}~. $$
Summing up, this shows that there exists some $\epsilon_6>0$
such that one of the terms on the left-hand side of equation
(\ref{eq:target}) is larger than $\epsilon_6$, and this proves
the proposition.
\end{proof}

\begin{proof}[Proof of Lemma \ref{lm:conf-ads}]
The lemma now follows from standard tools of functional analysis.
Starting from a couple $(g,h)$, it is possible to make a first
conformal transformation so as to obtain another couple, $(g',h)$
for which $K\in (-1,0)$, so that Proposition \ref{pr:ineg} applies.

One can then consider a sequence $(u_n)_{n\in \N}$ which is minimizing
the functional $F$, and the previous proposition shows that $(u_n)$
remains in a ball in $L^2(\Sigma)$. Therefore $(u_n)$ is
weakly converging to a limit $u_\infty$, and the usual arguments
involving elliptic regularity then show that $u_\infty$ is actually
a smooth function.
\end{proof}

\paragraph{GHMC AdS manifolds from maximal surfaces.}

The previous lemmas on maximal surfaces in AdS manifolds can be summed up in a
statement which is analogous to, but simpler than, Theorem \ref{tm:hyperbolic}.

\begin{thm} \label{tm:ads}
There exists a natural homeomorphism $\psi_{g,AdS}:M_{g,AdS}\rightarrow
T^*\cT_g$. Given a GHMC AdS metric $G$, it contains a unique, embedded,
space-like maximal
surface $S$, with induced metric $g$ and second fundamental form $h$. Then $h$
is the real part of a QHD $q$ on $S$. Thus, one obtains $\psi(G)$ that is the
element of 
$T^*\cT_g$ that is associated to $(c,q)$; here $c$ is the complex structure of
$g$. 
Conversely, for all $(c,q)\in T^*\cT_g$, we call $h$ the real part of $q$, and
there is then a unique metric $g$ in the conformal class defined by $c$ such
that $(g,h)\in \cH_{g,AdS}$, and then $g=I$ and $h=\II$ for a (unique) maximal
surface in a unique GHMC AdS metric $G$ on $\Sigma \times \R$.
\end{thm}

\begin{proof}
The existence of a unique embedded, space-like maximal surface in a GHMC AdS
manifold is a recent result of \cite{BBZ-cras} (see also
\cite{barbot-zeghib}). In the converse statement, the
existence and uniqueness of $g$ is established by Lemma \ref{lm:hqd-ads} and
Lemma \ref{lm:conf-ads}.  
\end{proof}

\paragraph{CMC surfaces.}

As in hyperbolic manifolds, the results given above on maximal surfaces in AdS
manifolds extend to CMC surfaces, with little differences. There is a direct
generalization of Lemma \ref{lm:hqd-ads}:

\begin{lemma} \label{lm:hqd-ads-cmc}
Let $g$ be a Riemannian metric on $\Sigma$, and let $h$ be a 
bilinear symmetric form on $T\Sigma$. Let $h_0=h-Hg$. Then:
\begin{enumerate}
\item The trace of $h_0$ with respect to $g$, $\tr_g(h)$, is zero if 
and only if $h_0$ is the real part of a quadratic differential $q$ over 
$\Sigma$.
\item If (1) holds, then $q$ is holomorphic if and only if $h$ satisfies the 
Codazzi equation, $d^\nabla h=0$.
\item If (1) and (2) hold, then $(g,h)$ is a CMC-$H$ surface in a germ
of AdS manifold if and only if the Gauss equation is satisfied, 
i.e. $K_g=(H^2-1) - \det_g(h_0)$. 
\end{enumerate}
\end{lemma}

The proof follows exactly the same arguments as the proof of Lemma
\ref{lm:hqd-ads}, we leave the details to the interested reader. 

Moreover, there is also a direct extension of Lemma \ref{lm:conf-ads}, whose
proof is again the same as in the maximal surface case.

\begin{lemma} \label{lm:conf-ads-cmc}
Let $(g,h)$ be as in Lemma \ref{lm:hqd-ads}, suppose that conditions (1) and 
(2) of that Lemma are satisfied. Let $g':=e^{2u}g$. Then conditions (1) 
and (2) are also satisfied for $(g',h)$, and condition (3) holds if and
only if:
\begin{equation} \label{eq:*ads-cmc}
\Delta u = (H^2-1) e^{2u} - K_g - e^{-2u} det_g(h_0)~. 
\end{equation}
Moreover, if $h_0$ is not identically zero and if $|H|<1$ then 
this equation has a unique solution $u$.
\end{lemma}

It follows that there is a natural description of the space of GHMC AdS
manifolds in terms of $T^*\cT_g$ analog to what is described in Theorem
\ref{tm:ads}, but using CMC-$H$ surfaces -- for $|H|<1$ -- rather than maximal
surfaces. This will not be used here in the AdS context, but will be used in
section 6 to study Minkowski manifolds, for which maximal surfaces are not
convenient. 

\paragraph{All GHMC AdS manifolds are ``almost-Fuchsian''.}

An elementary but interesting remark is that all GHMC AdS manifolds have the
property used in section 2 when studying hyperbolic manifolds: the maximal
surface that they contain has principal curvatures of absolute value less than
1. This can be considered as another instance of -- or an explanation of --
the relative simplicity of the AdS case, as compared to its hyperbolic
analog. It will also be important below to understand the Mess parameterization
of $\cM_{g,AdS}$ by two copies of Teichmüller space. 

\begin{lemma} \label{lm:curvatures-ads}
Let $S$ be a closed maximal surface in an AdS manifold. The principal
curvatures of $S$ have absolute value less than $1$.
\end{lemma}

\begin{proof}
The result is obvious if $S$ is totally geodesic, so we now suppose that this
is not the case.
Let $B$ be the shape operator of $S$, so that $B$ has eigenvalues $k$ and
$-k$, for $k\geq 0$, and $k$ has only isolated zeros. Let $(e,e')$ be the
orthonormal frame on $S$ (outside the points where $k=0$) with $e$ and $e'$
eigenvectors of $B$ of eigenvalues $k$ and $-k$, respectively. 

The Codazzi equation, $d^\nabla B=0$, translates as:
$$ \left\{ 
\begin{array}{ccc}
e^a\nabla_a k & = & -2k \omega(e') \\
e'{}^a\nabla_a k & = & 2k \omega(e)
\end{array} \right. $$
where $\omega$ is the connection 1-form of $(e,e')$, that is:
$\nabla_x e = \omega(x) e', \nabla_x e' = -\omega(x) e$. 
Setting $\chi:=\log(k)/2$, this can be written as:
$$ \left\{ \begin{array}{ccc}
e.\chi = e^a \nabla_a \chi & = & -\omega(e') \\
e'.\chi =e'{}^a\nabla_a\chi & = & \omega(e)
\end{array} \right. $$
It follows that: 
\begin{eqnarray*}
d\omega(e,e') & = & e.\omega(e') - e'.\omega(e) - \omega([e,e']) \\
& = & -e.e.\chi -
e'.e'.\chi +\omega(\omega(e)e +\omega(e')e') \\
& = & - e^a\nabla_a e^b \nabla_b\chi -
e'{}^a\nabla_a e'{}^b\nabla_b \chi +(e'{}^a\nabla_a \chi)^2 + (e^a\nabla_a \chi)^2 \\
& = & - (e^ae^b+e'{}^ae'{}^b)\nabla_a\nabla_b\chi = \Delta\chi~.
\end{eqnarray*}
But $d\omega(e,e') = -K$ by definition of the curvature, so that: 
$$ \Delta\chi = -K = 1-k^2 = 1-e^{4\chi}~. $$
Let us introduce: $f(\chi)=1-e^{4\chi}$. Then, $f(\chi)<0$ 
for all positive $\chi$. On the
other hand, 
at points where $\chi$ attains its maximum $\Delta\chi\geq 0$.
Thus, by the maximum principle, $\chi$ must be less then zero where it reaches
its maximum, and strictly less than $0$ unless it is identically zero.
Therefore it is less than zero everywhere on $S$, so $k<1$ everywhere on $S$.
\end{proof}

\subsection{Canonical maps between hyperbolic surfaces}

\paragraph{Bundle morphisms and changes of metric.}

The material of this paragraph deals with changes of metrics associated to
bundle morphisms on the tangent space of a surface. It will be
frequently used in what follows. The following proposition is
well-known, see \cite{L5}.

\begin{prop} \label{pr:metric-change}
Let $S$ be a surface, with a Riemannian metric $g$. Let $A:TS\rightarrow TS$ be
a smooth bundle morphism such that $A$ is everywhere invertible and 
$d^\nabla A=0$, where $\nabla$ is the
Levi-Civit\`a connection of $g$. Let $h$ be defined by:
$$ h(u,v) = g(Au,Av)~. $$
Then the Levi-Civit\`a connection of $h$ is given by:
$$ \nabla^h_uv = A^{-1}\nabla_u(Av)~,$$
and its curvature is 
$$ K^h = \frac{K^g}{\det(A)}~. $$
\end{prop}

\begin{proof}
Consider the connection $\nabla'$ on $S$ defined by: 
$\nabla'_uv = A^{-1}\nabla_u(Av)$.
It is torsion-free, because, for all vector fields $u,v$ on $S$:
\begin{eqnarray*}
\nabla'_uv - \nabla'_vu & = & A^{-1}\nabla_u(Av) - A^{-1}\nabla_v(Au) \\
& = & A^{-1}((\nabla_uA)v + A\nabla_uv - (\nabla_vA)u - A\nabla_vu) \\
& = & A^{-1} (d^\nabla A)(u,v) + \nabla_uv - \nabla_vu \\
& = & 0~.
\end{eqnarray*}
Moreover, $\nabla'$ is compatible with $h$, because, 
for all vector fields $u,v,w$ on $S$:
\begin{eqnarray*}
u.h(v,w) & = & u.g(Av,Aw) \\
& = & g(\nabla_u(Av),Aw) + g(Av,\nabla_u(Aw)) \\
& = & h(A^{-1}\nabla_u(Av),w) + h(v, A^{-1}\nabla_u(Aw)) \\
& = & h(\nabla'_uv,w) + h(v, \nabla'_uw)~.
\end{eqnarray*}
So $\nabla'$ is the Levi-Civit\`a connection of $h$.

Let $(e_1, e_2)$ be an orthonormal moving frame on $S$ for $g$, and let
$\beta$ be its connection 1-form, i.e.:
$$ \nabla_xe_1 = \beta(x)e_2, ~\nabla_xe_2 = -\beta(x)e_1~. $$
Then the curvature of $g$ is defined as: $d\beta = -K da$. 
Now let $(e'_1,e'_2) := (A^{-1}e_1,A^{-1}e_2)$; clearly it
is an orthonormal moving frame for $h$. Moreover the expression of
$\nabla'$ above shows that its connection 1-form is also $\beta$.
It follows that:
$Kda = -d\beta = K'da'$, where $K'$ is the curvature
and $da'$ is the area form of $h$, so that:
$$ K' = K \frac{da}{da'} = \frac{K}{\det(A)}~. $$ 
\end{proof}

\paragraph{Bundle morphisms and canonical diffeomorphisms.}

The changes of metrics described above are related to canonical
diffeomorphisms between hyperbolic surfaces (or more generally constant
curvature) metrics.

\begin{thm}[F. Labourie \cite{L5}] \label{tm:labourie}
Let $S$ be a closed surface of genus at least $2$, and let $g_+, g_-$ 
be two hyperbolic metrics on $S$. There exists a unique bundle morphism 
$b:TS\rightarrow TS$ such that:
\begin{itemize}
\item $b$ is self-adjoint for $g_+$, with positive eigenvalues. 
\item $d^{\nabla^+} b=0$, where $\nabla^+$ is the Levi-Civit\`a connection of
  $g_+$. 
\item $\det(b)=1$.
\item $g_-$ is the pull-back of $g_+(b\cdot, b\cdot)$ by a diffeomorphism
isotopic to the identity.
\end{itemize}
\end{thm}

The diffeomorphism obtained in this manner can be described geometrically (see
\cite{L5}): it is characterized by the fact that its graph is minimal
in the product of the two hyperbolic surfaces.

\subsection{GHMC AdS manifolds from space-like surfaces}

\paragraph{The Mess parameterization of $\cM_{g,AdS}$.}

G. Mess \cite{mess} has discovered a nice parameterization of $\cM_{g,AdS}$ by
the product of two copies of the Teichm\"uller space $\cT_g$. This map can be
described briefly as follows. Recall -- see \cite{mess} -- that $AdS^3$ has a
projective map to the interior of a solid 
torus $T^3$  in $\R P^3$, and that its
``boundary at infinity'' can be identified with the boundary of $T^3$, which
is a torus $T_\infty$. $T_\infty$ is foliated by two families of projective
lines, which we call $\cL_+$ and $\cL_-$ here.

Choose a fixed space-like totally geodesic plane
$P_0\subset AdS^3$, so that $P_0$, with the induced metric, is isometric to 
the hyperbolic plane. The boundary at infinity of $P_0$ is a closed curve in
$T_\infty$, which intersects exactly once each line in $\cL_-$ and each line
in $\cL_-$.

Let $P$ be any other space-like plane in $AdS^3$, so that the description made
for $P_0$ also applies to $P$. There are two natural maps $\pi_{\infty,+}$
and $\pi_{\infty, -}$ from $\dr_\infty P$ to $\dr_\infty P_0$, sending a point
$x\in \dr_\infty P$ to the intersection with $\dr_\infty P_0$ of the line of
$\cL_+$ (resp. $\cL_-$) containing $x$. It is not difficult to check that
$\pi_{\infty,+}$ and $\pi_{\infty, -}$ are projective, so that they are the
extension at infinity of hyperbolic isometries $\pi_+, \pi_-: P\rightarrow
P_0$. 

This defines two maps $\Pi_+, \Pi_-$ from the space of unit time-like
vectors tangent to $AdS^3$, say $U^1AdS^3$, 
to $P_0$: if $x\in AdS^3$ and $v\in T_xAdS^3$ is a
unit time-like vector, $v$ is orthogonal to a unique space-like plane $P$
containing $x$, and $\Pi_\pm(x)$ is defined as $\pi_\pm(x)$, where
$\pi_\pm$ is defined by reference to this plane $P$.

The parameterization used by Mess also uses two special, non-smooth surfaces,
which are the two boundary components of the convex core of a GHMC AdS
manifold, which we can call $S_+$ and $S_-$. Each of those surfaces ``lifts''
to a surface $S'_+$ (resp. $S'_-$) in $U^1AdS^3$, namely the set of unit
oriented normals to $S_+$ (resp. $S_-$). The pull-backs by $\Pi_+$
and by $\Pi_-$ of the hyperbolic metric on $P_0$ defines a hyperbolic metric
on $S'_+$, therefore two points in the Teichm\"uller space $\cT_g$. Those
hyperbolic metrics have a nice geometric interpretation: they are obtained
from the hyperbolic metric induced on $S_+$ by a right (resp. left) earthquake
along the measured bending lamination of $S_+$.

Using $S'_-$ instead of $S'_+$ leads to the same couple of points in
$\cT_g$. It is proved in  
\cite{mess} that this couple of points in $\cT_g$ uniquely determines the AdS
metric $G$, and that every point in $\cT_g\times \cT_g$ can be obtained.

\paragraph{Hyperbolic metrics on space-like surfaces.}

There is a differential-geometric version of the construction which has just
been briefly recalled, replacing the boundary components of the convex core by
any smooth surface.

\begin{df}
Let $S$ be a (smooth) embedded space-like surface in a GHMC AdS manifold $M$.
Let $I$ and $B$ be the induced metric and shape operator of $S$, and let $J$
be the complex structure on $S$ defined by $I$. We define symmetric bilinear
forms $I^\#_+$ and $I^\#_-$ on $S$ as:
$$ I^\#_\pm (x,y) := I((E\pm JB)x, (E\pm JB)y)~, $$
where $E$ is the identity.
\end{df}

It is quite clear that $I^\#_+$ and $I^\#_-$ are smooth metrics as soon 
as the eigenvalues of $B$ are in $(-1,1)$ -- however this is not necessary. 
By Lemma \ref{lm:curvatures-ads}, 
this holds in particular if $S$  is the maximal surface in $M$. Moreover, 
$I^\#_+$ and $I^\#_-$ are always hyperbolic metrics. 

\begin{lemma} \label{lm:hyper-metrics}
$I^\#_+$ and $I^\#_-$ have constant curvature equal to $-1$.
\end{lemma}

The proof is mostly a consequence of Proposition \ref{pr:metric-change}.

\begin{proof}
First note that $E+JB$ is a solution of the Codazzi equation, because, for all
$x\in S$ and all $u,v\in T_xS$:
\begin{eqnarray*}
d^\nabla(E+JB)(x,y) & = & (\nabla_x(E+JB))y - (\nabla_y(E+JB))x \\
& = & J(\nabla_xB)y - J(\nabla_yB)x \\
& = & J(d^\nabla B)(x,y) \\
& = & 0~. 
\end{eqnarray*}
Thus it follows from Proposition \ref{pr:metric-change}
that the curvature $K^+$ of $I^\#_+$ is equal to:
$$ K^+ = \frac{K}{\det(E+JB)} = \frac{-1-\det(B)}{1+\tr(JB)+\det(JB)} 
= \frac{-1-\det(B)}{1+\det(B)} = -1~. $$
The same arguments can be used for $I^\#_-$.
\end{proof}

\paragraph{The Mess parameterization from embedded surfaces.}

It should not come as a surprise to the reader that the hyperbolic metrics
$I^\#_+$ and $I^\#_-$ are the same as the hyperbolic metrics defined by Mess.

\begin{lemma}
Let $G$ be a GHMC AdS metric on $\Sigma\times \R$, and let $S$ be a closed 
embedded space-like surface in $(\Sigma\times \R, G)$. Then the metrics
$I^\#_+$ and $I^\#_-$ defined from $S$ are the same as the two hyperbolic
metrics associated to $G$  by the construction of Mess recalled above. In
particular, they do not depend on the choice of $S$.
\end{lemma}

\begin{proof}
Let $S$ be a space-like surface in $(\Sigma\times \R, G)$. We first define
two hyperbolic metrics $g_+$ and $g_-$ on $S$. 
Let $\St\subset AdS^3$ be the universal cover of $S$, and let $\rho$ be the 
canonical projection to $\St$ from its lift $\St'\subset U^1AdS^3$.
Then let $\phi_\pm:=\Pi_\pm\circ \rho^{-1}$, and let $g_\pm$ be 
the pull-backs to $\St$, by the maps $\phi_\pm$, of the
hyperbolic metric on $P_0$. The definition then
shows that $g_+$ and $g_-$ are invariant under the action of the fundamental
group of $S$ --- because replacing $P_0$ by another plane does not change
the pull-back metrics obtained --- so $g_\pm$ define hyperbolic metrics on $S$.

Given $x\in \St$, let $P$ be the totally geodesic plane tangent to $\St$ at $x$.
Then $T_x\St=T_xP$, and the differential $d_x\Pi_+:T_xP=T_x\St\rightarrow 
T_{\phi_+(x)}P_0$ can be used to identify $T_x\St$ with  
$T_{\phi_+(x)}P_0$. It is then interesting to compute the differential $d\phi_+(v)$,
for $v\in T_x\St$, it can be decomposed as follows: 
\begin{itemize}
\item one term corresponds to the parallel displacement of the unit normal $N$ to $\St$
at $x$, it is simply $d_x\Pi_+(v)$.
\item the other corresponds to the first-order variation of $N$ 
under a displacement on $\St$, corresponding to the fact that the
first-order variation of $N$ is equal to $Bv$. To understand it note
that the first-order variation of $P$ is a ``rotation'' along the geodesic $\Delta$ in 
$P$ orthogonal to $Bv$, induced by an AdS Killing vector
field which vanishes on the $\Delta$. The corresponding variation on $P_0$ 
vanishes on the image of $\Delta$, and it is a Killing vector field of $P_0$ since
$\pi_+:P\rightarrow P_0$ is always an isometry. So this Killing field is an 
infinitesimal translation of axis $\Delta$. Its translation speed can be computed
by taking for instance $P_0$ as a plane intersecting $P$ along $\Delta$ (since changing
$P_0$ is not going to change the resulting map $\phi_+$ up to isometries), a direct
computation then shows that it is equal to the norm of $v$. 
\end{itemize}
Adding those two terms, we find that, using the identification above: 
$$ d_x\phi_+(v) = v + JBv~, $$
so $g_+=I^\#_+$, as claimed. The same arguments also shows that $g_-=I^\#_-$. 

Now let $\gamma\in \pi_1\Sigma$, and let $\gamma'\in \mbox{Isom}_+(AdS^3)$ 
be the image of $\gamma$ by the holonomy representation of the AdS metric $G$.
Let $x\in \St$, so that $\gamma'(x)\in \St$ since $\St$ is invariant under the
action of $\pi_1\Sigma$ on $AdS^3$. The image of $\gamma'x$ by $\phi_+$ can
be obtained simply from the definition of $\phi_+$: it is the image of 
$\gamma' \phi_+(x)$ under the isometry $\Pi_{\gamma'}$ from $\gamma'P_0$ to $P_0$ which 
extends to $\dr_\infty (\gamma'P_0)$ as the projection on $\dr_\infty P_0$
along the lines in $\cL_+$. In other terms:
$$ \phi_+\circ \gamma' = \Pi_{\gamma'}\circ \gamma'\circ \phi_+~, $$
where $\Pi_{\gamma'}\circ\gamma'$ is considered as an isometry from $P_0$
to itself. This means that the holonomy of the hyperbolic metric $g_+$
under the action of $\gamma$ is given by $\Pi_{\gamma'}\circ\gamma'$, which
depends only on $\gamma'$ and not on the choice of the space-like surface $S$.
Therefore, $g_+$, as a hyperbolic metric on $\Sigma$ defined up to isotopy, 
does not depend on the choice of $S$. The same argument shows the same result
for $g_-$. (Another, more differential-geometric proof of the fact that
$I^\#_\pm$ do not depend on the 
choice of $S$ is given, in a more general context, in section 5.)

Finally, to prove that $g_+$ and $g_-$ correspond to the metrics defined
by Mess, it is sufficient to take the limit as $S$ goes to the upper (resp.
lower) boundary of the convex core of $(\Sigma\times \R, G)$.
\end{proof}

Therefore, we get a map:
\begin{equation}\label{mess}
\mess: \cM_{g,AdS} \to \cT_g\times \cT_g,
\end{equation}
sending a GHMC AdS metric on $\Sigma\times \R$ to the two hyperbolic metrics
on $\Sigma$ obtained as $I^\#_-$ and $I^\#_+$ for any choice of a space-like
surface in it. 

G. Mess \cite{mess} has proved that the map $\mess$ is one-to-one, thus
providing a 
parameterization of $\cM_{g,AdS}$ by $\cT_g\times \cT_g$. We provide here
another proof of this fact, using the canonical maps described above 
between hyperbolic surfaces,
and maximal surfaces in AdS manifolds. We will also recover as a by-product
the fact that any GHMC AdS manifold contains a unique maximal surface, by
methods very different from those used before. This approach
will be extended in section 5 to singular GHMC AdS manifolds.

\begin{thm} \label{tm:mess}
The map $\mess$ is one-to-one. Moreover, each GHMC AdS manifold contains a
unique embedded, space-like maximal surface.
\end{thm}

\begin{proof}
To prove the first point, we have to prove that, given $g_+,g_-\in \cT_g$,
there is a GHMC AdS metric $G$ on $\Sigma\times \R$ and a space-like surface
$S\subset \Sigma\times \R$ for that metric such that $I^\#_+=g_+$ and that
$I^\#_-=g_-$. We will actually show the existence of a maximal surface
$S\subset \cH_{g,AdS}$ such that $I^\#_+=g_+$ and that
$I^\#_-=g_-$, which will prove the second point at the same time.

By Theorem \ref{tm:labourie}, there is a unique bundle morphism $b:T\Sigma
\rightarrow T\Sigma$ such that:
\begin{itemize}
\item $b$ is self-adjoint for $g_+$, with positive eigenvalues.
\item $d^{\nabla^+}b=0$, where $\nabla^+$ is the Levi-Civit\`a connection of
  $g_+$.
\item $\det(b)=1$.
\item $g_-$ is isotopic to $g_+(b\cdot, b\cdot)$.  
\end{itemize}
Define a metric $g$ on $\Sigma$ by:
$$ 4g := g_+((E+b)\cdot, (E+b)\cdot)~, $$
and let $J$ be the complex structure on $\Sigma$ defined by $g$. 
Define another bundle morphism $B:T\Sigma\rightarrow T\Sigma$ by:
$$ JB:= (E+b)^{-1}(E-b)~. $$
Then $B$ is well-defined at all points since $b$ has positive eigenvalues. 
It follows from a simple computation that,
since $\det(b)=1$, $\tr(JB)=0$, so that $B$ is self-adjoint for $g$. 
Since $b$ has positive eigenvalues, the 
eigenvalues of $B$ are in $(-1,1)$, and another simple computation shows that: 
$$ b = (E+JB)^{-1}(E-JB)~, ~ ~ E+JB = 2(E+b)^{-1}~. $$
Note also that $JB$ is self-adjoint for $g$, indeed if $x,y$ are vectors
tangent to $\Sigma$ at the same point, then:
\begin{eqnarray*}
4g(JBx,y) & = & g_+((E+b)(E+b)^{-1}(E-b)x, (E+b)y) \\
& = & g_+((E-b)x, (E+b)y) \\
& = & g_+((E+b)x, (E-b)y) \\
& = & 4g(x, JBy)~,
\end{eqnarray*}
and it follows that $\tr(B)=0$. 

By Proposition \ref{pr:metric-change}, the Levi-Civit\`a connection $\nabla$
of $g$ is:
$$ \nabla_xy = (E+b)^{-1}\nabla^+((E+b)y)~. $$
Therefore: 
\begin{eqnarray*}
(d^\nabla JB)(x,y) & = & (E+b)^{-1}\nabla^+_x((E+b)JBy) -
(E+b)^{-1}\nabla^+_y((E+b) JBx) - JB[x,y] \\
& = & (E+b)^{-1}\nabla^+_x((E-b)y - (E+b)^{-1}\nabla^+_y((E-b)y) -
(E+b)^{-1}(E-b)[x,y] \\ 
& = & (E+b)^{-1}(d^{\nabla^+}(E+b))(x,y) \\
& = & 0~. 
\end{eqnarray*}
Moreover, $d^{\nabla}(JB)=Jd^\nabla B$, 
because $J$ is the complex structure of $g$, and it
follows that $d^\nabla B=0$. 

To check that $B$ verifies the Gauss equation, note that, by Proposition 
\ref{pr:metric-change} and since $g_+$ has curvature $-1$, 
its curvature is equal to: 
$$ K_g = \frac{-1}{\det((E+b)/2)} = - \det(E+JB) = -1-\det(B)~. $$

We have seen that $B$ is traceless, self-adjoint for $g$, and that it
satisfies the Codazzi and Gauss 
equation. So, setting $h:=g(B\cdot, \cdot)$, $(g,h)\in
\cH_{g,AdS}$. It is also quite clear that, if $x$ and $y$ are vectors tangent 
to $\Sigma$ at the same point, then:
\begin{eqnarray*}
I^\#_+(x,y) & = & g((E+JB)x,(E+JB)y) \\
& = & g(2(E+b)^{-1}x, 2(E+b)^{-1}y) \\
& = & g_+(x,y)~, 
\end{eqnarray*}
while:
\begin{eqnarray*}
I^\#_-(x,y) & = & g((E-JB)x,(E-JB)y) \\
& = & g(2(E+b)^{-1}bx, 2(E+b)^{-1}by) \\
& = & g_+(bx,by) \\
& = & g_-(x,y)~. 
\end{eqnarray*}
This shows that the map $\mess$ is one-to-one.

This also shows that any GHMC AdS manifold contains a maximal surface,
obtained by the construction just outlined from the two hyperbolic
metrics $I^\#_\pm$ induced on any space-like surface. We still have to
prove that this maximal surface is unique. But note that, given such a
maximal surface $S$, with its induced metric $g$ and shape operator $B$,
we can define: 
$$ b := (E-JB)(E+JB)^{-1}~. $$
Then $\det(b)=1$ because $\tr(JB)=0$. A direct check (as done above for $B$)
shows that $b$ is self-adjoint for $I^\#_+$, because $(E-JB)$ is self-adjoint
for $g$ (this follows from the fact that $B$ is traceless).  
Proposition \ref{pr:metric-change} shows that the Levi-Civit\`a connection
of $I^\#_+$ is: $\nabla^+_xy = (E+JB)^{-1}\nabla_x((E+JB)y)$, and another
direct computation shows then that $d^{\nabla^+}b=0$. So $b$ is the unique 
bundle morphism associated by Theorem \ref{tm:labourie} to $I^\#_+$ and
$I^\#_-$. Since $JB=(E+b)^{-1}(E-b)$ and $g$ is also uniquely determined by $g_+$
and $b$, there is at most one maximal surface with a given pair of hyperbolic
metrics as $I^\#_+$ and $I^\#_-$, which proves the theorem.
\end{proof}

\subsection{Another parameterization by two copies of Teichm\"uller space}

\paragraph{Definitions.}

In addition to the map discovered by Mess, there is another natural map
from $\cT_g\times \cT_g$ to $\cM_{g,AdS}$, defined in terms of the maximal 
surface in a GHMC AdS manifold and of two hyperbolic metrics defined on it.

\begin{df}
Let $G\in \cM_{g,AdS}$ be a GHMC AdS metric on $\Sigma\times \R$, and let 
$S$ be the space-like embedded maximal surface in $(\Sigma\times \R, G)$. We
will consider the metrics $I^*_\pm$ on $S$ defined by:
\begin{equation}\label{I*}
I^*_\pm(u,v) = I(u,v) \pm 2\II(u,v) + \III(u,v) = I((E\pm B)u, (E\pm B)v)~.
\end{equation}
\end{df}
We note that both of these metrics are hyperbolic. This is proved using
Proposition \ref{pr:metric-change}. We also note 
that the metrics $I^*_\pm$ are non-singular because, by Lemma 
\ref{lm:curvatures-ads}, the principal curvatures of $B$ are everywhere in
$(-1,1)$. 
Therefore, we get a map:
\begin{equation}\label{max}
\max: \cM_{g,AdS} \to \cT_g\times \cT_g.
\end{equation}
As is indicated in the name, this map is different from the one given by
(\ref{mess}). 

\paragraph{Canonical maps.}

We would now like to show that it is possible to construct the inverse of the
map (\ref{max}) in a very explicit fashion. To this end, we first note that
there is a simple relation between the two metrics $I^*_\pm$.
This follows from Theorem \ref{tm:labourie}, which can be 
applied to $I^*_\pm$ and states that there exists a 
bundle morphism which maps $I^*_+$ to $I^*_-$.
This morphism is given explicitly by the following lemma, similar to the
content of the previous subsection but with important little differences.

\begin{lemma} \label{lm:1}
Let $S$ be a closed maximal surface in an AdS manifold, and let $I^*_\pm$ be 
defined as above. Define $b:=(E+B)^{-1}(E-B)$. Then:
\begin{enumerate}
\item $b$ is self-adjoint for $I^*_+$ and has positive eigenvalues.
\item $\det(b)=1$.
\item $b$ is a solution of the Codazzi equation for $I^*_+$: 
$d^{\nabla^+}b=0$, where $\nabla^+$ is the Levi-Civit\`a connection of $I^*_+$.
\end{enumerate}
\end{lemma}

\begin{proof}
To prove the first point, consider 
two vector fields $u,v$ on $S$. Then:
$$ I^*_+(bu,v) = I((E-B)u,(e+B)v)=I((E+B)(E-B)u,v) = $$
$$ = I((E-B)(E+B)u,v) =
I((E+B)u,(E-B)v)=I^*_+(u,bv)~, $$
so that $b$ is self-adjoint for $I^*_+$. That $b$ is positive 
follows from the fact that the eigenvalues of $B$ are in $(-1,1)$
by an elementary computation. 

The second point is also an elementary computation:
$$ \det(b) = \frac{\det(E-B)}{\det(E+B)} = \frac{1+\det(B)}{1+\det(B)}=1~. $$

For the last point, let $u,v$ be vector fields on $S$; using 
Proposition \ref{pr:metric-change} for the expression of $\nabla^+$
yields that:
\begin{eqnarray*} 
d^{\nabla^+}b(u,v) & = & \nabla^+_u(bv) - \nabla^+_v(bu) - b[u,v] \\
& = & (E+B)^{-1}\nabla_u((E+B)bv) - (E+B)^{-1}\nabla_v((E+B)bu) 
- (E+B)^{-1}(E-B)[u,v] \\
& = & (E+B)^{-1}\nabla_u((E-B)v) - (E+B)^{-1}\nabla_v((E-B)u) 
- (E+B)^{-1}(E-B)[u,v] \\
& = & (E+B)^{-1} (d^\nabla(E-B))(u,v) \\
& = & 0~. 
\end{eqnarray*}
\end{proof}
Thus, $b$ satisfies all the properties of the morphism in Theorem
\ref{tm:labourie}, and it is easy to see that $I^*_-=I^*_+(b \cdot, b\cdot)$.

Now, instead of going from $I^*_+$ to $I^*_-$ one can go directly to the metric
$I$, which lies ``in between'', and the second fundamental form $B$. Once these
are found, it is easy to reconstruct the manifold. This will give a map
inverse to (\ref{max}). It is given by the following lemma:

\begin{lemma} \label{lm:2}
Let $g_+, g_-$ be two hyperbolic metrics that correspond to two points
$(c_+,c_-)\in 
\cT_\Sigma\times \cT_\Sigma$. Let $b$ be as in the theorem
\ref{tm:labourie}. Define $B:=(E+b)^{-1}(E-b)$, 
and $g:=g_+((E+B)^{-1}\cdot, (E+B)^{-1}\cdot)$, and call $\nabla$ the
Levi-Civit\`a connection of $g$. Then $B$ is self-adjoint for $g$, traceless,
with 
eigenvalues in $(-1,1)$, and it satisfies the Codazzi and Gauss equations: 
$d^\nabla B=0$ and $\det(B)=-1-K$, where $\nabla$ and $K$ are the Levi-Civit\`a
connection and the curvature of $g$.
\end{lemma}

Let us first note that $(E+B)^{-1}=(E+b)/2$. Therefore, the metric $g$ is
indeed  ``half way'' between $g_+$ and $g_-$. Note also that 
to find $g$ we could have used a morphism from $g_-$ to $g_+$ instead.

\begin{proof}
First note that, as follows from a very simple computation, 
$B$ is self-adjoint for $g_+$, with eigenvalues in $(-1,1)$, and that:
$$ b = (E+B)^{-1}(E-B)~. $$
The fact that $B$ is traceless follows directly from similar 
algebraic computations, using the fact that $\det(b)=1$. 

Since $b$ satisfies the Codazzi equation for $g_+$, it follows from 
Proposition \ref{pr:metric-change} that, if $u,v$ are two vector fields on
$S$, then: 
$$ g(Bu,v) = g_+((E+B)^{-1}Bu,(E+B)^{-1}v) = \frac{1}{4} 
g_+((E+b)(E+b)^{-1}(E-b)u, (E+b)v) = g_+((E-b)u, (E+b)v)~, $$
which is symmetric in $u$ and $v$ because $b$ is self-adjoint for $g_+$.
So $B$ is self-adjoint for $g$.

To check the Codazzi equation for $B$, note first that, because
$(E+B)^{-1}=(E+b)/2$ 
and since $b$ satisfies the Codazzi equation for $g_+$, Proposition 
\ref{pr:metric-change} shows that the Levi-Civit\`a connection $\nabla$
of $g$ is given by:
$$ \nabla_uv = (E+b)^{-1}\nabla^+_u((E+b)v)~. $$
Therefore:
\begin{eqnarray*} 
(d^\nabla B)(u,v) & = & \nabla_u(Bv) - \nabla_v(Bu) - B[u,b] =
(E+b)^{-1}\nabla^+_u((E+b)(E+b)^{-1}(E-b)v) \\
&-& (E+b)^{-1}\nabla^+_v((E+b)(E+b)^{-1}(E-b)u)- (E+b)^{-1}(E-b)[u,v] \\
& = & (E+b)^{-1}(\nabla^+_u((E-b)v) - \nabla^+_v((E-b)u) -  (E-b)[u,v]) \\
& = & (E+b)^{-1}(d^{\nabla^+}(E-b))(u,v) \\
& = & 0~.
\end{eqnarray*}
\end{proof}

The above results can be summarized in the following theorem.

\begin{thm} \label{tm:minparam}
The unique maximal surface in a globally hyperbolic AdS manifold $M$ defines
two  
metrics $I^*_\pm$ given by (\ref{I*}) that are hyperbolic. This, defines the
map 
(\ref{max}) from $\cM_{g, AdS}$ to $\cT_\Sigma\times \cT_\Sigma$. Conversely,
given 
two hyperbolic metrics $g_+$ and $g_-$ on $\Sigma$, there exists a unique
$M\in \cM_{g, AdS}$ such 
that $g_+,g_-$ are obtained from the unique maximal surface in $M$ via
(\ref{I*}). 
\end{thm}

\paragraph{Interpretations in terms of Teichm\"uller space.}

In view of the Theorems \ref{tm:ads} and Mess parameterization of the space 
$\cM_{g, AdS}$ we get a map:
\begin{equation}
\Phi_{AdS}: T^* \cT_g \to \cT_g\times \cT_g.
\end{equation}
Unlike a similar map in the hyperbolic situation, the AdS case map maps the
whole of the cotangent bundle to $\cT_g\times \cT_g$, not just a bounded
domain. 

Moreover, another such map $\Phi'_{AdS}: T^* \cT_g \to \cT_g\times \cT_g$ 
can be obtained from Theorem \ref{tm:ads}
and Theorem \ref{tm:minparam}. 

An interesting question is whether $\Phi_{AdS}$ or $\Phi'_{AdS}$ have a simple 
interpretation in terms of the geometry of Teichm\"uller space. We leave this
for future work.

\paragraph{Foliations by equidistant surfaces.}

Similarly to the hyperbolic situation, given a maximal surface $\Sigma$ in a
GHMC AdS manifold $M$, one can consider the foliation of a neighbourhood of
$\Sigma$  
by surfaces equidistant to $\Sigma$. Unlike the almost-Fuchsian situation,
however, 
the foliation only covers a neighbourhood of the minimal surface, and not the
whole manifold.

\begin{lemma}\label{lm:ads}
Let $S\subset AdS$ be a complete, oriented, smooth surface.
There exists a range of $r\in [-\pi/2,\pi/2]$, such that the set of
points $S_r$ at oriented distance $r$ from $S$ 
is a smooth embedded surface. The closest-point
projection defines a smooth map $u_r:S_r\rightarrow S$, 
and, identifying
the two surfaces by $u_r$, the induced metric $I_r$ on $S_r$ is:
\begin{equation}\label{metric-ads}
I_r(x,y) = I((\cos(r)E+\sin(r)B)x,(\cos(r)E+\sin(r)B)y)~,
\end{equation}
where $E$ is the identity operator. The shape operator of $S_r$ is: 
\begin{equation}\label{shape-ads}
B_r:=(\cos(r) E +\sin(r)B)^{-1}(\sin(r)E+\cos(r)B)~. 
\end{equation}
\end{lemma}

A proof of this lemma is analogous to the hyperbolic case and will not be repeated. Note
that the only change from the hyperbolic case is in hyperbolic trigonometric functions
being replaced by the usual ones. Note that the foliation by equidistant surfaces is smooth
only for $r\in [-r_{max},r_{max}]$, where
\begin{equation}
\frac{\cos(r_{max})}{\sin(r_{max})} = k_{max},
\end{equation}
where $k_{max}$ is the maximum of the principal curvature on $S$. 

\begin{cor}
Let $M$ be a GHMC AdS manifold, and
let $S\subset M$ be a minimal surface isotopic to $\Sigma$. 
Let $I$ and $B$ be the induced metric and
shape operator of $S$, respectively. Then there exists a neighbourhood of $S$ in $M$ that
is isometric to $\Sigma \times [-r_{max},r_{max}]$ with the metric: 
\be\label{metric-1-ads}
-dr^2 + I((\cos(r)E+\sin(r)B)\cdot,(\cos(r)E+\sin(r)B)\cdot)~. 
\ee
It is foliated by the smooth surfaces $\Sigma_r:=\Sigma\times \{r\}$.
\end{cor}

\paragraph{Fock-type metrics.}

Similarly to the hyperbolic case, one can describe the above metric in
the neighbourhood of the minimal surface in complex-analytic terms. To
this end we use the fact that $\II$ on $S$ is the real part of a
holomorphic quadratic differential. Thus, as before, we write the
metric on $\Sigma$ as: $g=e^\varphi |dz|^2$, where $z$ is the Fuchsian
uniformization coordinate, and the second fundamental form $h=tdz^2 + \bar{t}
d\bar{z}^2$.  
The metric takes the following simple form:
\begin{equation}\label{fock-ads}
ds^2 = -dr^2 + e^{\varphi} |\cos(r) dz + \sin(r) e^{-\varphi} \bar{t} d\bar{z}
|^2. 
\end{equation}
The equation (\ref{eq:*-ads}) takes the form:
\begin{equation}\label{eq:**-ads}
2\partial_{z\bar{z}} \varphi = e^{\varphi} - t\bar{t} e^{-\varphi}.
\end{equation}
Let us note that the metric (\ref{fock-ads}), as well as the equation
(\ref{eq:**-ads}), are formally obtained from the corresponding ones in 
the hyperbolic case by a simple ``analytic continuation'':
$r_{AdS}=i r_{hyp}, t_{AdS}=i t_{hyp}$. Both transformations
conform with one's intuition that the time and angular momentum
must be analytically continued when going from Lorentzian to Riemannian
space. The quadratic differential $t$ that induces a deformation 
from the Fuchsian case is also a measure of the amount of 
``rotation'' of a spacetime. This statement can be made more concrete,
but we will not pursue this in the present paper. Note that
the above analytic continuation, when considered from hyperbolic
to AdS case, is not surjective. Precisely for this reason it is
better to use another continuation, namely using the description of
both hyperbolic and AdS manifolds by two copies of Teichm\"uller space,
see the Introduction.

\section{Singular hyperbolic manifolds}

\subsection{Preliminaries: conical singularities}

\paragraph{Metrics with conical singularities.} 

Let us first define surfaces (metrics) with conical singularities. Let us
consider a  
compact Riemann surface $\Sigma$, let $p$ be a point on $\Sigma$, and $ds^2$ be
a metric on $\Sigma$, in the conformal class of $\Sigma$.
\begin{df} \label{df-cone} A point $p\in \Sigma$ is a conical singularity of
  angle $\theta\not=0$ 
of the metric $ds^2$ if there exists a non-singular conformal map $z: {\bf
  U}\to \C$ 
defined in the neighbourhood ${\bf U}$ of $x$ such that $z(p)=0$ and 
$ds^2 = \rho(z) |z|^{2(\theta/2\pi-1)} |dz|^2$ for some continuous positive
function $\rho(z)$. The value $\theta=0$ can also be allowed. However,
the metric is required to behave differently in this case: $ds^2 =
\rho(z)(\log|z|)^2 |dz|^2/|z|^2$. 
\end{df}

Another possible definition uses the asymptotic definition of the metric near
the singularities, rather than its conformal structure.
Some limiting cases deserve to be mentioned: $\theta=2\pi$ corresponds to no
singularity, whereas $\theta=0$ is the so-called puncture. 

\paragraph{Curvature of a metric with conical singularities.} 

The following result will be of importance for what follows.

\begin{thm}[M. Troyanov \cite{troyanov}] \label{tm:troyanov}
Let $\Sigma$ be a compact Riemann surface with conformal
structure $C$. Let $x_1, \cdots, x_n$ be marked points on $\Sigma$ 
and $\theta_1, \cdots, \theta_n>0$ be positive numbers. If
\begin{equation}\label{curv-sing}
2\pi(2-2g) + \sum_i^n (\theta_i-2\pi)<0
\end{equation}
then any smooth negative function $K$ on $\Sigma$ is the curvature of a unique
metric $g$, in the conformal class $c$, which has conical
singularities of angles $\theta_i$ at the points $x_i, 1\leq i\leq n$.
\end{thm}
We also note that the condition of
the theorem is always satisfied if $g\geq 2$ and $\theta\in[0,2\pi)$. We will
always 
assume the second condition in what follows. We explicitly exclude the no
singularity case $\theta=2\pi$ because this case should be treated using the
smooth setting methods of the previous section. Instead of
angle $\theta$ at the singularity, we will sometimes
use the notion of the {\it order}: $\alpha=(2\pi-\theta)/2\pi, \alpha\in(0,1]$.

\paragraph{Hyperbolic metrics with conical singularities.} 

Of special importance are metrics with conical singularities for which the
curvature function from Theorem \ref{tm:troyanov} is $-1$ everywhere on
$\Sigma$. According to Theorem \ref{tm:troyanov} there is a unique such metric 
whenever the sum of orders is greater than the Euler characteristics of the
surface: 
$\sum_i \alpha_i > 2-2g$.
 
We note that in the neighbourhood of any singular point $\theta\not=0$ the
hyperbolic metric can be written as:
\begin{equation}\label{cone}
ds^2 = dr^2 + \sinh^2(r) dt^2,
\end{equation}
where $r=0$ is the location of the singular point, and $t\in [0,\theta]$.

\paragraph{Teichm\"uller space with marked points: genus zero case.}

Our previous considerations suggest 
to consider the cotangent space of the Teichm\"uller space of conformal 
structures on $\Sigma$ with $n$ marked points. Let us remind the reader of some
well-known facts. It is instructive to consider the case of the sphere $g=0$
first.  
Thus, let us assume that $\sum_i \alpha_i > 2$ so that the condition of the
theorem \ref{tm:troyanov} is satisfied. We will also assume $\alpha_i\not=1$
for 
simplicity. However, everything can be extended to the case when punctures are
present as well. The following material is from \cite{TZ}.

By a conformal transformation any 3 points can be brought to a desired
position on $\C$. We choose 
the last 3 points $i=n-2,n-1,n$ and  put them at $0,1,\infty$
correspondingly. Thus,  
the Teichm\"uller space $\cT_{0,n}$ of the sphere with n marked points is
of complex dimension $n-3$ and can be parameterized by the positions of
the remaining $n-3$ points. The hyperbolic metric corresponding to a point in
$\cT_{0,n}$ is obtained by solving the Liouville equation:
$$\partial_{z\bar{z}} u = \frac{1}{4} e^{2u}$$
on the sphere $z\in\C$ with the conditions at the marked points $z_i,0,1$ given by:
\begin{eqnarray}
u(z) = - \alpha_i \log{|z-z_i|} + O(1) \qquad {\rm as}\,\,\, z\to z_i,
\alpha_i<1, 
\end{eqnarray}
and at $z=\infty$ given by:
\begin{equation}\label{fall-off-infty} u(z) = -(2-\alpha_n) \log{|z|}+O(1),
  \qquad {\rm as}\,\,\, z\to \infty, \alpha_n<1. 
\end{equation}
In order to see that this condition at infinity is equivalent to the
conditions at other points one 
just has to write the metric $g$ in the local coordinate $1/z$ near
infinity. Because 
the solution of the above Liouville equation is unique, the Teichm\"uller space
$\cT_{0,n}$ can 
also be viewed as the space of hyperbolic metrics on $\C$ with prescribed
conical singularities.

A Hermitian metric on the Teichm\"uller space $\cT_{0,n}$ is introduced
as follows. Let us introduce the following kernels:
$$ Q_i(z) = -\frac{1}{\pi} \left(  \frac{1}{z-z_i} + \frac{z_i-1}{z} - \frac{z_i}{z-1} \right),
\qquad i=1,\ldots,n-3.$$ 
The kernels are designed in such a way that they have a simple pole at the location of the
marked points $z=z_i$, as well as at $z=0,1$. There is no pole at $z=\infty$, actually
$Q_i(z)=O(1/z^3)$ as $z\to \infty$. 

\begin{lemma}
The functions $Q_i(z)$ are linearly independent and 
square integrable with respect to $e^{-2u}$ on $\C$.
\end{lemma}
\begin{proof}Indeed, $e^{-2u}\sim |z-z_i|^{2\alpha_i}$ at
$z\to z_i$, which makes $|Q_i(z)|^2\sim O(1/|z-z_i|^2)$ integrable for
$\alpha_i>0$. At infinity, one can similarly 
check that the behaviour (\ref{fall-off-infty}) guarantees that $|Q_i(z)|^2\sim
O(1/|z|^6)$ is integrable against 
$e^{-2u} |dz|^2$.
\end{proof}

Thus, one can define the scalar product of one-forms on $\cT_{0,n}$ by the
formula:
\begin{equation}
(dz_i,dz_j)= \frac{1}{2} \int_\C Q_i \overline{Q_j} e^{-2u} |dz|^2.
\end{equation}
This gives the matrix inverse to the Weyl-Peterson metric on $\cT_{0,n}$, and 
$Q_i(z)$ are identified with the basis in the fiber of the cotangent bundle
$T^* \cT_{0,n}$. 

\paragraph{Teichm\"uller space: general case.}

For a general Riemann surface $\Sigma$ the
tangent space $T\cT_{g,n}$ is identified with the Beltrami
differentials. The cotangent space $T^*\cT_{g,n}$, its dual, identified with
the space of meromorphic quadratic differentials on $\Sigma$ which are
holomorphic outside the $x_i$ and have at most simple poles at the $x_i$. This
condition at the $x_i$ is necessary because meromorphic quadratic
differentials having a second order pole (or worse) at one of the $x_i$ could
not be integrated against (bounded) Beltrami differentials. 

\subsection{Hyperbolic manifolds with conical singularities}

\paragraph{Cone-manifolds.}

Thurston \cite{thurston-notes} defined a general notion of cone-manifolds, the
main example being the hyperbolic cone-manifolds. We will consider here a
restricted notion, which can be defined in an elementary way. Consider first
the simplest example.

\begin{example}
Let $\theta>0$, we call $H^3_\theta$ the completion of 
$\R\times \R_{>0} \times (\R/\theta \Z)$, with the metric: 
\begin{equation}\label{cone-metric}
dz^2 + \cosh^2(z) (dr^2 + \sinh^2(r) dt^2)~, 
\end{equation}
where $z\in \R, r\in \R_{>0}$ and $t\in (\R/\theta \Z)$. 
We call $x_\theta$ the singular point corresponding to $z=r=0$.
\end{example}

A direct computation shows that $H^3_\theta$ is hyperbolic (i.e. has a smooth
metric with curvature $K=-1$) at points where $r>0$. It is not hard to see that
the above example is obtained by taking the hyperbolic plane with a 
conical singularity (\ref{cone}) and taking a ``warped product'' with the
$z$-axis. The procedure is the same as the one used in obtaining
a Fuchsian manifold from its $z=0$ section. Our main aim here is to
study how this procedure is generalized to non-Fuchsian singular manifolds.

Note that one can in addition allow the special value $\theta=0$. This case
is, actually, a bit simpler than a general angle deficit, for it
can be treated using the discrete groups acting on $H^3$. The 
hyperbolic manifold $H^3_0$ containing a single singularity of this type is
obtained as the quotient of $H^3$ by the group $\Gamma$ generated
by a single parabolic element of ${\rm PSL}(2,\C)$. Thus, we will
allow the value $\theta=0$ in what follows, assuming that (\ref{cone-metric})
is replaced by such $H^3_0$.

\begin{df}
A {\bf hyperbolic cone-manifold} is a metric space $M$ in which any point $x$
has a  
neighbourhood isometric to a subset of $H^3_\theta$, for some $\theta\geq 0$.
If $\theta$ can be taken equal to $2\pi$ then $x$ is a {\it smooth} point of
$M$,  
otherwise $\theta$ is uniquely determined and is the {\bf total angle} at $x$.
\end{df} 

It follows directly that the singular locus of a
hyperbolic cone-manifold -- following the definition given here -- is a
disjoint union of curves. To each of those curves is attached a number, which
is equal at each point to the number $\theta$ appearing in the definition, and
which we call the {\bf total angle} around the singular curve.

\paragraph{Quasi-Fuchsian cone-manifolds.}

We are interested here in a natural extension of the notion of quasi-Fuchsian
manifold, which allows for conical singularities along lines, with a total
angle around the singular lines which lies in the interval $[0,2\pi)$. The restriction
on $\theta$ from below does not need an explanation, $\theta$ is a positive number
by definition. The restriction from above does deserves an explanation. 
The angle $\theta$ can in principle be larger than $2\pi$, which would
correspond to having a negative angle deficit, as long as  the condition 
$ 2\pi(2-2g) + \sum_i^n (\theta_i-2\pi)<0$ is satisfied. However, as we shall
see below, a certain natural geometric condition requires $\theta\leq 2\pi$.
We will then explicitly exclude $\theta=2\pi$ because there is
no singularity in this case, so one should treat such situations as in the smooth
setting that we were considering so far.  

\begin{df}
A hyperbolic cone-manifold $M$ is {\bf quasi-Fuchsian} if:
\begin{itemize}
\item $M$ is complete.
\item $M=\Sigma \times \R$, where $\Sigma$ is a closed surface of genus at
 least $2$. 
\item The singular set of $M$ is the disjoint union of a non-compact curves
  going  
from $\Sigma \times \{ -\infty\}$ to $\Sigma \times \{ +\infty\}$.
\item The total angle $\theta$ around each singular curve is in $\theta\in[0,2\pi)$.
\item $M$ contains a compact subset $K$ which is convex, i.e. 
any geodesic segment with endpoints in $K$ is contained in $K$.
\end{itemize}
\end{df}

\paragraph{Surfaces in cone manifolds.}

We will assume $\theta>0$ in this paragraph. 
Let us consider a quasi-Fuchsian hyperbolic cone manifold $M$ and a
hypersurface  
$S\subset M$ that intersects each curve from the singular locus exactly
once. It is 
instructive to consider what are the possibilities for the induced metric
on $S$. The main question that we would like to address is whether the
metric on $S$ contains conical singularities and if yes what are the
corresponding angles. It is clear that to answer this question one only
has to consider the neighbourhood of the intersection of $S$ with
the singular locus. By definition, in the neighbourhood of every singular
curve the metric is given by (\ref{cone-metric}) and we just have to
understand what the induced metric on the surface intersecting the
$z=0$ singular curve can be. 

Let us first note that the metric 
(\ref{cone-metric}) can be obtained from the hyperbolic space $H^3$ by
the procedure of cutting out a wedge. Indeed, let us consider
the unit ball model of $H^3$, and consider a rotation around the
geodesic connecting the north and south poles by the angle of
$2\pi-\theta$. This 
rotation maps half a geodesic plane $P_1$ to another half of a geodesic plane
$P_2$, with $P_1, P_2$ intersecting along the vertical axis. One can now
identify 
the halves of geodesic planes $P_1, P_2$ to obtain the space
(\ref{cone-metric}),  
which is thus just the hyperbolic space $H^3$ with the wedge between $P_1,
P_2$ removed.  

Now consider a surface $S$ in (\ref{cone-metric}). Any such surface can be
obtained from a surface $S'$ in $H^3$ with the part of $S'$ in between $P_1$
and $P_2$ 
removed. It is clear however that not any surface $S'$ will do: 
the curves along which $S'$ intersects the half-planes $P_1, P_2$
should go one into the other under the $2\pi-\theta$ rotation.
Any surface $S$ obtained this way has an induced metric, and it is
clear that this metric contains a conical singularity in the sense
of definition \ref{df-cone}. Its angle, however, does not have
to be equal to the total angle around the curve, as we shall soon
see. The angle of the conical singularity of the metric on $S$ is
equal to the total angle around the singular line if and only if 
the surface $S$ intersects the singularity ``orthogonally''.
However, the notion of orthogonality requires a careful definition
in this case. We will use the following notion of orthogonality.
We first define a geodesic plane orthogonal to the singular line to
be the one isometric to the plane $\{z=0\}$ in (\ref{cone-metric}).

\begin{df}
Let $M$ be a hyperbolic cone-manifold, and let $x\in M$ be a singular
point. Then there exists a neighbourhood $U$ of $x$ in $M$ such that
the intersection of $U$ with $M$ is isometric to a neighbourhood $U'$ of 
the singular point $x_\theta$ in $H^3_\theta$, for some $\theta>0$. We define a
{\it geodesic plane orthogonal to the singular set} at $x$ as a subset $S\subset M$
such that: 
\begin{itemize}
\item $S\setminus \{ x\}$ is a totally geodesic surface in the regular
set of $M$.
\item $S\cap U$ corresponds, under the isometry with $U'$, to the intersection
  with $U'$ of the totally
geodesic plane $\{ z=0\}$.
\end{itemize}
\end{df}

We now define a plane intersecting the singular line orthogonally to be
one that infinitesimally close to the singular point 
intersects it as an orthogonal geodesic plane.

\begin{df}\label{df-orth}
Let $S$ be a subset of $M$ which is a smooth surface outside the singular set
of $M$. $S$ is orthogonal to the singular set if, for each $x\in
S$ which is a singular point of $M$, there exists a neighbourhood $U$ of $x$
in $M$ and a plane $P$ orthogonal to the singular set at $x$ such that:
 $$ \lim_{y\in S, y\rightarrow x} \frac{d(y,P)}{d(y,x)} = 0~. $$
\end{df}

It is clear from these definitions that each plane intersecting a conical
singularity 
line of angle $\theta$ orthogonally 
has an induced metric containing the conical singularity of
angle exactly $\theta$. Let us now consider what happens when the intersection
is generic. Of special interest are non-orthogonal geodesic planes. These can
be 
defined similarly to orthogonal geodesic planes. Thus, consider a neighbourhood
$U$ of the intersection point $x$, with $U$ being isometric to $U'$, subset of
the 
cone space (\ref{cone-metric}). Consider a totally geodesic plane $S'$ in
$H^3$ that 
intersects the planes $P_1, P_2$ along two curves that go into one another
under the rotation. It is clear that this plane $S'$ descends to the cone space
after the wedge between $P_1, P_2$ is removed. However, the plane $S$ obtained
is now rather special: the wedge on $S'$ that is removed to obtain $S$ is
located 
either ``behind'' or ``in front'' of the point particle. In other words, the
conical singularity line, which is the vertical axis in $H^3$, when projected
onto $S'$ lies exactly in between the lines of intersection of $S'$ with
$P_1, P_2$. Having this way constructed a geodesic plane intersecting the
conical line 
at an angle, one can easily convince oneself that the induced metric on $S$ is
that with a conical singularity, but the angle of this singularity is 
different from $\theta$. It is also clear that the geodesic surface
$S$ is embedded into the conical space (\ref{cone-metric}) in a non-smooth way.
Indeed, unless $S$ is orthogonal to the singular line, there is a fold on $S$
extending away from the particle, and coinciding with the line which
came from two lines on $P_1, P_2$ that were identified. Let us note that
there is a simple relation between the conical angle on the surface,
the total angle around the singular curve, and the bending angle.
It can be easily worked out, but we won't need the precise relation here. 

It is now clear that infinitesimally close to the intersection point
a generic intersection of $S$ with the singular line is like the one with
a geodesic plane. It is thus clear that a plane
intersecting the conical singularity line in a generic way contains
a conical singularity with the angle not equal to $\theta$ (unless the
intersection is orthogonal), and with a line of folding extending from
the conical defect. This picture admits a nice interpretation in physics
terms: the deficit angle as measured on the surface $S$ is the 
particle's energy,
while the bending angle along the line of the fold is a measure
of the particle's momentum. The two always satisfy a relation that connects
them to the total angle around the singular line, which is interpreted
as particle's mass. 

\paragraph{Minimal surfaces in hyperbolic cone-manifolds.}

The main goal of this section is to extend the results and the methods of
section 2 from hyperbolic manifolds to hyperbolic cone-manifolds. The first
step is therefore to define the correct notion of minimal surfaces in
hyperbolic cone-manifolds. 

As we have seen in the previous paragraph, a generic surface $S$ intersects 
the conical singularity line in such a way that a fold is formed on it. One
of the principal curvatures blows up along the fold line, the other stays
finite 
(zero if the surface is geodesic). For this reason it is easy to see that
a minimal surface cannot intersect the conical defect lines in any other way
but orthogonally. Indeed, the first variation of the area of a surface is
given by:
$$ A' = 2\int uH da~. $$
This cannot vanish unless $S$ intersects the conical defect lines orthogonally.
Indeed, if the intersection is not orthogonal, then there is always a fold
(or a system of folds) on $S$, and the first variation cannot vanish for
a general infinitesimal variation of the surface: there will always be
a contribution to the above integral from the fold lines. As the
above discussion was heuristic rather than 
rigorous, we build the condition of orthogonality
into the definition of the minimal surface.

\begin{df}
Let $M$ be a quasi-Fuchsian cone-manifold, and let $S\subset M$. $S$ is a {\bf
  minimal surface} in $M$ if:
\begin{itemize}
\item $S$ is a smooth surface of zero mean curvature
outside the singular set of $M$,
\item it is orthogonal to the singular set.
\end{itemize}
\end{df}

The following statement is analogous to what we had in the non-singular case.

\begin{remark} \label{rk:second-sing}
Let $S$ be a minimal surface in $M$, with principal curvatures less than $1$
at each (non-singular) point. Then $S$ is a local minimum of the area among
surfaces orthogonal to the singular locus.  
\end{remark}

\begin{proof}
Consider a normal deformation of $S$ of the form $uN$, where $N$ 
is the unit normal to $S$. The classical formula for the second 
variation of the area of $S$ is (see e.g. \cite{spivak}): 
$$ A'' = \frac{1}{2} \int u\Delta u + 2(1-k^2)u^2 da~. $$
Integrating by parts and using the fact that $S$ remains smooth all the way 
to the singular locus shows that: 
$$ A'' = \frac{1}{2}\int |du|^2 + 2(1-k^2)u^2 da~, $$
and it follows that $A''$ is positive if $|k|<1$ and $u$ is not identically 
$0$.
\end{proof}

With the definition of minimal surfaces
used here it would be interesting to know whether the classical
result on the existence of minimal surfaces in quasi-Fuchsian manifolds can be
extended to quasi-Fuchsian cone-manifolds: 

\begin{question} \label{q:minimal}
Does every quasi-Fuchsian cone-manifold contain a closed minimal surface ?  
\end{question}

It is possible that one has to assume that the angle around the singular curves is less than $\pi$ 
to get a positive answer to this question. The main results of this section are that many
almost-Fuchsian cone-manifolds can be constructed from minimal surfaces; a
positive answer to question \ref{q:minimal} would show that all
almost-Fuchsian cone-manifolds can be obtained in this manner.

\paragraph{Almost-Fuchsian cone-manifolds.}

The definition of minimal surfaces in quasi-Fuchsian cone-manifolds leads
naturally to the definition of almost-Fuchsian cone-manifolds, as for smooth
manifolds. The uniqueness result for minimal surfaces in almost-Fuchsian
manifolds also extends to the singular setting.

\begin{df}
An {\bf almost-Fuchsian} cone-manifold is a quasi-Fuchsian cone-manifold which
contains a closed, embedded minimal surface, which has principal curvatures in
$(-1,1)$. 
\end{df}

\begin{remark} As we shall see in the next subsection, the total angles of
an almost-Fuchsian cone manifolds are necessarily in $[0,2\pi]$. 
\end{remark}

\subsection{Quasi-Fuchsian cone-manifolds from minimal surfaces}

\paragraph{Minimal surfaces in germs of cone-manifolds.}

As in the non-singular case, it is helpful to consider a notion of 
minimal surface based on considering the possible induced metrics and second
fundamental forms, but now on a surface with marked points.

\begin{df}
A {\bf minimal surface in a germ of hyperbolic cone-manifold} is a couple
$(g,h)$ on $\Sigma$ with $n$ marked points $x_1, \cdots, x_n$, such that:
\begin{itemize}
\item $g$ is a smooth metric on $\Sigma$ with conical singularities at the
  $x_i$, with total angle $\theta_i\in[0,2\pi)$
  at $x_i, 1\leq i\leq n$.
\item $h$ is a bilinear symmetric form on $T\Sigma$ defined outside the $x_i$,
which is traceless with respect to $g$.

\item the determinant of $h$ with respect to $g$ remains bounded.
\item $h$ satisfies the Codazzi and Gauss equations with respect to $g$:
$$ d^\nabla h=0~ , ~ ~ K_g = -1 + {\rm det}_g(h)~, $$
where $\nabla$ is the Levi-Civit\`a connection of $g$ and $K_g$ is its
curvature.  
\end{itemize}
\end{df}

Let us first note that, because $g,h$ satisfy the Gauss equation, the
singularity of 
$h$ at $x_i$, which we assumed to be that of a pole, cannot be of arbitrary
order. 

\begin{lemma} Let $(g,h)$ be a minimal surface in a germ of cone-manifold. The
second fundamental form is the real part of a meromorphic quadratic 
differential $q$ which has at most a first order pole at $x_i$. 
If the angle $\theta_i$ at $x_i$ is in $(\pi, 2\pi)$ then $q$ has no
pole at $x_i$.
\end{lemma}

\begin{proof}
The integral of the left hand side of the Gauss equation exists, and is equal
to 
(\ref{curv-sing}). Thus, the right hand side must also be integrable. We know
that 
the metric $g$ behaves as $r^{-2\alpha}$ at the singular points, where
$\alpha$ is the 
order, and $r$ is the distance to the singularity measured using some
non-singular,  
e.g. flat, metric. Because $h$ is traceless, it can be written as the real
part of a quadratic differential t: $h=tdz^2 + \bar{t} d\bar{z}^2$. Then
near the singularities ${\rm det}_g(h) da$, where $da$ is the area form of
$g$,  
behaves as $t\bar{t}\, r^{2\alpha} r dr$. This is bounded for all
$\alpha\in(1/2,1]$ iff $t$ has at most a first order pole at the 
singularity, and is bounded for all $\alpha\in(0,1]$ if and only 
if $t$ has no pole at at the singularity.
\end{proof}

\begin{cor} \label{lm:cone-hqd}
Let $(g,h)$ be a minimal surface in a germ of cone-manifold. There exists
$c\in \cT_{g,n}$ and $t\in T^*_c\cT_{g,n}$ such that $g$ is in the conformal
class defined by $c$ and $h=Re(t)$.
\end{cor}

\begin{proof}
As the above lemma shows, the quadratic differential $t$ has at most a first
order pole at the singular points. It is moreover holomorphic with respect 
to the complex structure defined by $g$ due to Codazzi equation. Thus, the
metric $g$ gives us $c\in \cT_{g,n}$, and together with $t$ this gives
$(c,t)\in T^* \cT_{g,n}$.
\end{proof}

Let us also explain why it is natural to restrict the range of angles $\theta$
from  
above. It is natural to consider the third fundamental form of the surface $S$,
and require that the area form it defines is integrable on $\Sigma$. This area
form 
is equal to ${\rm det}_g(h) da$, where $da$ is the area form of $g$. We already
know that $h$ is the real part of a meromorphic quadratic differential 
$q$ which can have poles at singular points, and those poles are at most
first order. We will require the quantity ${\rm det}_g(h) da$ to be 
integrable even when $h$ has such first order poles. This is equivalent
to the condition that $r^{2\alpha-2} rdr$ is integrable, which requires
$\alpha>0$. 

\begin{remark} It is easy to see that when a total angle at one of the
singularities is in $(\pi,2\pi)$ then the principal curvature at 
that point diverges as soon as the meromorphic quadratic differential
$q$ has a pole.
Indeed, then the principal curvature behaves as $k\sim r^{2\alpha - 1}$. 
This diverges
for $\alpha\in(0,1/2)$ or for total angles larger than $\pi$.
\end{remark}

\paragraph{Singular minimal surfaces from $T^*\cT_{g,n}$.}

Conversely, let $(c,t)\in T^* \cT_{g,n}$ and consider 
the question whether $(c,t)$ is
associated to one -- or several -- minimal surfaces with singularities. 
This question is related to the problem of existence/uniqueness of solutions of the same 
partial differential equation as in section 2. Given $\alpha_1, \cdots, \alpha_n\in (0,1)$, 
by Theorem \ref{tm:troyanov} there exists a metric $g_0$ on $\Sigma$ with curvature in 
$(-1,0)$, with conical singularities of orders $\alpha_i$ at $x_i, 1\leq i\leq n$.

\begin{lemma}
Let $g=e^{2u}g_0$, then $(g,h)$ is a minimal surface in a germ of hyperbolic
cone-manifolds, with orders $\alpha_i$ at $x_i, 1\leq i\leq n$, if
and only if $u$ is bounded and is a solution on $\Sigma\setminus \{ x_1,
\cdots, x_n\}$ of equation (\ref{eq:*}).
\end{lemma}

\begin{proof}
The proof of Lemma \ref{lm:confchange} shows that $g$ and $h$ satisfy the
Gauss and Codazzi equations on the complement of the $x_i$ if and only if $u$
is a solution of (\ref{eq:*}). If $u$ is bounded, then $g$ is clearly a
cone-manifold with the same singular curvatures as $g_0$. Conversely, if 
$g$ has the same singular curvatures as $g_0$, then it has the same 
asymptotic behaviour at each of the $x_i$ with respect to a smooth conformal
reference metric, so that $u$ is bounded near each of the $x_i$.
\end{proof}

In case when there exist a solution of (\ref{eq:*}), and the principal curvatures are in $(-1,1)$, 
this solution is unique. This gives an extension of 
lemma \ref{unique} from smooth to singular minimal surfaces setting.

\begin{lemma} \label{lm:unique-cone}
For each $(c,t)\in T^*\cT_{g,n}$, there is at most one 
almost-Fuchsian minimal surface in a germ of hyperbolic cone-manifold,
with given angles $\theta_1, \cdots, \theta_n\in (0,\pi)$ at the 
singular points,
which corresponds to $(c,q)$ under Lemma \ref{lm:cone-hqd}. 
\end{lemma}

\begin{proof}
Suppose that there are two such surfaces, say $(g,h)$ and $(g',h')$. Then
by definition $h=h'$, while $g'=e^{2u}g$, for some function
$u:\Sigma\rightarrow \R$. Then, as in the proof of Lemma \ref{unique}, $u$ is
a solution of the equation:
$$ \Delta u = (1-e^{2u}) + k^2(1-e^{-2u})~. $$ 
Moreover, since $g$ and $g'$ are two metrics with the same angles at the
conical points, $u$ is bounded (see \cite{troyanov}). And the argument in the
proof of Lemma \ref{unique} clearly shows that $u$ can not have a positive
maximum except perhaps at one of the conical points. 

Let $x_0$ be one of the conical points of $\Sigma$, and suppose that $u$ has a
maximum at $x_0$. There is a neighbourhood of $x_0$ in which the metric can be
written as:
$$ g = |z|^{-2\alpha} |dz|^2~, $$
where $\alpha$ is the order of the singularity at $x_0$. If $\Delta_0$ denotes
the Laplace operator of the flat metric, it follows from the classical
conformal invariance of the Laplace operator on surfaces that:
$$ \Delta =  |z|^{-2\alpha} \Delta_0~, $$
and, taking into account the equation which is satisfied by $u$:
$$ \Delta_0 u = |z|^{2\alpha}((1-e^{2u}) + k^2(1-e^{-2u}))~. $$ 
Since the left-hand side is again non-positive when $u\geq 0$, it follows from
the maximum principle that $u(x_0)\leq 0$. So $u$ is non-positive on
$\Sigma$. Repeating the same argument after exchanging $g$ and $g'$, as in the
proof of Lemma \ref{unique}, therefore shows that $u=0$ on $\Sigma$, so that
$g=g'$. 
\end{proof}

\paragraph{Almost-Fuchsian manifolds from minimal surfaces.}

Minimal surfaces with conical singularities, just like smooth minimal
surfaces in section 2, can be used to construct hyperbolic metrics. The
only difference is that the metrics constructed in this way have conical
singularities, and they correspond precisely to the ``quasi-Fuchsian
cone-manifolds'' defined above.

\begin{lemma} \label{lm:414}
Let $(g,h)$ be an almost-Fuchsian minimal surface in a germ of hyperbolic 
cone-manifold. Let $B:T\Sigma\rightarrow T\Sigma$ be the corresponding shape 
operator, i.e. $h(\cdot, \cdot)=g(B\cdot, \cdot)$. Define a metric $G$ on
$\Sigma\times \R$ as:
$$ G = dt^2 + g((\cosh(t)E+\sinh(t)B)\cdot,(\cosh(t)E+\sinh(t)B)\cdot)~. $$
Then $(\Sigma\times \R, G)$ is an almost-Fuchsian cone-manifold, with conical
singularities on the curves $\{ x_i\}\times \R$, containing $\Sigma
\times \{ 0\}$ as a minimal surface.
\end{lemma}

Without the assumption that $(g,h)$ is almost-Fuchsian, the conclusion is only
valid for $(\Sigma\times (-r,r), G)$, where $r$ is determined by the maximum
of the eigenvalues of $B$ on $\Sigma$.

\begin{proof}
It follows from the arguments in the proof of Lemma \ref{lm:27} that $G$ is a 
(smooth) hyperbolic metrics on $(\Sigma\setminus \{ x_1, \cdots, x_n\})\times 
\R$. In addition, the total angle around each of the curves $\{ x_i\}\times \R$
is easily seen to be equal to the total angle around $x_i$ in $\Sigma$, so 
that $\{ x_i\}\times \R$ has a conical singularity of total angle $\theta_i$. 
\end{proof}

Clearly, every almost-Fuchsian hyperbolic cone-manifold is of this form, since
by definition it contains an almost-Fuchsian minimal surface. 

\begin{lemma} \label{lm:unique-sing}
An almost-Fuchsian cone-manifolds contains only one closed, embedded minimal
surface. 
\end{lemma}

\begin{proof}
The argument given for (smooth) almost-Fuchsian works just as well here, since
the surfaces $\Sigma \times \{ r\}$ have mean curvature vector towards 
$\Sigma\times \R$; since minimal surfaces are required to be orthogonal to 
the singularities, the maximum of $r$ on a minimal surface has to be 
non-positive by the maximum principle, while its minimum has to be
non-negative, 
so that the only embedded minimal surface is $\Sigma \times \{ r\}$.
\end{proof}

We can now sum up the results of this subsection as follows.

\begin{thm} \label{tm:hyper-cone-main}
Choose $\theta_1, \cdots, \theta_n\in (0,\pi)$, and let $x_1, \cdots, x_n$ 
be distinct points on $\Sigma$.
The space of almost-Fuchsian metrics on $\Sigma\times \R$, with conical 
singularities along the lines $\{ x_i\}\times \R$ of angles 
$\theta_1, \cdots, \theta_n$, for $1\leq i\leq n$,
considered up to isotopies fixing the singular lines,
is parametrized by an open subset of $T^*\cT_{g,n}$.
The point $(c,q)\in T^*\cT_{g,n}$ associated to an almost-Fuchsian
metric $G$ is defined as follows: $c$ is the conformal structure of the induced
metric, and $q$ corresponds to the second fundamental form, of the unique
minimal surface contained in $(\Sigma\times \R, G)$.
\end{thm}

\subsection{The conformal structure at infinity}

Let $M$ be a quasi-Fuchsian cone-manifold. Each convex subset of $M$ 
which is disjoint from the singular locus is isometric to a subset of
hyperbolic 3-space, so that, if it is infinite, it has a well-defined
boundary at infinity. The boundaries at infinity of the infinite convex
subsets can be glued together to obtain a natural boundary at infinity
of $M$, which is made of two disjoint copies of a surfaces of genus $g$,
with $n$ points removed, corresponding to the endpoints of the singular
curves in $M$. In addition it is natural to assign to each point which
is removed in the boundary at infinity an angle, namely the total angle
around the corresponding component of the singular locus of $M$. 
One obtains in this way two Riemann surfaces, each with $n$ marked
points and with an angle in $[0, 2\pi)$ attached to each marked
point. The numbers attached to 
corresponding points on the two  boundary components are the same.

This construction defines a natural map:
$$ \Phi_{g,n}:\cM_{g,n}\rightarrow \cT_{g,n}\times\cT_{g,n}~,
$$ 
which is an extension to Riemann surfaces with marked points of the 
Bers map sending a (smooth) quasi-Fuchsian manifold to the conformal
structures on the two connected components of its boundary at infinity.

For almost-Fuchsian cone-manifolds, the conformal structure at infinity
has a simple expression in terms of the minimal surface.

\begin{remark}
Let $M$ be an almost-Fuchsian cone-manifold, containing an embedded minimal
surface $S$ with induced metric $I$ and shape operator $B$. Then the 
conformal structures at infinity of $M$ contain the metrics: 
$$ I^*_\pm = I((E\pm B)\cdot, (E\pm B)\cdot)~, $$
with the marked points corresponding to the intersections of $S$ with 
the singular curves in $M$ and the angles at the marked points corresponding
to the angle around the corresponding singular curves.
\end{remark}

\begin{proof}
The conformal class of the boundary at infinity of $M$ is obtained by 
considering the metrics on the surfaces at distance $r$ from $S$. 
Those surfaces are smooth by Lemma \ref{lm:414}, and the induced metric on
those surfaces is also explicit from Lemma \ref{lm:414}. The result then 
follows from this explicit expression and a direct normalization argument.
\end{proof}

It is therefore natural to ask whether the analog of the Ahlfors-Bers
theorem holds.

\begin{question} \label{q:ab-sing}
Is the map $\Phi_{g,n}$ is a homeomorphism ?
\end{question}

One of the key points towards an answer would be to understand whether
$\Phi_{g,n}$ is locally injective, i.e. whether its differential is
an isomorphism. It is possible to tackle this question using extensions
to a non-compact setting of the rigidity argument of Hodgson and Kerckhoff
\cite{HK}, but this entails some technical difficulties (it is the subject
of a work currently in progress by S. Moroianu and the second author).

\section{Singular AdS manifolds}

\subsection{AdS cone-manifolds and maximal surfaces}

\paragraph{AdS cone-manifolds.}

It is quite straightforward to define GHMC AdS cone-manifolds, by analogy with
both the quasi-Fuchsian hyperbolic cone-manifolds and the non-singular GHMC
AdS manifolds. As in the hyperbolic case, the first step is to define the
local structure near a conical singularity. 

\begin{example}
Let $\theta>0$, we call $AdS^3_\theta$ the completion of 
$\R\times \R_{>0} \times (\R/\theta \Z)$, with the metric: 
\begin{equation}\label{cone-metric-ads}
- dz^2 + \cos^2(z) (dr^2 + \sinh^2(r) dt^2)~, 
\end{equation}
where $z\in \R, r\in \R_{>0}$ and $t\in (\R/\theta \Z)$. 
We call $x_\theta$ the singular point corresponding to $z=r=0$.
\end{example}

Clearly the metric is Lorentzian outside the singular points. 
As in the hyperbolic case, a straightforward computation shows that it has
constant curvature $-1$ outside the singular line, so that it is locally
modelled on the Anti-de Sitter space. Actually this space is obtained from the
complete hyperbolic surface with a single conical singularity of angle
$\theta$, by taking a warped product with $\R$, in the same manner as the AdS
space is obtained from the hyperbolic plane.

\begin{df}
An {\bf AdS cone-manifold} is a metric space $M$ in which any point $x$
has a  
neighbourhood isometric to a subset of $AdS^3_\theta$, for some $\theta\geq 0$.
If $\theta$ can be taken equal to $2\pi$ then $x$ is a {\it smooth} point of
$M$, otherwise $\theta$ is uniquely determined and is the {\bf total angle} at
$x$. 
\end{df} 

Note that, as in the hyperbolic case, more general definitions could be given,
in particular it would be conceivable to allow singularities also along
space-like curves; this would correspond to generalized ``black holes''. 
The definition given here is kept as simple as possible,
given the applications that we have in mind.

\paragraph{GMHC cone-manifolds.}

We are interested here in a particular class of AdS cone-manifolds, defined as
in the corresponding smooth case. The first point is to define a notion of
surface orthogonal to the singular locus. As in the hyperbolic case, the space
$AdS^3_\theta$ can be obtained from $AdS^3$ by removing a wedge. 
There is therefore a natural notion of totally geodesic
plane in $AdS^3$ orthogonal to the singular line. Moreover, each point in the
singular line is contained in exactly one such orthogonal totally geodesic
plane. Given such a plane, say $P$, there is a natural notion of {\it
  distance} to $P$ in $AdS^3_\theta$, defined as the length of the (time-like)
geodesic going from a point in $AdS^3_\theta$ to $P$ and orthogonal to $P$.
One can then define a general surface $S$ orthogonal to the singularity line
to be one that intersects the line like a geodesic plane.

\begin{df}
Let $S\subset AdS^3_\theta$ be a space-like surface, intersecting the singular
line at a point $x$. $S$ is {\it orthogonal to the singular locus} at $x$ if
the distance to the geodesic plane $P$ orthogonal to the singular locus at $x$
is such that: 
$$ \lim_{y\rightarrow x, y\in S}\frac{d(y,P)}{d_S(x,y)} = 0~, $$
where $d_S(x,y)$ is the distance between $x$ and $y$ along $S$.

If now $S$ is a space-like surface in an AdS cone-manifold $M$, intersecting the
singular locus at a point $x'$, it is orthogonal to the singular locus at $x'$
if there exists a neighbourhood $U$ of $x'$ in $M$ which is isometric to a
neighbourhood of a singular point in $AdS^3_\theta$, the isometry sending
$S\cap U$ to a surface orthogonal to the singular locus in $AdS^3_\theta$.
\end{df}

The definition of a GHMC AdS manifold --- where GHMC stands for ``globally
hyperbolic maximal compact'' --- is analog to the corresponding definition in
the non-singular case. 

\begin{df}
An AdS cone-manifold is $M$ is {\it GHMC} if: 
\begin{itemize}
\item $M$ contains a closed space-like surface $S$ orthogonal to the singular
  locus. 
\item $S$ intersects each time-like geodesic, disjoint from the singular
  locus, at exactly one point.
\item If $N$ is another AdS cone-manifold satisfying the previous
assumptions and $\phi:M\rightarrow N$ is an
  injective isometry then $\phi$ is one-to-one.
\end{itemize}
\end{df} 

We will need some specific notations for the space of GHMC AdS 
cone-manifolds with
given singular angles. These are provided by the following definition.

\begin{df}
Let $g\geq 2$ and $n\geq 1$, and let $\theta_1,\cdots, \theta_n>0$.
Consider a closed surface of genus $g$ with $n$ marked points $x_1, \cdots, 
x_n$. We call $\cM'_{g,n,AdS}(\theta_1, \cdots, \theta_n)$ the space
of GHMC AdS metrics on $\Sigma\times \R$, with conical singularities at the
lines $\{ x_i\}\times \R$, and with angle $\theta_i$ around the
corresponding singular curve. We call 
$\cM_{g,n,AdS}(\theta_1, \cdots, \theta_n)$
the universal cover of the quotient of 
$\cM'_{g,n,AdS}(\theta_1, \cdots, \theta_n)$ by the diffeomorphisms 
sending $\{ x_i\}\times \R$ to $\{ x_i\}\times \R$, for $1\leq i\leq n$.
\end{df}

Elements of $\cM_{g,n,AdS}(\theta_1, \cdots, \theta_n)$ will sometimes
be called ``GHMC AdS cone-manifolds'', although there is a clear abuse of
notation here, since there are many elements of  
$\cM_{g,n,AdS}(\theta_1, \cdots, \theta_n)$ which correspond to the
same GHMC AdS cone-manifold (they differ by diffeomorphisms).

We are particularly interested here in the situation where $S$ has genus at
least equal to $2$, and where the angle around each singular line is less than 
$2\pi$ --- and sometimes even less than $\pi$.

It is not difficult to check that if $M$ is a GHMC AdS manifold, containing a
closed space-like surface $S$, then $M$ is topologically the product of $S$ by
an interval. 

\paragraph{Maximal surface.}

The notion of maximal surfaces that we consider includes the condition of
orthogonality to the singular locus. The reasons for including the
orthogonality condition into the definition are basically the
same as those in the hyperbolic case, and won't be repeated here.

\begin{df}
A closed, space-like, embedded surface $S$ in an AdS cone-manifold is 
{\it maximal} if: 
\begin{itemize}
\item it is smooth and has mean curvature $0$ outside the singular locus,
\item it is orthogonal to the singular locus,
\item its principal curvatures are bounded.
\end{itemize}
\end{df}

\paragraph{Two questions.}

There are two interesting questions which arise concerning maximal surfaces in
GHMC cone-manifolds. The first concerns the existence of such a surface when
the angles around the singular angles are small enough.

\begin{question} \label{q:q1}
Does every GHMC AdS cone-manifold, with singular angle less than $\pi$
(resp. $2\pi$) at each
singular line, contain a maximal surface ? Is this maximal surface unique ?
\end{question}

It appears conceivable that a positive answer to the existence part of the
question could follow from arguments
close to those standard in the theory of maximal surfaces in Lorentzian
manifolds. A positive answer to the uniqueness part of the question also
appears to be within reach, using informations on the geometry of GHMC AdS
cone-manifolds and the maximum principle. A positive answer to both parts of
the question would also follow directly from a positive
answer to question \ref{q:q3} below. 

The second question is about the possibility to get, from a maximal surface, a
CMC foliation of the whole manifold. This is possible for non-singular GHMC
AdS manifolds, see \cite{BBZ-cras}.

\begin{question} \label{q:q2}
Let $M$ be a GHMC AdS cone-manifold, containing a maximal surface $S$. Is
there a foliation of $M$ by space-like CMC surfaces, one of them being $S$?
\end{question}

Again it appears conceivable that a proof can be obtained by those same 
methods which are used in the non-singular case.

\paragraph{Example: Fuchsian GHMC cone-manifolds.}

As in the hyperbolic setting, the basic example of cone-manifolds is
provided by what can be called ``Fuchsian'' manifolds. All other
manifolds have to be thought of as arising via a deformation of
the Fuchsian case. The ``Fuchsian'' GHMC AdS cone-manifolds 
 are built from a closed surface
$\Sigma$ with $n$ marked points $x_1, \cdots, x_n$, and with a hyperbolic
metric $g$ with conical singularities at the $x_i$. At this point 
any value (in $[0,2\pi]$) of the total angle is possible. However,
as we shall see later, the cases $\theta\in[0,\pi]$ and
$\theta\in[\pi,2\pi]$ are actually treated very differently.

The ``Fuchsian'' GHMC metric is defined on $\Sigma \times (-\pi/2,\pi/2)$ as 
$G:=-dt^2 + \cos^2(t) g$. It is not difficult to check that it is indeed a GHMC
AdS metric --- indeed it follows from the definition that the metric is
everywhere of the prescribed form, including along the singular lines $\{
x_i\}\times (-\pi/2,\pi/2)$. Clearly the surface $\{ t=0\}$ is totally
geodesic, whence maximal. It is also quite clear that the surfaces $\{
t=t_0\}$, for $-\pi/2<t_0<\pi/2$, are umbilic, and therefore have constant
mean curvature, and that their mean curvature varies from $-\infty$ to
$+\infty$. 

\begin{remark}
$(\Sigma \times (-\pi/2,\pi/2), G)$ contains no other space-like, closed,
embedded maximal surfaces apart from $\Sigma\times \{0\}$.
\end{remark}

\begin{proof}
It follows directly from the maximum principle, applied using the CMC foliation
by umbilic surfaces.  
\end{proof}

\paragraph{Maximal surfaces in germs of AdS cone-manifolds.}

It is quite natural to extend to the singular case the definition given above
for maximal surfaces in germs of AdS manifolds. 

\begin{df}
Let $g\geq 2$ and $n\geq 1$, and let $\theta_1,\cdots, \theta_n>0$.
A maximal surface of genus $g$ with $n$ marked points and cone 
angles $\theta_i$ in a germ of AdS cone-manifold is a couple 
$(g,h)\in \cH_{g,n,AdS}(\theta_1, \cdots, \theta_n)$, where
\begin{itemize} 
\item $g$ is a smooth 
metric on a closed surface $\Sigma$ of genus $g$, with $n$ conical
singularities where the angles are $\theta_1, \cdots,\theta_n$,
\item $h$ is a symmetric bilinear form on $T\Sigma$ defined outside the
marked points,
\end{itemize}
such that: 
\begin{enumerate}
\item $h$ is traceless relative to $g$,
\item $h$ satisfies the Codazzi equation: $d^\nabla h=0$, with respect to $g$, 
\item the Gauss equation holds: $K_g=-1-\det_g(h)$,
\item $\det_g(h)$ is bounded.
\end{enumerate}
We call $\cH_{g,n,AdS}$ the universal cover of the space of maximal surfaces
in germs of AdS cone-manifolds, considered up to diffeomorphisms of the
underlying surface. 
\end{df}

There is a natural map:
$$ \ext: \cH_{g,n,AdS}(\theta_1, \cdots, \theta_n) \rightarrow
\cM_{g,n,AdS}(\theta_1, \cdots, \theta_n)~, $$
obtained by constructing the (unique) GHMC AdS manifold containing a given
maximal surface, i.e. ``extending'' the metric on the surface to a maximal
AdS metric. By construction, the image of $\ext$ is the space of GHMC AdS
manifold which contain a maximal surface. Question \ref{q:q1} above is
equivalent to asking whether $\ext$ is one-to-one. 

\subsection{Maximal surfaces from QHD}

The main point of this section is that, as in the non-singular case, maximal
surfaces in germs of AdS
cone-manifolds can be parametrized by the cotangent bundle of Teichm\"uller
space. 

\paragraph{Angles between $0$ and $\pi$.}

As in the hyperbolic case, one has to distinguish between the
cases $\theta\in [0,\pi]$ and $\theta\in[\pi,2\pi]$. We will first 
consider the case of small total angles.

\begin{thm} \label{tm:main-ads-sing}
Let $\theta_1, \cdots, \theta_n\in (0,\pi)$ be fixed. There is a natural
homeomorphism $\Psi_{g,n,AdS, \theta_1, \cdots, \theta_n}$ 
between $T^*\cT_{g,n}$ and $\cH_{g,n,AdS}(\theta_1, \cdots, 
\theta_n)$. The image of a point $(c,h)\in T^*\cT_{g,n}$ --- where $c\in
\cT_{g,n}$ and $q$ is a QHD with at most simple poles at the marked points ---
is the unique maximal surface $(g,h)$ such that $c$ is the conformal class of 
$g$ and $h=Re(q)$.  
\end{thm}

\begin{proof}
The key point concerning the relation between QHD and maximal surfaces remains
Lemma \ref{lm:hqd-ads}, which remains valid in the context of maximal surfaces
with conical singularities. Lemma \ref{lm:conf-ads} also remains true under
the condition that $g$ has conical singularities with angles in $(0,\pi)$ and
that the HQD $q$ has at most simple poles at the singular points of $g$.
Clearly equation (\ref{eq:*-ads}) still corresponds to the Gauss equation for
the conformal metric $g'$, and the fact that it has a unique solution can be
obtained as in the non-singular case, since the function $\det_g(h)$ goes to
$0$ at the cone singularities. 

Conversely, given a maximal surface in a germ of AdS cone-manifold, with
angles in $(0,\pi)$ at the singular points and principal curvatures going to
$0$ near the singular points, one can consider a ``reference'' flat
metric $g_0$ conformal to $g$, near one of the singular points, say $x$.
The metric $g$ then grows as $r^{(\theta-2\pi)/2\pi}$ near
$x$ (where $r$ is the distance to $x$ for $g_0$), so that $h$ is the real part
of a QHD with at most a simple pole at $x$. 

The outcome is that, as in the non-singular
case, the map from maximal surfaces to QHD with at most simple poles at the
singular points is a homeomorphism.
\end{proof}

\paragraph{Angles between $\pi$ and $2\pi$.}

It is also relevant, from a mathematical or a physical viewpoint, to consider
conical singularities with angles between $\pi$ and $2\pi$, so that the angle
defect at the singular lines is in $(0,\pi)$. It is easy to construct
examples of cone-manifolds with such cone angles using a warped product of
the corresponding 2-metric with the time axis. We would like to understand
to which extent such manifolds can be deformed. Our first observation is
as follows.

\begin{remark}
Let $\theta_1, \cdots, \theta_n\in (0, 2\pi)$. Let $(c,q)\in T^*\cT_{g,n}$ be
such that the QHD $q$ has no pole at $x_i$ for each $i$ such that $\theta_i\in
[\pi, 2\pi)$. Then the construction in Theorem \ref{tm:main-ads-sing} can
still be used, and provides us with an AdS cone-manifold with cone angles 
$\theta_1, \cdots, \theta_n$ along time-like line.
\end{remark}

This remark can be used for instance for data lying in $(c,q)\in T^*\cT_g$, since then
$q$ has no pole. However, it is intuitively clear that only 
``spacetimes with particles
with no momentum'' can be obtained this way. 

A probably more interesting description of the case of large total angles 
uses the notion of duality. Namely,
there is a notion of duality for surfaces in AdS. The duality interchanges the
first and third fundamental forms of the surface, so that the metric induced 
on the dual surface is precisely the third fundamental form of the original one.
This notion of duality appears to extend, to some
extent, to maximal surfaces in AdS cone-manifolds, and this is what we propose
to use to describe the case of large total angles. Our proposal is based
on the following lemma.

\begin{lemma} \label{lm:duality-ads}
Let $\theta_1, \cdots, \theta_n>0$, and let $(g,h)\in \cH_{g,n,AdS}(\theta_1,
\cdots, \theta_n)$ be a maximal surface in a germ of AdS cone-manifold. Let
$B$ be its shape operator, i.e. $B$ is self-adjoint for $g$ and
$h=g(B\cdot,\cdot)$, and let $g^*$ be its third fundamental form, so that
$g^*(x,y) := g(Bx,By)$. 
Then $(g^*,h)$ is a maximal surface in a germ of AdS cone-manifold, with
cone angles $2\pi-\theta_1, \cdots, 2\pi-\theta_n$ at the $x_i$. 
Moreover, $g$ is the third fundamental form of $(g^*,h)$. 
\end{lemma}

Note that the duality used here can not be understood in a global way, this is
quite apparent because the two dual surfaces do not have the same cone angles,
so that they do not ``live'' in the same cone-manifold. It should rather be
thought as a ``local'' construction, which by a local reconstruction yields a
maximal surface in ``another'' AdS cone-manifold. We call $(g^*,h)$ the
maximal surface dual to $(g,h)$.

\begin{proof}
By definition of a maximal surface in a germ of AdS cone-manifold, 
$B$ is a solution of Codazzi equation so that Proposition \ref{pr:metric-change} 
applies. So the Levi-Civit\`a connection of $g^*$ is:
$$ \nabla^*_xy = B^{-1}\nabla_x(By)~. $$
A simple computation then shows that $B^{-1}$ is a solution of the Codazzi
equation for $g^*$: $d^{\nabla^*}B^{-1}=0$. This can also be written as: 
$d^{\nabla^*}h=0$. 

Moreover, by the Gauss equation for $(g,h)$, $K=-1-\det(B)$, so that
$K^*=K/\det(B)=-K/(K+1)=-1-\det(B^{-1})$. Thus, $(g^*,h)$ is maximal surface in a
germ of AdS cone-manifold. The fact that conical singularities of $g$ of total angle
$\theta_i$ become conical singularities of $g^*$ of angles $2\pi -\theta_i$ is easily verified.
\end{proof}

Thus, as the above lemma shows, given
the data $(c,q)\in T^* \cT_{g,n}$ one can reconstruct the ``dual'' maximal surface
$(g^*,h)$ in a germ of AdS cone-manifold by solving the corresponding Gauss
equation. One gets the dual metric $g^*$ with conical angles equal to $2\pi-\theta_i$,
and with its corresponding principal curvatures going to zero at the singular
points. One then gets the corresponding AdS cone-manifold
by our previous results. Thus, what we propose is that it is
this ``dual'' AdS cone-manifold that should be used to describe particles
with large total angles. 

However there is a serious problem in this description, which is somewhat 
hidden in the way Lemma  \ref{lm:duality-ads} is stated: the metric $g^*$
is in general not smooth outside the singular points of $g$, since $B$
might be singular (i.e. non-invertible) at some points. More precisely,
$g^*$ has ramification points --- conical singularities with angle an
integer multiple of $2\pi$ --- at the zeros of $h$.

So at this stage the suggestion of using the duality is only an idea, because
it is not clear to us at the moment how the ``original'' 3-manifold can be
reconstructed from the ``dual'' one. We leave these interesting issues for future work.
 
Let us also remark that there is a possible analog of this ``particle duality'' construction for 
quasi-Fuchsian hyperbolic cone-manifolds with large cone angles. The main difference
with the AdS case, however, is that the dual surface in this case lies in a different space:
here we have the familiar hyperbolic - dS duality. However, as it will become clear from
section 6, in the dS setting it is natural to use not zero, but constant mean curvature
surfaces. Thus, it is not clear how the above ``particle duality'' applies to these two settings.
We leave this to future work.

\paragraph{Deformations of maximal surfaces in AdS cone-manifolds.}

Although we will not use it formally, it is interesting to remark that
maximal surfaces have no non-trivial deformation in GHMC AdS cone-manifolds. 
Another way to formulate it is that the map $\ext$ defined below has a nice
local behaviour: it is a local homeomorphism between the space of maximal
surfaces in germs of AdS cone-manifolds and the space of GHMC AdS
cone-manifolds (for fixed values of the angles in $(0,\pi)$).

\begin{lemma} \label{lm:defos-ads-cone}
Let $S$ be a closed, space-like maximal surface in a GHMC AdS cone-manifold
$M$. There is no first-order deformation of $S$ among maximal surfaces in $M$.
\end{lemma}

In other words, there is no first-order deformation of $S$ in $M$ such that
the mean curvature of $S$ remains $0$ (at first order) and $S$ remains
orthogonal to the singular lines. 

\begin{proof}[Sketch of the proof.] 
The second variation of the area, given by equation (\ref{sec-var-ads}), is 
still valid for maximal surfaces with conical singularities. Under a
first-order deformation such that the surface remains orthogonal to the
singular locus, $df$ goes to zero at the singular points, so that an
integration by part can be performed, showing that $A''$ is negative definite,
and the result follows.
\end{proof}

\subsection{Towards two parameterizations by two copies of Teichm\"uller space} 

We have seen in section 3 that, in the non-singular case, there are two
parameterization of the space of GHMC AdS manifolds by the product of two
copies of Teichm\"uller space, one by the construction due to Mess
\cite{mess}, translated in terms of the metrics $I^\#_\pm$ on any space-like
closed surface, the other using the metrics $I^*_\pm$ on the unique maximal
surface in the manifold considered.

It would be quite satisfactory to have a similar parameterization, when the
singular angles are fixed and at most equal to $\pi$, by
the copies of two products of the Teichm\"uller space with $n$ marked points. 
The results presented here are however limited; we only remark that,
given a GHMC AdS cone-manifold $M$, one can define from it a pair of hyperbolic
metrics $^\#_\pm$ with conical singularities, in a way 
which generalizes the Mess
construction which was already recalled in section 3 in the non-singular
case. Moreover the angles at the conical singularities of the hyperbolic
metrics are the same as the angles at the singular lines of $M$. 

Under the hypothesis that $M$ contains a maximal surface $S$, it is 
also possible
to associated to it two hyperbolic metrics with conical singularities in a way
which extends the definition of the metrics $I^*_\pm$ of section 3.
The angles at the singular points of those metrics are again the same as
the angles at the singular lines of $M$.

We do not know yet whether either of those construction determines a map from
the space of GHMC AdS cone-manifolds (with given angles at the singular lines)
and $\cT_{g,n}\times \cT_{g,n}$. It is interesting to remark, however, that
both those construction lead to a natural question concerning an extension of
Theorem \ref{tm:labourie} to the setting of hyperbolic metrics with conical
singularities, and that a positive answer to this question would imply that
the two maps to $\cT_{g,n}\times \cT_{g,n}$ constructed here are one-to-one
--- and also that question \ref{q:q1} above has a positive answer.

\paragraph{Principal curvatures of maximal surfaces.}

Before we can define the metrics $I^\#_\pm$ and $I^*_\pm$, we have to check
that the principal curvatures of maximal surfaces in AdS cone-manifolds are
bounded by $1$.

\begin{lemma} \label{lm:almost-ads-sing}
Let $\theta_1, \cdots, \theta_n\in (0,\pi)$, and let $(g,h)\in
\cH_{g,n,AdS}(\theta_1, \cdots, \theta_n)$ be a maximal surface in a germ of
AdS cone-manifold. Then the principal curvatures of $(g,h)$ are everywhere in
$(-1,1)$.   
\end{lemma}

\begin{proof}
The argument used in the proof of Lemma \ref{lm:curvatures-ads}, based on the
maximal principle for the function $k$ equal to the largest principal
curvature of $(g,h)$, can still be applied since $k\rightarrow 0$ near the
singular points. 
\end{proof}

\paragraph{Hyperbolic metrics from space-like surfaces.}

Let $M$ be a GHMC AdS cone-manifold, and let $S$ be a smooth space-like
surface in $M$ which is orthogonal to the singular lines, and such that 
its principal curvatures go to $0$ near the singular lines. We define two
metrics on $S$ as in section 3: given two vectors $x,y$ tangent to $S$ at the
same point: 
$$ I^\#_\pm(x,y) = I((E\pm JB)x, (E\pm JB)y)~, $$
where $B$ is the shape operator of $S$, $E$ is the identity, and $J$ is the
complex structure of $g$ on $S$. Lemma \ref{lm:hyper-metrics} still applies,
and shows that $I^\#_+$ and $I^\#_-$ are hyperbolic metrics. 

Moreover, since
$S$ is orthogonal to the singular lines in $M$ and its principal curvatures go
to $0$ there, it is quite clear that $I^\#_+$ and $I^\#_-$ have conical
singularities at the intersection of $S$ with the singular lines of $M$, with
angle equal to the angle of $M$ at those singular lines. 

\paragraph{The metrics $I^\#_\pm$ do not depend on the space-like surface.}

As in the non-singular case, the metrics $I^\#_\pm$ do not depend on the
choice of the space-like surface $S$ from which they are ``built''. It is
however not possible to use the same proof as in section 3. We will
instead describe explicitly how those metrics change under a deformation of
$S$ of the form $fN$, where $N$ is a unit orthogonal vector to $S$. To express
the result, we consider two tangent vector fields $x,y$ on $S$; we extend the
vector field $fN$ as a smooth vector field on $M$ in a neighbourhood of $S$,
and extend $x$ and $y$ as vector fields defined in a neighbourhood of $S$,
such that $[fN,x]=[fN,y]=0$. 

\begin{lemma} \label{lm:defo-s}
Let $V:=(E+JB)^{-1}JDf$, where $Df$ is the gradient of $f$ on $S$ for its
induced metric. Then, under the first-order deformation $fN$, the metric
$I^\#_+$ varies as: 
$$ \dot{I}^\#_+(x,y) = 2(\delta^*_{I^\#_+}V)(x,y)~, $$
where $\delta^*_{I^\#_+}V$ is defined in the usual way (see \cite{Be}) as:
$$ 2(\delta^*_{I^\#_+}V)(x,y) = I^\#_+(\nabla^{I^\#_+}_xV, y) +
I^\#_+(\nabla^{I^\#_+}_yV, x)~. $$
\end{lemma}
A similar statement applies to $I^\#_-$. 

\begin{proof}
The first-order variations of the induced metric and the shape operator of $S$
under a normal deformation are well-known, and can be found by an elementary
computation as is classical in the Riemannian context (see e.g. \cite{spivak}):
$$ \dot{I} = 2fI~, ~~ \dot{B} = \hess(f) - f(E+B^2)~. $$  
It then follows from the expression of the variation of $I$ that: 
$$ \dot{J} = f(JB-BJ)~. $$
An expression of the first-order variation of $I^\#_+$ follows: 
$$ \dot{I}^\#_+(x,y) = \dot{I}((E+JB)x(E+JB)y) 
+ I((\dot{J}B+J\dot{B})x, (E+JB)y)
+ I((E+JB)x, (\dot{J}B+J\dot{B})y) = $$
$$ 2f\II((E+JB)x,(E+JB)y) + f I((JB^2 - BJB -J -JB^2)x, (E+JB)y) 
+ fI((E+JB)x, (JB^2-BJB-J-JB^2)y) $$
$$ + I(J\nabla_xDf, (E+JB)y) + I((E+JB)x, J\nabla_yDf)~. $$
A tedious but straightforward computation then shows that most terms cancel
out, leaving that: 
\begin{eqnarray*}
\dot{I}^\#_+(x,y) & = & I(J\nabla_xDf, (E+JB)y) + I((E+JB)x, J\nabla_yDf) \\ 
& = & I(\nabla_x(JDf), (E+JB)y) + I((E+JB)x, \nabla_y(JDf)) \\
& = & I^\#_+((E+JB)^{-1}\nabla_x(JDf), y) + I^\#_+(x,
(E+JB)^{-1}\nabla_y(JDf))~.  
\end{eqnarray*}
But Proposition \ref{pr:metric-change} applies to this situation, and shows
that the Levi-Civit\`a connection $\nabla^{I^\#_+}$ of $I^\#_+$ is given by:
$ \nabla^{I^\#_+}_xy = (E+JB)^{-1}\nabla_x((E+JB)y)$. So we obtain that: 
$$ \dot{I}^\#_+(x,y) = I^\#_+(\nabla^{I^\#_+}_x((E+JB)^{-1}JDf), y) + I^\#_+(x,
\nabla^{I^\#_+}_y((E+JB)^{-1}JDf))~, $$
as needed.
\end{proof}

\begin{cor} 
The metric $I^\#_\pm$ do not depend on the choice of $S$. 
\end{cor}

\begin{proof}
It is sufficient to show that, under a small deformation of $S$, those metrics
only change by the action of a diffeomorphism. But the previous lemma
precisely shows that, under a normal deformation of $S$, the metric $I^\#_+$
varies as under the action of the vector field $V$, and the result
follows. The same proof can be used for $I^\#_-$.  
\end{proof}

It follows that, as in the non-singular case, there is a well-defined map
$\mess$ from the space of GHMC AdS manifolds, with angles $\theta_1, \cdots,
\theta_n$ at the singular lines, to the space of pairs of
hyperbolic metrics with conical singularities of angles  $\theta_1, \cdots,
\theta_n$ on $\Sigma$.

\paragraph{Hyperbolic metrics from maximal surfaces.}

Suppose now that $M$ contains a maximal surface $S$.
Lemma \ref{lm:almost-ads-sing} then allows us to define another pair of
hyperbolic metrics defined on $S$, as in section 3, by:
$$ I^*_\pm(x,y) = I((E\pm B)x, (E\pm B)y)~. $$
The argument used in section 3, based on Proposition \ref{pr:metric-change},
still shows that $I^*_\pm$ are hyperbolic metrics. By the definition of a
maximal surface used here, $S$ is orthogonal to the singular locus and its
principal curvatures go to $0$ there, and it follows quite directly that 
$I^*_+$ and $I^*_-$ have conical singularities at the intersections of $S$
with the singular lines in $M$, and that their angles at those singular points
are equal to the angles around the corresponding singular lines. 

It follows that, for fixed angles $\theta_1, \cdots, \theta_n$, 
there is a map $\max:\cM_{g,n,AdS}(\theta_1, \cdots, \theta_n)\rightarrow
\cT_{g,n}\times \cT_{g,n}$, sending a maximal surface with conical
singularities to a pair of points in the Teichm\"uller space with $n$ marked
points. 

\paragraph{Diffeomorphisms between hyperbolic metrics with conical
  singularities.}

The two constructions made above, of the hyperbolic metrics $I^\#_\pm$ and
$I^*_\pm$, lead naturally to ask whether the analog of Theorem
\ref{tm:labourie} holds for hyperbolic surfaces with conical singularities, in
particular when the cone angles are in $(0,\pi)$ (or perhaps even in
$(0,2\pi)$. 

\begin{question}\label{q:q3}
Let $S$ be a closed surfaces, with $n$ marked points $x_1, \cdots, x_n$, and
let $\theta_1, \cdots, \theta_n\in (0,\pi)$. Let $g,g'$ be two hyperbolic
metrics on $S$, with conical singularities at the $x_i$, with angle $\theta_i$
at $x_i$. Is there a unique bundle morphism $b:TS\rightarrow TS$, defined
outside the $x_i$, such that: 
\begin{itemize}
\item $b$ is self-adjoint for $g$, with positive eigenvalues,
\item $d^\nabla b=0$, where $\nabla$ is the connection of $g$,
\item $\det(b)=1$,
\item $g(b\cdot, b\cdot)$ is isotopic to $g'$ ?
\end{itemize}
\end{question}

It is quite conceivable that methods close to those used in \cite{L5} could
lead to a positive answer, the first case to check being when the angles are
equal to $0$ (i.e. for hyperbolic with cusps rather than conical
singularities). Considering this question here, however, would take us too far
from main idea of the paper.

As mentioned above, a positive answer to question \ref{q:q3} would imply a
positive answer to question \ref{q:q1}. Indeed the machinery developed in
section 3 for the metrics $I^\#_\pm$ already shows that, if $g_+$ and $g_-$
are two metrics on $\Sigma$ with conical singularities of angles $\theta_1,
\cdots, \theta_n$ at the points $x_1, \cdots, x_n$, and if $g_+=I^\#_+$ and
$g_-=I^\#_-$ for some GHMC AdS cone-manifold containing a maximal surface $S$,
then there exists a bundle morphism $b:T\Sigma\rightarrow T\Sigma$, defined
outside the $x_i$, with the property in the statement of question \ref{q:q3};
$b$ is constructed as: 
$$ b:= (E+JB)^{-1}(E-JB)~, $$
where $B$ is the shape operator of $S$. Conversely, if $g_+$ and $g_-$ 
are as in the statement of question \ref{q:q3} and a bundle morphism $b$
exists, then setting: 
$$ B := -J(E+b)^{-1}(E-b)~ , ~ ~ g:=g_+((E+JB)^{-1}\cdot, (E+JB)^{-1}\cdot) $$
yields the induced metric and shape operator of a maximal surface in a
GHMC AdS manifold for which $I^\#_+=g_+$ and $I^\#_-=g_-$.

Then, knowing that any GHMC AdS cone-manifold contains a unique maximal
surface (by the construction just outlined) it would also follow that 
the metrics $I^*_+$ and $I^*_-$ also provide another parameterization of the
space of GHMC AdS cone-manifolds, with fixed angles $\theta_1, \cdots,
\theta_n$, by the product of two copies of $\cT_{g,n}$.

\paragraph{Interpretations in Teichm\"uller theory.}

The maps defined above can be composed in various ways. In particular, for
each choice of angles $\theta_1, \cdots, \theta_n$, we obtain 
two maps from $T^*\cT_{g,n}$ to $\cT_{g,n}\times\cT_{g,n}$, one constructed
from the map $\mess:\cH_{g,n,AdS}(\theta_1, \cdots, \theta_n)\rightarrow
\cT_{g,n}\times\cT_{g,n}$, the other using the map
$\max:\cH_{g,n,AdS}(\theta_1, \cdots, \theta_n)\rightarrow 
\cT_{g,n}\times\cT_{g,n}$. Both those maps send the zero section in
$T^*\cT_{g,n}$ --- which corresponds to totally geodesic maximal surfaces ---
to the diagonal in $\cT_{g,n}\times\cT_{g,n}$. It would be interesting to know
whether those maps have an interpretation in terms of Teichm\"uller theory.

\section{Minkowski and de Sitter manifolds}

The previous sections have dealt with hyperbolic and AdS manifolds, possibly with
conical singularities. As is mentioned in the introduction, some of the results
described in this paper also apply, with appropriate modifications, to 3-manifolds modelled on 
Minkowski and de Sitter spaces. The purpose of this section is to collect together the
corresponding statements. We first consider de Sitter manifolds, for which things 
turn out to be quite similar to AdS manifolds. The Minkowski manifolds, where most differences occur,
is considered last. Most of the results presented in this section are more or less obvious
generalizations of the earlier statements. Thus, in most cases the proofs are
omitted or only sketched.

\subsection{De Sitter manifolds}

\paragraph{Constant Gauss curvature foliations and duality.}

In this subsection we consider globally hyperbolic maximally compact de Sitter
manifolds. By definition, these are manifolds (i) locally modelled on the de Sitter space, denoted
here by $dS^3$; (ii) containing a closed, embedded space-like surface which
intersect each time-like curve exactly once; (iii) maximal under the above
conditions. 

As is recalled in the introduction (see also section 2), there is a natural duality 
between quasi-Fuchsian hyperbolic manifolds and (some) 
pairs of de Sitter GHMC manifolds. To each
quasi-Fuchsian hyperbolic manifold $M$ 
the duality associates two connected GHMC
dS manifolds, each corresponding to one connected component of the space
of totally geodesic planes in $M$ which do not intersect its convex core. We
will call those dS manifolds $M^*_+$ and $M^*_-$, respectively.

Labourie \cite{L6} proved that each end of a quasi-Fuchsian manifold ---
i.e. each connected component of the complement of the convex core --- has a
unique foliation by constant Gauss curvature surfaces. Those surfaces
are convex, so that their tangent planes do not intersect the convex core of
$M$. Moreover, the Gauss curvature of those surfaces is increasing as one
moves away from the convex core, and it varies from $-1$ --- near the boundary
of the convex core --- to $0$. We call $k(p), p\in M$ the Gauss curvature of 
the constant Gauss curvature surface $S_p$ at a point $p\in M$. Thus, $k$ is a smooth function 
defined on the complement of the convex core of $M$.

By definition of the dual of $M$, each of those tangent planes corresponds to
a point in $M^*_+$ (for surfaces in the ``upper'' end, which we can call
$M_+$, of $M$) or of $M^*_-$
(for surfaces in the other end of $M$, called $M_-$). 
So each surface $S$ in the complement of the convex core of $M$ has a dual
surface $S^*$, contained either in $M^*_+$ or in $M^*_-$. 
Moreover, a well-known property of the
hyperbolic-de Sitter duality (see e.g. \cite{Ri,RH,shu}) is that
$S^*$ is space-like and convex, and that its curvature is equal to $K/(K+1)$,
where $K$ is the curvature of $S$ 
at the corresponding point. Therefore, if $S$ has constant Gauss curvature, 
so does $S^*$.

Let $x\in M^*_+$. Then $x$ is dual to a totally geodesic plane $X^*\subset M$
contained in $M_+$ and not intersecting the 
convex core of $M$. Let $p\in X^*$ be a critical point of the function $k$ introduced above.
At this point the constant Gauss curvature surface $S_p$ through $p$ is tangent to
$X^*$. Therefore the point $x$ lies in the plane $S_p^*$ dual to $S_p$. 
Moreover, since the constant Gauss curvature surfaces 
are convex, any critical point of the restriction of $k$ to $X^*$
corresponds to a local maximum. But as one moves towards infinity on $x^*$ the
distance to the convex core goes to $\infty$, so that $k\rightarrow 0$. It
follows that the restriction of $k$ to $x^*$ has exactly one critical point,
so that $x$ is contained in the dual of exactly one constant Gauss curvature
in $M_+$, and therefore the duals of the constant Gauss curvature of $M_+$ are
the leaves of a foliation of $M^*_+$ by constant Gauss curvature
surfaces. 

This simple duality argument shows that GHMC dS manifolds that are dual to
hyperbolic ends of hyperbolic 3-manifolds are foliated by convex, 
constant Gauss curvature surfaces (cf also \cite{barbot-zeghib}).  
The curvature varies between $-\infty$ --- 
near the boundary corresponding to the dual
of the convex core --- and $0$ --- near the
boundary at infinity. 

Now, a simple argument can be given to shown that all GHMC dS manifolds are
dual to hyperbolic ends of 3-manifolds modelled on $H^3$. The argument is
roughly as follows. Given a GHMC dS manifold, say $M_+$ one gets a set of data that
can be used to reconstruct a hyperbolic end. This set of data is as follows.
One piece of the data is the conformal structure of the boundary of
$M_+$. This conformal structure must then be the same as the one of
the hyperbolic end in question. The second piece of data is what
can be called a ``constant'' part of the metric on $M_+$ near infinity.
Indeed, the lemma \ref{lm:27} of Section 2 can be
easily adapted to the dS context. It is then easy to show that the metric on any
GHMC dS manifold, sufficiently close to the conformal boundary, can be
written as:
\begin{equation}
ds^2 = -dt^2 + e^{2t} g + 2h + e^{-2t} h g^{-1} h.
\end{equation}
Here $g$ is a metric in the conformal class of the boundary of 
$M_+$ that corresponds to $t\to\infty$, and $h$ is what is
referred to as the ``constant'' term in the metric. The data
$(g,h)$ are not arbitrary and are satisfy the Gauss and Codazzi
equations, but we will not need to spell out all the details
for the purposes of the present paper. The metrics of the
above form can be referred to as those of foliations ``equidistant to infinity'',
in contrast to the foliations equidistant to constant mean curvature
surfaces considered in the present paper. The present authors 
are planning to describe such foliations ``from infinity'' in
a future publication. The reader is referred to it for more details.

To finish the argument, given a GHMC dS manifold, one gets a set
of data $(g,h)$. The $g$ part of the data is of course not uniquely determined;
any metric in the conformal class of the boundary will do. One can then
take this data and use them to reconstruct a part of the hyperbolic
end in question. The metric in (a part of) the hyperbolic end is
obtained as:
\begin{equation}
ds^2 = dr^2 + e^{2r} g + 2h + e^{-2r} h g^{-1} h,
\end{equation}
where the same data $(g,h)$ are to be used. The metric above does
not in general cover the whole of the hyperbolic end, but it
does cover a ``large enough'' portion of it so that the whole
end can be reconstructed. The argument presented here is not
a proof, and is presented for motivation purposes only. We hope
to give more details in a future publication.

\paragraph{CMC foliations of de Sitter manifolds.}

In addition to the constant Gauss curvature surfaces mentioned above, the GHMC dS
manifolds also have a unique foliation by constant mean curvature surfaces,
see \cite{BBZ-cras,barbot-zeghib}. 
The mean curvature of those surfaces varies between $1$ 
(near infinity) and $\infty$ (near the boundary inside the manifold);
in particular it contains no maximal surface. The existence of a
foliation by maximal surfaces, with monotonous mean curvature, makes it
possible to use the maximal principle to show that, for any GHMC dS manifold
$M$ and for any $H\in (1,\infty)$, there is a unique closed, space-like,
embedded surface in $M$ with mean curvature $H$. We shall refer to this
surface as the unique CMC-H surface.

\paragraph{Construction from $T^*\cT_g$.}

As in the AdS case, there is a natural parameterization of the GHMC dS
manifolds of genus $g$ by the cotangent bundle of Teichm\"uller space. It works in
the same way as in the AdS case except that one has to consider constant mean
curvature rather than maximal surfaces. We call $\cM_{g,dS}$ the space of GHMC
dS manifolds of genus $g$.

\begin{lemma}\label{lm:ds-cotangent}
For each $H\in (1,\infty)$, there is a natural map
$\psi_{g,dS}^H:\cM_{g,dS}^H\rightarrow T^*\cT_g$. Namely, each $M\subset \cM_{g,dS}^H$
contains a unique CMC-$H$ surface $S$, its image is the point $(c,q)\in
T^*\cT_g$ such that $c$ is the conformal class of the induced metric of $S$
and the real part of $q$ is its second fundamental form. Moreover, the above
map $\psi_{g,dS}^H$ is one-to-one.
\end{lemma}

\begin{proof}[Sketch of the proof.]
Following previously used notations, we can call $\cH_{g,dS}^H$ the space of
``CMC-$H$ surfaces in germs of dS manifolds'', i.e. the space of couples
$(g,h)$, where $g$ is a metric on $\Sigma$, $h$ is a bilinear symmetric form
on $T\Sigma$, and $g$ and $h$ satisfy the Gauss and Codazzi equation for
CMC-$H$ surfaces: $d^\nabla h=0$, $\tr_g(h)=2H$, and $K_g=1-{\rm det}_g(h)$, where
$\nabla$ and $K_g$ are the Levi-Civit\`a connection and curvature of $g$
correspondingly.

But, given $g$ and $h$, we can set $h_0=h-Hg$. Then:
\begin{enumerate}
\item $\tr_g(h)=2H$ if and only if $h_0=Re(q)$, where $q$ is a bilinear form
  on the complexified bundle of $T\Sigma$,
\item if (1) is satisfied, then $d^\nabla h=0$ if and only if $q$ is a HQD,
\item if (1) and (2) are satisfied and $g'=e^{2u}g$, then the Gauss equation
  is satisfied by $(g',h)$ if and only if: 
\begin{equation} \label{eq:*-ds}
\Delta u = (1-H^2) e^{2u} - K_g - e^{-2u} det_g(h_0)~. 
\end{equation}
\end{enumerate}
We see that, given a CMC-$H$ surface in $M$ we get a point in $T^*\cT_g$.

Let us now show that the converse is also true. To this end, 
let us show that there is a unique choice of $u$ such that $(e^{2u}g,h)$ satisfies the Gauss equation.
The argument used in the proof of Lemma \ref{lm:conf-ads} can be applied 
here as well. Namely, let us introduce a functional:
\begin{equation}
F(u) = \frac{1}{2}\int_\Sigma \left( \| \nabla u\|^2 + (H^2-1)e^{2u} + 2K_gu + k^2 e^{-2u}
   \right) da~.
\end{equation}
As is easy to see, this functional is constructed so that the equation that follows from it under
the variation with respect to $u$ is exactly (\ref{eq:*-ds}).
Since $H>1$ everywhere in $M$, the above functional is a sum of four strictly
convex functionals. Therefore, if there is a solution to equation (\ref{eq:*-ds}),
it must be unique. The fact that (\ref{eq:*-ds}) actually has a solution
can be shown by methods similar to those used in the proof of Lemma \ref{lm:conf-ads}.

Thus, each element of
$T^*\cT_g$ corresponds to a unique CMC-$H$ surface in a germ of dS
manifold. Moreover, since each GHMC dS manifold contains a unique
CMC-$H$ surface, we get the result that each element of $T^*\cT_g$
corresponds to a unique GHMC dS manifold $M$. An explicit metric on 
(a portion of) this $M$ will be written down below.
\end{proof}

Thus, in spite of the fact that there is a duality between the hyperbolic and 
de Sitter settings, the things in de Sitter case work very much like they
do in the AdS case. Namely, one obtains all of the cotangent bundle 
over Teichm\"uller space as the image of the map $\psi$, and this map
is invertible. The main difference with the AdS situation is that now
there is not one, but infinitely many such maps, parameterized by
$H\in(1,\infty)$. 

\paragraph{Construction from $\cT_g\times \cT_g$.}

As we just saw, GHMC dS manifolds behave exactly as GHMC AdS manifolds
as far at the parameterization by $T^* \cT_g$ is concerned. This
analogy extends much further, as we shall presently see.
As we saw in the previous section, there are two possible ways
to parameterize the space $\cH_{g,AdS}$ of GMHC AdS manifolds by
two copies of $\cT_g$; one is by the pair of hyperbolic metrics
called $I^\#_\pm$ associated to any space-like surface, the other by the pair
of hyperbolic metrics $I^*_\pm$ associated to the unique maximal surface in a
GHMC AdS manifold.

In the dS case, there appears to be no analog of the pair of 
metrics $I^\#_\pm$. However, there is a natural analog of the  
metrics $I^*_\pm$. The only complication that arises is that, since it is not possible 
to use maximal surfaces, the construction will depend on a number $H\in (1, \infty)$. 
Thus, given such a number, we can consider the unique CMC-$H$ surface in a given GHMC dS manifold,
and obtain two metrics $I^*_\pm$, and thus two points in $\cT_g$. 
Conversely, given
a number $H\in(1,\infty)$ and a pair of points in $\cT_g$, we can reconstruct the manifold $M$.

Before we define the metrics $I^*_\pm$, we need a lemma 
analogous to Lemma \ref{lm:curvatures-ads}. 

\begin{lemma} \label{lm:curvatures-ds}
Let $S$ be a closed CMC-$H$ surface in a dS manifold, and let $B$ be its shape
operator. The eigenvalues of $B_0 = B-HE$ have absolute values bounded by
$\sqrt{H^2-1}$. 
\end{lemma}

\begin{proof}[Sketch of the proof.]
The proof is basically the same as the proof of Lemma
\ref{lm:curvatures-ads}. We call $k$ the non-negative eigenvalue of $B-HE$. 
We set $\chi:=\log(k)/2$. Then the Codazzi equation shows that $\Delta\chi = -K$. 
But, according to the Gauss equation:
$$ K =1 - \det(B) = 1-\det(HE+B_0) = 1-H^2 -\det(B_0) = 1-H^2+k^2~, $$
so that $\Delta\chi = H^2-1 - k^2$. At the maximum of $\chi$, $\Delta\chi\geq
0$, so that $k\leq \sqrt{H^2-1}$. The same argument can be used for the
smallest eigenvalue of $B-HE$, showing that it is larger than $-\sqrt{H^2-1}$.
\end{proof}

It is therefore possible to define the metrics $I^*_\pm$ analogously to
the AdS case. The metrics obtained will be smooth as is guaranteed by the above lemma.
We shall keep the same notation for these two metrics, hoping that it will
not lead to any confusion. 

\begin{df}
Let $H\in (1,\infty)$. For any GHMC dS manifold $M$, let $S$ be the unique
closed, space-like, embedded CMC-$H$ surface in $M$. We define the metrics
$I^*_\pm$ as: 
$$ I^*_\pm = I((\sqrt{H^2-1}E\pm B_0)\cdot, (\sqrt{H^2-1}E\pm B_0)\cdot)~, $$
where $B_0$ is the traceless part of the shape operator of $S$ and $I$ is its
induced metric.
\end{df}

Proposition \ref{pr:metric-change} can be used in this setting, it shows that
the Levi-Civit\`a connections of $I^*_+$ and $I^*_-$ are given by: 
$$ \nabla^\pm_xy = (\sqrt{H^2-1}E\pm B_0)^{-1}\nabla_x(((\sqrt{H^2-1}E\pm
B_0)y)~, $$ 
and that their curvatures are:
$$ K^\pm = \frac{K}{\det(\sqrt{H^2-1}E- B_0)} =
\frac{1-\det(HE+B_0)}{H^2-1+\det(B_0)} = \frac{1-H^2 -
\det(B_0)}{H^2-1+\det(B_0)} = -1~. $$

It is then possible to define: 
$$ b := (\sqrt{H^2-1}E+B_0)^{-1}(\sqrt{H^2-1}E - B_0)~, $$
then the same argument as those used in the proof of Theorem \ref{lm:1} show
that: 
\begin{enumerate}
\item $b$ is self-adjoint for $I^*_+$,
\item $I^*_- = I^*_+(b\cdot, b\cdot)$,
\item $\det(b)=1$,
\item $d^{\nabla^+}b=0$.
\end{enumerate}
Conversely, Theorem \ref{tm:labourie} shows that $b$ is uniquely determined by
the pair $I^*_+$ and by $I^*_-$. One can then set: 
$$ B:=\sqrt{H^2-1}(E-b)(E+b)^{-1}~, ~~ g = \frac{1}{4(H^2-1)} g_+((E+b)\cdot,
(E+b)\cdot)~.$$
It is then easy to show that $g$ and $B$ are the induced metric and
shape operator of a CMC-$H$ surface in a dS manifold. Then, for this surface,
the metrics $I^*_\pm$ coincide with the ones we started with, 
and this shows that the 
map sending a GHMC dS manifold to the metrics $I^*_+$ and 
$I^*_-$ on its CMC-$H$ surface is both injective and surjective. 

To summarize, we have two parameterization of the space $\cM_{g,dS}$ of
GHMC dS manifolds. Both parameterizations depend on a number
$H\in(1,\infty)$. The first parameterization is by the cotangent bundle
over the Teichm\"uller space. The second parameterization is by two
points in Teichm\"uller space. In both cases the relevant maps are one-to-one.
One can combine these two maps and obtain yet another map
$T^* \cT_g \to \cT_g \times \cT_g$. It would be of great interest
to clarify the significance of all these maps for Teichm\"uller theory.

\subsection{Equidistant foliations}

In this subsection we would like to write down explicitly the metric on a GHMC dS manifold
that results from the above described parameterization by $T^* \cT_g$. This metric
is of the same general type already encountered by us in various parts of this paper,
see e.g. (\ref{metric-1}). Namely, it is the one that arises by considering
an equidistant foliation. The only thing that changes in the dS context is the
signature. The resulting metric is given by:
\begin{equation}\label{metric-ds}
ds^2 = -dt^2 + I((\cosh(t) E + \sinh(t) B)\cdot,(\cosh(t) E + \sinh(t) B)\cdot).
\end{equation}
This metric describes a foliation of a dS manifold by surfaces equidistant to
a given one $(t=0)$. Let us choose this surface to be the one of constant 
mean curvature $H$. Let us rewrite the resulting metric in the complex-analytic
form, similar to (\ref{fock}). We get:
\begin{equation}
ds^2 = -dt^2 + e^\varphi | (\cosh(t)+ H\sinh(t)) dz + \sinh(t) e^{-\varphi} \bar{t} d\bar{z}|^2.
\ee  
An important difference with (\ref{metric-1}) is the appearance of $H$ in the term
proportional to $dz$, and, of course, the signature.
The ``Liouville'' field $\varphi$ satisfies its Gauss equation, which takes the following form:
\be
2 \partial_{z\bar{z}} \varphi = (H^2-1) e^\varphi - e^{-\varphi} t\bar{t}.
\ee

The metric (\ref{metric-ds}) covers the part of the space between the CMC-$H$ surface
$t=0$ and infinity completely. It becomes singular for $t=-T$, where
\begin{equation}
\frac{\cosh{T}}{\sinh{T}} = (k_{max}+H).
\end{equation}
Thus, as in all previous cases, the equidistant foliation typically does not cover
all of the space. It would be of interest to characterize the
position of the singularity of the foliation with respect to the
inner boundary of $M$. 

An alternative explicit expression for the metric on $M$ can be obtained
from the parameterization by two copies of the Teichm\"uller space. This is
a straightforward exercise which is left to the reader.

\subsection{De Sitter cone-manifolds}

We have seen that the results for smooth de Sitter manifolds are pretty similar to 
those for Anti de Sitter ones. In this subsection we generalize things to 
the situation when point particles are present.

\paragraph{Definitions.}

The notion of a dS cone-manifold can be defined similarly to what was done in the 
AdS cone-manifold case in section 5. GHMC AdS cone-manifolds are defined in the same manner.
As in the AdS case, most of the results of the non-singular case
can be generalized to the cone-manifolds. The main difference
with the AdS case is that now everything depends on a choice of
a number $H\in (1,\infty)$. As in the cases considered previously, 
one should treat differently the situation when all the total angles are less than $\pi$ 
and the one when the total angles are $\theta\in[\pi,2\pi)$.

\paragraph{What about canonical foliations ?}

The analogy with the non-singular case suggests that any GHMC dS
cone-manifold should have a foliation by constant
Gauss curvature surfaces. Such statement, if exists, may require the
restriction for all angles to lie in $[0,\pi]$.

\begin{question}
Let $M$ be a GMHC dS cone-manifold, with angles in $(0,\pi)$. Does $M$ have a
foliation by constant Gauss curvature surfaces orthogonal to the
singular lines ? Is it unique ?
\end{question}

The existence of such a foliation would follow from the existence of a similar
foliation for hyperbolic ends with conical singularities along lines, under
the same angle 
conditions. The argument used to go from the hyperbolic to the dS statement
uses the duality in a ``local'' way, since the presence of the singularities
makes it impossible to apply globally the hyperbolic-de Sitter duality. Given a
constant Gauss curvature foliation of a hyperbolic end with conical
singularities, one can consider the dual surfaces and use them to ``build'' a
GHMC dS cone-manifold.

In a similar way, and as in the AdS case, it is natural to wonder whether GHMC
dS cone-manifold have a CMC foliation. Again, one may require the total
angles to be less than $\pi$.

\begin{question}
Let $M$ be a GMHC dS cone-manifold, with angles in $(0,\pi)$. Does $M$ have a
foliation by CMC surfaces orthogonal to the singular lines ?
\end{question}

Given the existence of such a foliation, its uniqueness would follow from an
elementary argument based on the maximum principle, as already used above.

\paragraph{Description from $T^*\cT_{g,n}$.}

The construction of GHMC dS manifolds as parameterized by points in the cotangent 
bundle of the Teichm\"uller space still works in the dS setting. As in
the AdS case, we do not know whether all GHMC dS cone-manifolds can be
constructed in this way --- we only know that the construction provides a
parameterization by $T^*\cT_{g,n}$ of the space of GHMC dS cone-manifolds
containing a CMC-$H$ surface. Moreover the construction depends on the choice
of $H\in (1,\infty)$. 

As in other settings, given $n$ angles $\theta_1, \cdots, \theta_n\in
[0,2\pi)$, we denote by $\cM_{g,dS,n}(\theta_1, \cdots, \theta_n)$ the space of
dS metrics on $\Sigma \times \R$, with conical singularities along the lines
$\{ x_i\}\times\R$ (up to isotopies) 
for given points $x_1, \cdots, x_n\subset \Sigma$. 

Given such a GHMC dS manifold $M$, and a CMC-$H$ surface $S\subset M$ which is
orthogonal to the singular lines, one can consider the
conformal structure $c$ of the induced metric $I$ on $S$, as well as the QHD 
$q$ such that $\II-HI$ is the real part of $q$. Then the same
arguments as those used in sections 4, 5 show that $q$ has at most simple poles at
the singular points of $S$, so that $(c,q)$ corresponds to a point in
$T^*\cT_{g,n}$. 

Conversely, given a point $(c,q)\in T^*\cT_{g,n}$, one can choose a metric $g$
in the conformal class of $c$, with conical singularities of angles $\theta_1,
\cdots,\theta_n$ at the points $x_1, \cdots, x_n$, and set $h_0:=Re(q)$. Then
$h_0$ is traceless and satisfies the Codazzi equation with respect to
$g$. Considering another metric of the form $g'=e^{2u}g$ shows that $g'$ is
the induced metric and $h_0+Hg'$ is the second fundamental form 
if and only if equation (\ref{eq:*-ds}) is satisfied. But, as we have already seen in the
smooth case, this equation has a unique solution. We thus 
obtain a parameterization by $T^*\cT_{g,n}$ of those GHMC dS manifold
containing a CMC-$H$ surface orthogonal to the singular lines.

Everything said above applies in general. However, this description is essentially
useless when the total angles are large. Indeed, in that case the principal
curvatures, i.e., eigenvalues of $B_0= B - E H$ blow up at 
the singular points. To avoid this, we shall treat the case when all
angle as in $[0,\pi]$ and the case when angles are in $[\pi,2\pi)$ differently.
Let us consider the first case first. In this case, the principal curvatures
go to zero at the singular points. If this is the case, 
Lemma \ref{lm:curvatures-ds}
applies, and guarantees that the principal curvatures are bounded by
$\sqrt{H^2-1}$ everywhere on $S$. This then allows to reconstruct (a large portion of) the
manifold $M$ by considering e.g. the metric (\ref{metric-ds}). 

In the other case of large total angles, the principal curvatures of $B_0=B- E H$
diverge at the singularities. However, one can introduce another shape operator,
namely the one associated to the third fundamental form. Thus, let us as before
define: $B'=\III^{-1} \II$ and $B'_0=B'-E H'$, where $H'=(1/2){\rm Tr}(B')$. The ``principal 
curvatures'', i.e. the eigenvalues of $B'_0$ tend to zero at the singular points. 
One can then repeat all of the story above, but with the third fundamental form and
the ``shape operator'' $B'$ used instead of $I$ and $B$. As in the case of 
small total angles, these data are
parameterized by the cotangent bundle over the Teichm\"uller space, and this
parameterization map is one-to-one. Given a point in $T^* \cT_{g,n}$ one 
can reconstruct the data $\III, \II$. As in the AdS case of section 5,
the question is whether such a ``dual'' description is any useful. We will
not pursue this further in the present paper.

\paragraph{Construction from $\cT_{g,n}\times \cT_{g,n}$.}

It should not come as a surprise to the reader that the construction made in
the non-singular case extends, partially, to provide a map from the space
of ``CMC-$H$ surfaces in germs of dS cone-manifolds'' --- defined as in the
AdS case --- to the product of two copies of $\cT_{g,n}$. Provided that each
GHMC dS manifold (possible with a condition on angles) contains a unique CMC-$H$
surface, this would provide a map from the space of such GHMC dS
cone-manifolds to $\cT_{g,n}\times \cT_{g,n}$. This map is one-to-one if
question \ref{q:q3} has answer in the affirmative. Indeed, when the
total angles are less than $\pi$ the principal curvatures are bounded and
$I^*_\pm$ are hyperbolic metrics with conidal singularities
that can be used to reconstruct $M$. 

\subsection{Globally hyperbolic Minkowski manifolds}

\paragraph{First properties.}

We finish considerations of this paper by considering GHMC Minkowski manifolds.
These are manifolds locally modelled on the Minkowski 
3-dimensional space, which are
(i) globally hyperbolic; (ii) contain a closed space-like surface; 
(iii) maximal. 
They have a unique foliation by CMC-$H$ surfaces, 
with the mean curvature varying between
$0$ and $\infty$, see \cite{barbot-zeghib}.

\paragraph{Construction from $T^*\cT_g$.}

As in other cases considered in this paper, GHMC Minkowski manifolds 
are parametrized by $T^*\cT_g$. The construction in this case is well-known
and goes back to the paper by Moncrief \cite{Moncrief}. Thus, we will be
much more brief than in the previous cases. The construction is based on the use of CMC
surfaces --- the parameterization itself depending on the choice of a mean
curvature $H>0$. What makes it possible is the fact that
each GHMC Minkowski manifold contains a unique CMC-$H$ surface, for any $H>0$. 

As before, we will introduce the notion of a ``CMC-$H$ surface in a germ
of Minkowski manifold''. The fact that
each GHMC Minkowski manifold contains a unique CMC-$H$ surface means that
a CMC-$H$ surface in a germ of $M$ actually determines the whole $M$.
Thus, we are led to the problem of understanding couples $(g,h)$, 
where $g$ is a Riemannian metric on $\Sigma$ and $h$ is a bilinear 
symmetric 2-form on $T\Sigma$, such that: 
\begin{itemize}
\item $h_0=h-Hg$ is traceless with respect to $g$,
\item $d^\nabla h_0=0$,
\item $h$ satisfies the Gauss equation for surfaces in $R^2_1$:
$$ K_g = -{\rm det}_g(h) = -H^2 -{\rm det}_g(h_0)~. $$
\end{itemize}

Given such a couple $(g,h)$, one gets a point $(c,q)\in T^*\cT_g$. Indeed, 
$c\in \cT_g$ is the conformal structure of $g$, while $q$ is the QHD such that $h_0=Re(q)$. 
Conversely, it is not hard to show that any $(c,q)\in
T^*\cT_g$ corresponds to a unique couple $(g,h)$. Indeed,
as before, the only thing to consider is the Gauss equation. It
can be solved by finding a conformal factor $u$ that satisfies:
\begin{equation} \label{eq:*-minko}
\Delta u = -H^2 e^{2u} - K_g - e^{-2u} det_g(h_0)~. 
\end{equation}
This equation has a unique solution, reasons being the same 
as in the dS and AdS cases. 

\paragraph{One hyperbolic metric --- only.}

As we have seen in the previous cases, one can parameterize
the 3-manifolds by two copies of the Teichm\"uller space. In contrast, 
the ²Minkowski setting gives rise to a single
hyperbolic metric, not two such metrics. In this sense, the flat setting
appears as the ``degenerate'' version of the AdS and dS ones.

The natural hyperbolic metric that can be associated to a GHMC Minkowski
manifold $M$ is simply the third fundamental form of a surface in $M$. 
As in the AdS case with its two metrics $I^\#_\pm$, this
metric is independent of the space-like surface which is chosen. 
Note, however, that, in the AdS case, it was possible to recover the holonomy 
of the manifold from the two hyperbolic metrics $I^\#_+$ and $I^\#_-$. 
The hyperbolic metric of the Minkowski case --- the third
fundamental form --- only determines the linear part of the holonomy.

Indeed, recall that the group of orientation and time-orientation preserving
isometries of the Minkowski space is $\isom(\R^{2,1})=PSL(2,R) \rtimes R^3$, the semi-direct
product of the subgroup of isometries fixing the origin by the translations. 
Given a space-like surface $S$ in $\R^{2,1}$, its third fundamental form is the
pull-back of the metric on $H^2$ by the ``Gauss map'' $G:S\rightarrow H^2$
sending a point $x\in S$
to the unit normal to $S$ at $x$, considered as a point in $H^2$. If $S$ is
the image of the universal cover $\tilde{S'}$ of a closed surface $S'$ under a map
$\phi:\tilde{S'}\rightarrow \R^{2,1}$, 
equivariant under an action of $\rho:\pi_1(S')\rightarrow \isom(\R^{2,1})$, 
then
$G\circ \phi:\tilde{S'}\rightarrow H^2$ is equivariant under the linear part
$\rho_L:\pi_1(S')\rightarrow PSL(2,\R)$ of $\rho$. It is thus quite clear that
the third fundamental form of a space-like surface in a GHMC Minkowski
manifold uniquely determines the linear part of its holonomy, and conversely.

\subsection{Minkowski cone-manifolds}

\paragraph{Definitions.}

The definition of a Minkowski cone-manifold follows the same lines as the
corresponding definitions in the AdS and dS case, so we do not repeat it
here.

\paragraph{Parameterization by the cotangent bundle of Teichm\"uller space.}

By essentially the same reasons as in the hyperbolic case, a 
CMC-$H$ surface $S$ in a Minkowski cone-manifold must be orthogonal to the
singularity locus. One then considers the first $g$ and second $h$
fundamental forms of $S$, and introduces $h_0 = h- Hg$. These data
determine a point $(c,q)\in T^*\cT_{g,n}$, where $n$ is the number
of conical singularities. Indeed, the induced metric determines a
point $c\in \cT_{g,n}$. The quantity $h_0$ is the real part of a
holomorphic quadratic differential $t$ by Codazzi equation. Moreover,
the requirement that the right hand side of the Gauss equation
is integrable over the surface implies that $t$ has at most simple poles
at the singular points. Thus, $(c,t)\in T^*\cT_{g,n}$.

Conversely, each point in $T^*\cT_g$ is obtained in a unique way from a couple $(g,h)$
corresponding to the induced metric and second fundamental form of a ``CMC-$H$
surface in a germ of Minkowski cone-manifold''. The proof of this is based on the fact
that equation (\ref{eq:*-minko}) has exactly one solution.

As in the previous cases, there are two situations
to consider: when the total angles are less than $\pi$, and the
case of total angles $\theta\in[\pi,2\pi)$. Let us consider the first
case first. In this case the ``principal curvatures'', i.e.
the eigenvalues of $B_0 = B-HE$ tend to zero as one approaches
the singular points.  This is the ``good'' case, when the
description by $T^* \cT_{g,n}$ works. When there are 
conical singularities with 
total angles greater than $\pi$, the principal curvatures diverge unless 
the second fundamental form is a real part of a meromorphic
quadratic differential with no pole at the corresponding point, and
the above given description is less interesting. To remedy the
situation one may, as before, consider some dual description
instead. We will not pursue this issue further in the present paper.

\paragraph{The third fundamental form.}

As for non-singular Minkowski manifolds, one can consider the third
fundamental form of a space-like surface in a GHMC Minkowski cone-manifold
$M$. This form does not depend on the space-like surface chosen. However, as in the
non-singular case, it is only possible to recover the linear
part of the holonomy from $\III$.

\section{Hamiltonian Formulation}

In this section we show that the fact that the moduli space of constant
curvature 3-manifolds  
can often be parameterized by the cotangent bundle over the Teichm\"uller space
is complemented 
by the fact that the symplectic form on this moduli space that is induced by
the gravity action 
is also that of $T^* \cT_\Sigma$. We do this in the setting of globally
hyperbolic AdS manifolds, 
where this result is of direct relevance to physics. Other settings can be
treated analogously 
and will be commented upon below.

\subsection{Hamiltonian formulation of gravity in AdS setting}

Thus, we consider a GHMC AdS manifolds $M$. As we have discussed in Section 3,
each such manifold 
contains a unique maximal surface $S$. Moreover, it can be foliated by
surfaces equidistant to 
$S$, at least in the neighbourhood of $S$. The Hamiltonian formulation we give
is based on this equidistant foliation. 

The gravity action is given by:
$$ S = \int_M (R- 2\Lambda) dv_g, $$
where $R$ is the scalar curvature (trace of the Ricci tensor) of $g$,
$\Lambda$ is the cosmological constant, and $dv_g$ is the volume element of
$g$.
In our considerations above we have worked in ``units'' in which $|\Lambda|=1$,
and we shall continue to do so here. To arrive at the Hamiltonian formulation
one chooses a foliation of $M$ by spatial hypersurfaces. In our case these
will be the surfaces equidistant to the maximal surface $S$. One then chooses
a vector field that will play the role of time. We will use the vector field
normal to the equidistant surfaces. Thus, in particular, the shift function is
equal to zero. 

The conjugate variables of the Hamiltonian formulation are the metric $I_t$
induced on the surfaces $S_t$, and 
the momentum $\pi_t$. The action in the Hamiltonian formulation is given by:
$$ S = \int_\R \int_{S_t} \left(\langle \pi_t, \dot{I}_t\rangle_{I_t} - \cH\right) da_t dt, $$
where $\cH$ is the Hamiltonian density and the momentum $\pi_t$ is given by:
\begin{equation}
\pi_t = \II_t - 2H I_t .
\end{equation}
Here $I_t$ is the induced metric, $H$ is the mean curvature, and $\II_t$ is the second fundamental
form of $S_t$. The Hamiltonian density is given by:
\begin{equation}
\cH = \langle\pi_t, \dot{I_t}\rangle_{I_t} - {\cal L},
\end{equation}
where $\cal L$ is the Lagrangian density.

In this formulation, the reduced phase space of the theory is just the space
of classical solutions, in our case the moduli space of constant curvature
3-manifolds. To find the symplectic structure on this phase space it is
sufficient to compute the pre-symplectic one-form: 
$$ \Theta = \int_{S_t} \langle \pi_t, \delta I_t\rangle_{I_t} da_t~. $$
To this end, we shall use the complex-analytic description. Because the
pre-symplectic one-form does not depend on which surface $S_t$ is used for its
calculation, it is convenient to use the maximal surface itself. 

The metric on
the maximal surface is $I_0=e^\varphi |dz|^2$. The second fundamental form is
the real part  $\II_0=(t dz^2+\bar{t} d\bar{z}^2)/2$ of a holomorphic quadratic
differential $t$. It is easy to compute the conjugate momentum. Because $H=0$
for the maximal surface, one only has to find the tensor $\II$, which is an
elementary computation. We express the result as the momentum tensor ``with both of
its indices raised'': 
$$ \pi = e^{-\varphi} \left( 
\bar{t} \left( \frac{\dr}{\dr z} \right)^2 + t \left( \frac{\dr}{\dr \bar{z}}
\right)^2 \right)~. $$ 
To compute the pre-symplectic one-form we also need to compute the variation
of the metric. This can come from two sources: one can either vary the
quadratic differential $t$ or  vary the complex structure on $S$. Under both of
these variations the Liouville field $\varphi$ that satisfies the Gauss
constraint changes in a complicated way. Fortunately, we do not need to know
any explicit expression. Indeed, under a finite deformation the metric
becomes:  $I'_0=e^{\varphi'} |df|^2$, where $f$ is the new complex coordinate
defining the new complex structure. This metric can be pulled back to the
surface $S$ as $h'=e^{\varphi'\circ f} |f_z|^2 |dz+ \mu d\bar{z}|^2$, where
$f$ is now assumed to satisfy the Beltrami equation $f_{\bar{z}} = \mu
f_z$. The above formula can be used to obtain an expression for the first
variation of the metric. One gets: 
$$ \delta I_0 = (\delta\varphi  + \varphi_z \delta f  + \varphi_{\bar{z}}
\delta 
\bar{f} + \delta f_z + \delta \overline{f_z} + \delta\mu + \delta
\bar{\mu})I_0 +  e^\varphi (dz^2 \delta\bar{\mu}+d\bar{z}^2 \delta\mu)~. $$
It is clear that only the last term will give the contribution to the
pre-symplectic one form and 
that we do not need to express all other terms in terms of variations $\delta
\mu, \delta t$ and 
their complex conjugates. It is now trivial to see that:
\begin{equation}
\Theta = \int (t\delta\mu + \bar{t}\delta\bar{\mu}).
\end{equation}
This proves that the symplectic structure coming from gravity is just the
usual symplectic 
structure on $T^* \cT_\Sigma$.

It is of interest to compute the Hamiltonian as well. The Hamiltonian does not
have to 
be time-independent, and in the formulation of \cite{Moncrief} one does get an
explicitly time-dependent answer, so one might suspect that the same happens
in our setting. Thus, we cannot compute it using the maximal surface only. To
compute the Hamiltonian we use: $\dot{I_0} = 2\II_0$. The
Hamiltonian density is equal to:
\begin{equation}
\cH = \tr_{I_0}(\III_0) - (\tr_{I_0}(\II_0)^2 = - 2\, {\rm det}_{I_0}(\II_0)~.
\end{equation}
However, as we know from Lemma \ref{lm:ads}, the shape operator of $S_t$ is
given by: 
$$ B_t = (E\cos(t) + B \sin(t))^{-1} (E\sin(t) + B\cos(t)), $$
where $B=B_0$ is the shape operator of the maximal surface. 
This immediately gives:
\begin{equation}
\cH da_t = -( \sin^2(t) + \det B \cos^2(t))da.
\end{equation}
Thus, the time dependent Hamiltonian is given by:
\begin{equation}
\int_{S} (\cos^2(t) k^2 - \sin^2(t))da,
\end{equation}
where $k$ is the principal curvature of the maximal surface. This can be
further 
simplified in terms of the data on the maximal slice. One gets:
\begin{equation}
A_0 \cos(2t)  - 2\pi(2g-2) \cos^2(t),
\end{equation}
where we have used the Gauss constraint equation (\ref{eq:*-ads}).
Thus, the Hamiltonian is indeed explicitly time dependent, similarly to what
happens in the case of constant mean curvature slicing of \cite{Moncrief}. 

\subsection{Other settings}

Other situations are treated similarly, the only change is in the
trigonometric functions used and some signs. 
In particular, as it is clear from the computation of the pre-symplectic
one-form in the 
previous subsection, the symplectic structure is {\it always} that of the
cotangent 
bundle over the Teichm\"uller space.

\section*{Acknowledgements}

The second author would like to thank Thierry Barbot, Fran\c{c}ois B\'eguin
and  Philippe Eyssidieux for helpful
conversations during the writing of this paper. The first author is supported 
by the EPSRC advanced fellowship.


\begin{thebibliography}{BdCS97}

\bibitem[Ahl66]{ahlfors}
L.~V. Ahlfors.
\newblock {\em Lectures on quasiconformal mappings}.
\newblock D. Van Nostrand Co., Inc., Toronto, Ont.-New York-London, 1966.
\newblock Manuscript prepared with the assistance of Clifford J. Earle, Jr. Van
  Nostrand Mathematical Studies, No. 10.

\bibitem[BB05]{Benedetti}
Riccardo Benedetti and Francesco Bonsante.
\newblock Canonical wick rotations in 3-dimensional gravity.
\newblock math.DG/0508485, 2005.

\bibitem[BBZ03]{BBZ-cras}
Thierry Barbot, Fran{\c{c}}ois B{\'e}guin, and Abdelghani Zeghib.
\newblock Feuilletages des espaces temps globalement hyperboliques par des
  hypersurfaces \`a courbure moyenne constante.
\newblock {\em C. R. Math. Acad. Sci. Paris}, 336(3):245--250, 2003.

\bibitem[BdCS97]{berard-docarmo-santos}
P.~B{\'e}rard, M.~do~Carmo, and W.~Santos.
\newblock The index of constant mean curvature surfaces in hyperbolic
  {$3$}-space.
\newblock {\em Math. Z.}, 224(2):313--326, 1997.

\bibitem[Bes87]{Be}
A.~Besse.
\newblock {\em Einstein Manifolds}.
\newblock Springer, 1987.

\bibitem[Bro04]{bromberg2}
K.~Bromberg.
\newblock Rigidity of geometrically finite hyperbolic cone-manifolds.
\newblock {\em Geom. Dedicata}, 105:143--170, 2004.

\bibitem[BZ04]{barbot-zeghib}
Thierry Barbot and Abdelghani Zeghib.
\newblock Group actions on {L}orentz spaces, mathematical aspects: a survey.
\newblock In {\em The Einstein equations and the large scale behavior of
  gravitational fields}, pages 401--439. Birkh\"auser, Basel, 2004.

\bibitem[Foc04]{Fock}
V.~Fock.
\newblock Talk given at ``quantum hyperbolic geometry'' workshop,
  {AEI-Potsdam}.
\newblock 2004.

\bibitem[Gra90]{gray}
A.~Gray.
\newblock {\em Tubes}.
\newblock Addison-Wesley, 1990.

\bibitem[Gro81]{gromov-bourbaki}
Michael Gromov.
\newblock Hyperbolic manifolds (according to {T}hurston and {J}\o rgensen).
\newblock In {\em Bourbaki Seminar, Vol. 1979/80}, volume 842 of {\em Lecture
  Notes in Math.}, pages 40--53. Springer, Berlin, 1981.

\bibitem[HK98]{HK}
Craig~D. Hodgson and Steven~P. Kerckhoff.
\newblock Rigidity of hyperbolic cone-manifolds and hyperbolic {Dehn} surgery.
\newblock {\em J. Differential Geom.}, 48:1--60, 1998.

\bibitem[Hod05]{hodge-these}
Thomas W.~S. Hodge.
\newblock {\em Hyperk{\"a}hler geometry and Teichm{\"u}ller space}.
\newblock PhD thesis, Imperial College, 2005.

\bibitem[Hop51]{hopf}
Heinz Hopf.
\newblock \"{U}ber {F}l\"achen mit einer {R}elation zwischen den
  {H}auptkr\"ummungen.
\newblock {\em Math. Nachr.}, 4:232--249, 1951.

\bibitem[Kap01]{kapovich}
Michael Kapovich.
\newblock {\em Hyperbolic manifolds and discrete groups}, volume 183 of {\em
  Progress in Mathematics}.
\newblock Birkh\"auser Boston Inc., Boston, MA, 2001.

\bibitem[Kra00]{Holography}
Kirill Krasnov.
\newblock Holography and {R}iemann surfaces.
\newblock {\em Adv. Theor. Math. Phys.}, 4(4):929--979, 2000.

\bibitem[Kra02]{Continuation}
Kirill Krasnov.
\newblock Analytic continuation for asymptotically {A}d{S} 3{D} gravity.
\newblock {\em Classical Quantum Gravity}, 19(9):2399--2424, 2002.

\bibitem[Lab91]{L6}
Fran{\c{c}}ois Labourie.
\newblock Probl\`eme de {M}inkowski et surfaces \`a courbure constante dans les
  vari\'et\'es hyperboliques.
\newblock {\em Bull. Soc. Math. France}, 119(3):307--325, 1991.

\bibitem[Lab92]{L5}
F.~Labourie.
\newblock Surfaces convexes dans l'espace hyperbolique et {CP1}-structures.
\newblock {\em J. London Math. Soc., II. Ser.}, 45:549--565, 1992.

\bibitem[Mes90]{mess}
G.~Mess.
\newblock Lorentz spacetimes of constant curvature.
\newblock Preprint I.H.E.S./M/90/28, 1990.

\bibitem[Mon89]{Moncrief}
Vincent Moncrief.
\newblock Reduction of the {E}instein equations in {$2+1$} dimensions to a
  {H}amiltonian system over {T}eichm\"uller space.
\newblock {\em J. Math. Phys.}, 30(12):2907--2914, 1989.

\bibitem[O'N83]{O}
B.~O'Neill.
\newblock {\em Semi-Riemannian Geometry}.
\newblock Academic Press, 1983.

\bibitem[Ota96]{otal-hyperbolisation}
Jean-Pierre Otal.
\newblock Le th\'eor\`eme d'hyperbolisation pour les vari\'et\'es fibr\'ees de
  dimension 3.
\newblock {\em Ast\'erisque}, 235:x+159, 1996.

\bibitem[RH93]{RH}
Igor Rivin and Craig~D. Hodgson.
\newblock A characterization of compact convex polyhedra in hyperbolic 3-space.
\newblock {\em Invent. Math.}, 111:77--111, 1993.

\bibitem[Riv86]{Ri}
Igor Rivin.
\newblock {\em Thesis}.
\newblock PhD thesis, Princeton University, 1986.

\bibitem[Rub]{rubinstein}
J.H. Rubinstein.
\newblock Minimal surfaces in geometric 3-manifolds.
\newblock to appear in the Global theory of minimal surfaces, Proceedings of
  the Clay Institute workshop at MSRI, edited by D. Hoffman.

\bibitem[Sch98]{shu}
Jean-Marc Schlenker.
\newblock M\'etriques sur les poly\`edres hyperboliques con\-vexes.
\newblock {\em J. Differential Geom.}, 48(2):323--405, 1998.

\bibitem[Spi75]{spivak}
M.~Spivak.
\newblock {\em A comprehensive introduction to geometry, Vol.I-V}.
\newblock Publish or perish, 1970-1975.

\bibitem[Tau04]{taubes}
Clifford~Henry Taubes.
\newblock Minimal surfaces in germs of hyperbolic 3--manifolds.
\newblock {\em GEOM. TOPOL. MONOGR.}, 7:69, 2004.

\bibitem[Thu80]{thurston-notes}
William~P. Thurston.
\newblock Three-dimensional geometry and topology.
\newblock Recent version available on
  http://www.msri.org/publications/books/gt3m/, 1980.

\bibitem[Tro91]{troyanov}
Marc Troyanov.
\newblock Prescribing curvature on compact surfaces with conical singularities.
\newblock {\em Trans. Amer. Math. Soc.}, 324(2):793--821, 1991.

\bibitem[TT03]{takhtajan-teo}
Leon~A. Takhtajan and Lee-Peng Teo.
\newblock Liouville action and {W}eil-{P}etersson metric on deformation spaces,
  global {K}leinian reciprocity and holography.
\newblock {\em Comm. Math. Phys.}, 239(1-2):183--240, 2003.

\bibitem[TZ01]{TZ}
L.~Takhtajan and P.~Zograf.
\newblock Hyperbolic 2-spheres with conical singularities, accessory parameters
  and kaehler metrics on $\mathcal{M}_{0,n}$.
\newblock math.cv/0112170, 2001.

\bibitem[Uhl83]{uhlenbeck}
Karen~K. Uhlenbeck.
\newblock Closed minimal surfaces in hyperbolic {$3$}-manifolds.
\newblock In {\em Seminar on minimal submanifolds}, volume 103 of {\em Ann. of
  Math. Stud.}, pages 147--168. Princeton Univ. Press, Princeton, NJ, 1983.

\bibitem[Wit89]{witten-jones}
Edward Witten.
\newblock  (2+1)-Dimensional Gravity As An Exactly Soluble System.
\newblock {\em  Nucl. Phys. B}, 311:46, 1988.

\end{thebibliography}

\def\cprime{$'$}

\end{document}